%% file: preprint.tex
\newcount\Comments  
\documentclass[letterpaper,11pt]{scrartcl}

\usepackage[utf8]{inputenc}
\usepackage[english]{babel}
\usepackage{amsmath,amstext,amsbsy,amsopn,amscd,amsxtra,upref,amssymb}
\usepackage{amsfonts,bm}
\usepackage{amsthm}
\usepackage{color}
\usepackage{graphicx}
\usepackage{sidecap}
\usepackage{multirow}
\usepackage{booktabs}
\usepackage{bm}
\usepackage{fullpage}

\usepackage{tabularx}
\usepackage[font=small]{caption}
\usepackage{subfig}

\usepackage{longtable}
\usepackage{rotating}

\usepackage{algorithm}
\usepackage{algpseudocode}
\extrafloats{100}

\usepackage{lineno,hyperref}

\graphicspath{{figures/}}
\DeclareGraphicsExtensions{.pdf,.eps,.png,.jpg,.jpeg}

\newcommand{\bfalpha}{\boldsymbol \alpha}
\newcommand{\bfbeta}{\boldsymbol \beta}

\newcommand{\bff}{\boldsymbol f}
\newcommand{\bfg}{\boldsymbol g}

\newcommand{\Dcal}{\mathcal{D}}

\newcommand{\Vcal}{\mathcal{V}}

\newcommand{\bfmu}{\boldsymbol \mu}

\newcommand{\bfe}{\boldsymbol e}

\newcommand{\bfp}{\boldsymbol p}
\newcommand{\bfq}{\boldsymbol q}

\newcommand{\bfs}{\boldsymbol s}

\newcommand{\bfC}{\boldsymbol C}

\newcommand{\bfF}{\boldsymbol F}
\newcommand{\bfG}{\boldsymbol G}

\newcommand{\bfP}{\boldsymbol P}
\newcommand{\bfR}{\boldsymbol R}

\newcommand{\bfV}{\boldsymbol V}

\newcommand{\bfQ}{\boldsymbol Q}

\newcommand{\bfS}{\boldsymbol S}

\newcommand{\nh}{N}
\newcommand{\nr}{n}
\newcommand{\nrp}{m}
\newcommand{\nrs}{m_s}

\newcommand{\win}{w}
\newcommand{\winit}{w_{\text{init}}}
\newcommand{\nz}{z}

\newcommand{\tbff}{\tilde{\bff}}

\newcommand{\Qcal}{\mathcal{Q}}

\newcommand{\tbfq}{\tilde{\bfq}}
\newcommand{\hbfq}{\hat{\bfq}}
\newcommand{\cbfS}{\breve{\bfS}}

\graphicspath{{figures/}}

\newenvironment{keywords}%
   {\begin{trivlist}\item[]{\bfseries\sffamily Keywords:}\ }
   {\end{trivlist}}

\theoremstyle{definition}

\ifpdf
\hypersetup{
  pdftitle={Reduced models with nonlinear approximations of latent dynamics for model premixed flame problems}, 
  pdfauthor={Wayne Isaac Tan Uy and Christopher R. Wentland and Cheng Huang and Benjamin Peherstorfer} 
}
\fi

\makeatletter 
\renewcommand\Huge{\@setfontsize\Huge{20pt}{18}}
\renewcommand\huge{\@setfontsize\huge{15pt}{18}}
\makeatother

\begin{document}

\title{Reduced models with nonlinear approximations of latent dynamics for model premixed flame problems\thanks{This work was supported in part by the Air Force Center of Excellence on Multi-Fidelity Modeling of Rocket Combustor Dynamics, award FA9550-17-1-0195. The fourth author was additionally partially supported by NSF-2012250 and by the Air Force Office of Scientific Research (AFOSR) award FA9550-21-1-0222 (Dr.~Fariba Fahroo).}}

\author{Wayne Isaac Tan Uy\footnote{Courant Institute of Mathematical Sciences, New York University, New York NY} \and Christopher R. Wentland\footnote{University of Michigan, Ann Arbor MI} \and Cheng Huang\footnote{University of Kansas, Lawrence KS}  \and Benjamin Peherstorfer${}^1$}

\maketitle              

\begin{abstract}
Efficiently reducing models of chemically reacting flows  
is often challenging because their characteristic features such as sharp gradients in the flow fields and couplings over various time and length scales lead to dynamics that evolve in high-dimensional spaces.
In this work, we show that online adaptive reduced models that construct nonlinear approximations by adapting low-dimensional subspaces over time can predict well latent dynamics with properties similar to those found in chemically reacting flows. The adaptation of the subspaces is driven by the online adaptive empirical interpolation method, which takes sparse residual evaluations of the full model to compute low-rank basis updates of the subspaces. 
Numerical experiments with a premixed flame model problem show that reduced models based on online adaptive empirical interpolation accurately predict flame dynamics far outside of the training regime and in regimes where traditional static reduced models, which keep reduced spaces fixed over time and so provide only linear approximations of latent dynamics, fail to make meaningful predictions. 
\end{abstract}

\input{maintext}

\end{document}

%% file: maintext.tex
\section{Introduction}
Even with advances in modern computational capabilities, high-fidelity, full-scale simulations of chemically reacting flows in realistic applications remain computationally expensive \cite{HUANGMPLSVT,doi:10.1080/10618562.2014.911848,doi:10.2514/6.2018-4675,doi:10.2514/6.2018-1183}.
Traditional model reduction methods \cite{RozzaPateraSurvey,SIREVSurvey,SerkanInterpolatory,PWG17MultiSurvey} that seek reduced solutions in low-dimensional subspaces fail for problems that involve chemically reacting advection-dominated flows because the strong advection of sharp gradients in the solution fields lead to high-dimensional features in the latent dynamics; see \cite{Notices} for an overview of the challenges of model reduction of strongly advecting flows and other problems with high-dimensional latent dynamics.
We demonstrate on a model premixed flame problem \cite{chris_wentland_2021_5517532}---which greatly simplifies the reaction and flow dynamics but preserves some of the model reduction challenges of more realistic chemically reacting flows---that adapting subspaces of reduced models over time \cite{Peherstorfer15aDEIM,P18AADEIM,CKMP19ADEIMQuasiOptimalPoints} can help to provide accurate future-state predictions with only a few degrees of freedom.

Traditional reduced models are formulated via projection-based approaches that seek approximate solutions in lower dimensional subspaces of the high-dimensional solution spaces of full models; see \cite{RozzaPateraSurvey,SIREVSurvey,SerkanInterpolatory} for surveys on model reduction. Mathematically, the traditional approximations in subspaces lead to linear approximations in the sense that the parameters of the reduced models that can be changed over time enter linearly in the reduced solutions. 
It has been observed empirically that for certain types of dynamics, which are found in a wide range of science and engineering applications, including chemically reacting flows, the accuracy of such linear approximations grows slowly with the dimension $\nr$ of the reduced space.
Examples of such dynamics are flows that are dominated by advection. In fact, for solutions of the linear advection equation, it has been shown that under certain assumptions on the metric and the ambient space the best-approximation error of linear approximations in subspaces cannot decay faster than  $1/\sqrt{\nr}$.
This slowly decaying lower bound is referred to as the Kolmogorov barrier, because the Kolmogorov $\nr$-width of a set of functions is defined as the best-approximation error obtained over all subspaces; see \cite{Ohlberger16,Greif19,10.1093/imanum/dru066} for details.

A wide range of methods has been introduced that aim to circumvent the Kolmogorov barrier; see \cite{Notices} for a brief survey.
There are methods that introduce nonlinear transformations and nonlinear embeddings to recover low-dimensional structures. Examples are transformations based on Wasserstein metrics \cite{ehrlacher19}, deep network and deep autoencoders \cite{LEE2020108973,KIM2022110841,https://doi.org/10.48550/arxiv.2203.00360,https://doi.org/10.48550/arxiv.2203.01360}, shifted proper orthogonal decomposition and its extensions \cite{doi:10.1137/17M1140571,PAPAPICCO2022114687}, quadratic manifolds \cite{https://doi.org/10.48550/arxiv.2205.02304,BARNETT2022111348}, and other transformations \cite{OHLBERGER2013901,TaddeiShock,https://doi.org/10.48550/arxiv.1911.06598,Cagniart2019}.
In this work, we focus on online adaptive reduced models that adapt the reduced space over time to achieve nonlinear approximations. In particular, we build on online adaptive empirical interpolation with adaptive sampling (AADEIM), which adapts reduced spaces with additive low-rank updates that are derived from sparse samples of the full-model residual \cite{Peherstorfer15aDEIM,P18AADEIM,CKMP19ADEIMQuasiOptimalPoints}.
The AADEIM method builds on empirical interpolation \cite{barrault_empirical_2004,deim2010,QDEIM}. We refer to \cite{SAPSIS20092347,doi:10.1137/050639703,ZPW17SIMAXManifold,doi:10.1137/16M1088958,doi:10.1137/140967787,Musharbash2020,refId0Hesthaven} for other adaptive basis and adaptive low-rank approximations. The idea of evolving basis functions over time has a long history in numerical analysis and scientific computing, which dates back to at least Dirac~\cite{dirac_1930}.

We apply AADEIM to construct reduced models of a model premixed flame problem with artificial pressure forcing. Our numerical results demonstrate that reduced models obtained with AADEIM provide accurate predictions of the fluid flow and flame dynamics with only a few degrees of freedom. In particular, the AADEIM model that we derive predicts the dynamics far outside of the training regime and in regimes where traditional, static reduced models, which keep the reduced spaces fixed over time, fail to provide meaningful prediction. 

The manuscript is organized as follows. We first provide preliminaries in Section~\ref{sec:Prelim} on traditional, static reduced modeling. We then recap AADEIM in Section~\ref{sec:AADEIM} and highlight a few modifications that we made compared to the original AADEIM method introduced in \cite{Peherstorfer15aDEIM,P18AADEIM}. The reacting flow solver PERFORM \cite{chris_wentland_2021_5517532} and the model premixed flame problem is discussed in Section~\ref{sec:PERFORM}. Numerical results that demonstrate AADEIM on the premixed flame problem are shown in Section~\ref{sec:NumExp} and conclusions are drawn in Section~\ref{sec:Conc}.

\section{Static model reduction with empirical interpolation}\label{sec:Prelim}
We briefly recap model reduction with empirical interpolation \cite{barrault_empirical_2004,grepl_efficient_2007,deim2010,QDEIM} using reduced spaces that are fixed over time. We refer to reduced models with fixed reduced spaces as static models in the following sections.

\subsection{Static reduced models}
Discretizing a system of partial differential equations in space and time can lead to a dynamical system of the form
\begin{equation}\label{eq:FOM}
\bfq_k(\bfmu) = \bff(\bfq_{k-1}(\bfmu); \bfmu)\,,\qquad k = 1, \dots, K\,,
\end{equation}
with state $\bfq_k(\bfmu) \in \mathbb{R}^{\nh}$ at time step $k = 1, \dots, K$ and physical parameter $\bfmu \in \Dcal$. The function $\bff: \mathbb{R}^{\nh} \times \Dcal \to \mathbb{R}^{\nh}$ is vector-valued and nonlinear in the first argument $\bfq$ in the following.
A system of form \eqref{eq:FOM} is obtained, for example, after an implicit time discretization of the time-continuous system.
The initial condition is $\bfq_0(\bfmu) \in \Qcal_0 \subseteq \mathbb{R}^{\nh}$ and is an element of the set of initial conditions~$\Qcal_0$.

Consider now training parameters $\bfmu_1, \dots, \bfmu_M \in \Dcal$ with training initial conditions $\bfq_0(\bfmu_1), \dots, \bfq_0(\bfmu_M)$. Let further $\bfQ(\bfmu_1), \dots, \bfQ(\bfmu_M) \in \mathbb{R}^{\nh \times (K + 1)}$ be the corresponding training trajectories defined as
\[
\bfQ(\bfmu_i) = [\bfq_0(\bfmu_i), \dots, \bfq_K(\bfmu_i)]\,,\qquad i = 1, \dots, M\,.
\]
From the snapshot matrix $\bfQ = [\bfQ(\bfmu_1), \dots, \bfQ(\bfmu_M)] \in \mathbb{R}^{\nh \times M(K+1)}$, a reduced space $\Vcal$ of dimension $\nr \ll \nh$ with basis matrix $\bfV \in \mathbb{R}^{\nh \times \nr}$ is constructed, for example, via proper orthogonal decomposition (POD) or greedy methods \cite{RozzaPateraSurvey,SIREVSurvey}.

The static Galerkin reduced model is
\[
\tilde{\bfq}_k(\bfmu) = \bfV^T\bff(\bfV\tilde{\bfq}_{k - 1}(\bfmu); \bfmu)\,,\qquad k = 1, \dots, K\,,
\]
with the initial condition $\tbfq_0(\bfmu) \in \tilde{\Qcal}_0 \subseteq \mathbb{R}^{\nr}$ for a parameter $\bfmu \in \Dcal$ and reduced state $\tbfq_k(\bfmu) \in \mathbb{R}^{\nr}$ at time $k = 1, \dots, K$.
However, evaluating the function $\bff$ still requires evaluating $\bff$ at all $\nh$ components. Empirical interpolation \cite{barrault_empirical_2004,grepl_efficient_2007,deim2010,QDEIM} provides an approximation $\tbff: \mathbb{R}^{\nr} \times \Dcal \to \mathbb{R}^{\nr}$ of $\bff$ that can be evaluated $\tbff(\tbfq; \bfmu)$ at a vector $\tbfq \in \mathbb{R}^{\nr}$ with costs that grow as $\mathcal{O}(\nr)$.
Consider the matrix $\bfP = [\bfe_{p_1}, \dots, \bfe_{p_\nrp}]\in \mathbb{R}^{\nh \times \nrp}$ that has as columns the $\nh$-dimensional unit vectors with ones at $\nrp$ unique components $p_1, \dots, p_{\nrp} \in \{1, \dots, \nh\}$. It holds for the number of points $\nrp \geq \nr$. The points $p_1, \dots, p_{\nrp}$ can be computed, for example, with greedy \cite{barrault_empirical_2004,deim2010}, QDEIM \cite{QDEIM}, or oversampling algorithms \cite{PDG18ODEIM,doi:10.2514/6.2021-1371}.
We denote with $\bfP^T\bff(\bfq; \bfmu)$ that only the component functions $f_{p_1}, \dots, f_{p_{\nrp}}$ corresponding to the points $p_1, \dots, p_{\nrp}$ of the vector-valued function $\bff = [f_1, \dots, f_{\nh}]$ are evaluated at $\bfq \in \mathbb{R}^{\nh}$ and $\bfmu \in \Dcal$. The empirical-interpolation approximation of $\bff$ is 
\[
\tbff(\tbfq; \bfmu) = (\bfP^T\bfV)^{\dagger}\bfP^T\bff(\bfV\tilde{\bfq}; \bfmu)
\]
where $(\bfP^T\bfV)^{\dagger}$ denotes the Moore--Penrose inverse (pseudoinverse) of $\bfP^T\bfV$.

Based on the empirical-interpolation approximation $\tbff$, we derive the static reduced model
\[
\tilde{\bfq}_k(\bfmu) = \tilde{\bff}(\tilde{\bfq}_{k - 1}(\bfmu); \bfmu)\,,\qquad k = 1, \dots, K\,,
\]
with the reduced states $\tbfq_k(\bfmu)$ at time steps $k = 1, \dots, K$.

\subsection{The Kolmogorov barrier of static reduced models}\label{sec:ProblemFormulation}
The empirical-interpolation approximation $\tilde{\bff}$ depends on the basis matrix $\bfV$, which is fixed over all time steps $k = 1, \dots, K$. This means that the reduced approximation $\tilde{\bfq}_k$ at time $k$ depends linearly on the basis vectors of the reduced space $\Vcal$, which are the columns  of the basis matrix $\bfV$. 
Thus, the lowest error that such a static reduced model can achieve is related to the Kolmogorov $\nr$-width, i.e., the best-approximation error in any subspace of dimension $\nr$. We refer to \cite{Notices} for an overview of the Kolmogorov barrier in model reduction and to \cite{10.1093/imanum/dru066,MADAY2002289,Ohlberger16} for in-depth details. It has been observed empirically, and in some limited cases proven, that systems that are governed by dynamics with strong advection and transport exhibit a slowly decaying $\nr$-width, which means that linear, static reduced models are inefficient in providing accurate predictions.

\section{Online adaptive empirical interpolation methods for nonlinear model reduction}\label{sec:AADEIM}
In this work, we apply online adaptive model reduction methods to problems motivated by chemically reacting flows, which are often dominated by advection and complex transport dynamics that make traditional static reduced models inefficient.
We focus on reduced models obtained with AADEIM, which builds on the online adaptive empirical interpolation method \cite{Peherstorfer15aDEIM} for adapting the basis and on the adaptive sampling scheme introduced in \cite{P18AADEIM,CKMP19ADEIMQuasiOptimalPoints}. We utilize the one-dimensional compressible reacting flow solver PERFORM \cite{chris_wentland_2021_5517532}, which provides several benchmark problems motivated by combustion applications. We will consider a model premixed flame problem with artificial pressure forcing and show that AADEIM provides reduced models that can accurately predict the flame dynamics over time, whereas traditional static reduced models fail to make meaningful predictions.

For ease of exposition, we drop the dependence of the states $\bfq_k(\bfmu), k = 1, \dots, K$ on the parameter $\bfmu$ in this section.

\subsection{Adapting the basis}\label{sec:ADEIMAdaptBasis}
To allow adapting the reduced space over time, we formally make the basis matrix $\bfV_k$ and the points matrix $\bfP_k$ with points $p_1^{(k)}, \dots, p_{\nrp}^{(k)} \in \{1, \dots, \nh\}$  depend on the time step $k = 1, \dots, K$.
In AADEIM, the basis matrix $\bfV_k$ is adapted at time step $k$ to the basis matrix $\bfV_{k + 1}$ via a rank-one update
\[
\bfV_{k + 1} = \bfV_k + \bfalpha_k\bfbeta_k^T\,,
\]
given by $\bfalpha_k \in \mathbb{R}^{\nh}$ and $\bfbeta_k \in \mathbb{R}^{\nr}$.
To compute an update $\bfalpha_k\bfbeta_k^T$ at time step $k$, we introduce the data matrix $\bfF_k \in \mathbb{R}^{\nh \times \win}$, where $\win \in \mathbb{N}$ is a window size. First, similarly the empirical-interpolation points, we consider $\nrs \leq \nh$ sampling points $s_1^{(k)}, \dots, s_{\nrs}^{(k)} \in \{1, \dots, \nh\}$ and the corresponding sampling matrix $\bfS_k = [\bfe_{s_1^{(k)}}, \dots, \bfe_{s_{\nrs}^{(k)}}] \in \mathbb{R}^{\nh \times \nrs}$. We additionally consider the complement set of sampling points $\{1, \dots, \nh\} \setminus \{s_1^{(k)}, \dots, s_{\nrs}^{(k)}\}$ and the corresponding matrix $\cbfS_k$. We will also need the matrix corresponding to the union of the set of sampling points $\{s_1^{(k)}, \dots, s_{\nrs}^{(k)}\}$ and the points $\{p_1, \dots, p_{\nrp}\}$, which we denote with $\bfG$, and its complement $\breve{\bfG}$.

The data matrix $\bfF_k$ at time step $k$ is then given by 
\[
\bfF_k = [\hbfq_{k - w + 1}, \dots, \hbfq_{k}] \in \mathbb{R}^{\nh \times w}\,,
\]
where we add the vector $\hbfq_k$ that is defined as 
\begin{equation}
\bfG^T_k\hbfq_k = \bfG^T_k\bff(\bfV_k\tbfq_{k-1})\,,\quad \breve{\bfG}^T_k\hbfq_k = \breve{\bfG}^T_k\bfV_k(\bfG_k^T\bfV_k)^{\dagger}\bfG_k^T\bff(\bfV_k\tbfq_{k-1})\,.
\label{eq:FillRHSMatrixVector}
\end{equation}
The state $\tbfq_k$ used in \eqref{eq:FillRHSMatrixVector} is the reduced state at time step $k$.
The vector $\hbfq_k$ at time step $k$ serves as an approximation of the full-model state $\bfq_k$ at time $k$. This is motivated by the full-model equations \eqref{eq:FOM} with the reduced state $\tbfq_{k-1}$ as an approximation of the full-model state $\bfq_{k-1}$ at time step $k - 1$; we refer to \cite{P18AADEIM} for details about this motivation.

The AADEIM basis update $\bfalpha_k\bfbeta_k^T$ at time step $k$ is the solution to the minimization problem
\begin{equation}
\min_{\bfalpha_k \in \mathbb{R}^{\nh},\, \bfbeta_k \in \mathbb{R}^{\nr}}\, \left\|(\bfV_k + \bfalpha_k\bfbeta_k^T)\bfC_k - \bfF_k\right\|_F^2\,,\label{eq:ADEIMUpdate}
\end{equation}
where the coefficient matrix is
\begin{equation}
\bfC_k = \bfV_k^T\bfF_k. 
\label{eq:CoeffMat}
\end{equation}
The matrix $\bfP_k$ is adapted to $\bfP_{k + 1}$ by applying QDEIM \cite{QDEIM} to the adapted basis matrix $\bfV_{k + 1}$

We make two modifications compared to the original AADEIM approach. First, we sample from the points given by $\bfG$, which is the union of the sampling points $\{s_1^{(k)}, \dots, s_{\nrs}^{(k)}\}$ and the points $\{p_1, \dots, p_{\nrp}\}$. This comes with no extra costs because the full-model right-hand side function needs to be evaluated at the points corresponding to $\bfS$ and $\bfP$ even in the original AADEIM approach. 
Second, as proposed in \cite{ChengADEIMImprovements}, we adapt the basis at all components of the residual in the objective, rather than only at the sampling points given by $\bfS$. This requires no additional full-model right-hand side function evaluations but comes with increased computational costs when solving the optimization problem \eqref{eq:ADEIMUpdate}. However, typically solving the optimization problem is negligible compared to sparsely evaluating the full-model right-hand side function.

\subsection{Adapting sampling points}\label{sec:ADEIM:AdaptSamplingPoints}
When adapting the sampling matrix $\bfS_{k-1}$ to $\bfS_k$ at time step $k$, we evaluate the full-model right-hand side function $\bff$ at all $\nh$ components to obtain
\begin{equation}\label{eq:SamplingUpdateQk}
\hbfq_k = \bff(\bfV_k\tbfq_{k-1})
\end{equation}
and put it as column into the data matrix $\bfF_k$. We then compute the residual matrix
\[
\bfR_k = \bfF_k - \bfV_k(\bfP_k^T\bfV_k)^{\dagger}\bfP_k^T\bfF_k\,.
\]
Let $r_k^{(i)}$ denote the 2-norm of the $i$-th row of $\bfR_k$ and let $i_1, \dots, i_{\nh}$ be an ordering such that
\[
r^{(i_1)}_k \geq \dots \geq r^{(i_{\nh})}_k\,.
\]
At time step $k$, we pick the first $\nrs$ indices $i_1 = s_1^{(k)}, \dots, i_{\nrs} = s_{\nrs}^{(k)}$ as the sampling points to form $\bfS_k$, which is subsequently used to adapt the basis matrix from $\bfV_k$ to $\bfV_{k + 1}$.

Two remarks are in order. First, the sampling points are quasi-optimal with respect to an upper bound of the adaptation error \cite{CKMP19ADEIMQuasiOptimalPoints}. Second, 
adapting the sampling points requires evaluating the residual at all $i = 1, \dots, \nh$ components, which incurs computational costs that scale with the dimension $\nh$ of the full-model states. However, we adapt the sampling points not every time step, but only every $\nz$-th time step as proposed in \cite{P18AADEIM}.

\begin{algorithm}[t]
\caption{AADEIM algorithm}\label{alg:ABAS}
\begin{algorithmic}[1]
\Procedure{AADEIM}{$\bfq_0, \bff, \bfmu, \nr, \winit, \win, \nrs, \nz$}
\State Solve full model for $\winit$ time steps $\bfQ = \texttt{solveFOM}(\bfq_0, \bff, \bfmu)$\label{alg:ABAS:SolveFOM}
\State Set $k = \winit+1$\label{alg:ABAS:StartROMInit}
\State Compute $\nr$-dimensional POD basis $\bfV_k$ of $\bfQ$\label{alg:ABAS:PODBasisConstruction}
\State Compute QDEIM interpolation points $\bfp_k = \texttt{qdeim}(\bfV_k)$
\State Initialize $\bfF = \bfQ[:, k-\win+1:k-1]$ and $\tilde{\bfq}_{k - 1} = \bfV_k^T\bfQ[:, k - 1]$\label{alg:ABAS:EndROMInit}
\For{$k = \winit + 1, \dots, K$}\label{alg:ABAS:ROMLoop}
\State Solve $\tilde{\bfq}_{k - 1} = \tilde{\bff}(\tilde{\bfq}_k; \bfmu)$ with DEIM, using basis matrix $\bfV_k$ and points $\bfp_k$\label{alg:ABAS:SolveROM}
\State Store $\bfQ[:, k] = \bfV_k\tilde{\bfq}_k$
\If{$\operatorname{mod}(k, \nz) == 0 || k == \winit + 1$}\label{alg:ABAS:IfSampling}
\State Compute $\bfF[:, k] = \bff(\bfQ[:, k]; \bfmu)$\label{alg:ABAS:AdaptSamplingPointsStart}
\State $\bfR_k = \bfF[:, k - \win + 1:k] - \bfV_k(\bfV_k[\bfp_k, :])^{-1}\bfF[\bfp_k, k - \win + 1:k]$
\State $[\sim, \bfs_k] = \texttt{sort}(\texttt{sum}(\bfR_k.\widehat{~~}2, 2), \text{'descend'})$
\State Set $\breve{\bfs}_k = \bfs_k[\nrs+1:\text{end}]$ and $\bfs_k = \bfs_k[1:\nrs]$\label{alg:ABAS:AdaptSamplingPointsEnd}
\Else
\State Set $\bfs_k = \bfs_{k-1}$ and  $\breve{\bfs}_k = \breve{\bfs}_{k - 1}$ 
\State Take the union of points in $\bfs_k$ and $\bfp_k$ to get $\bfg_k$ and complement $\breve{\bfg}_k$
\State Compute $\bfF[\bfg_k, k] = \bff(\bfQ[\bfg_k, k]; \bfmu)$\label{alg:ABAS:EvalResidualAtSamplingPoints}
\State Approximate $\bfF[\breve{\bfg}_k, k] = \bfV_k[\breve{\bfg}_k, :](\bfV_k[\bfg_k, :])^{-1}\bfF[\bfg_k, k]$\label{alg:ABAS:ApproxResidualAtSamplingPoints} 
\EndIf
\State Compute update $\bfalpha_k, \bfbeta_k$ by solving \eqref{eq:ADEIMUpdate} for $\bfF[:, k - \win + 1:k]$ and $\bfV_k$\label{alg:ABAS:CompBasisUpdate}
\State Adapt basis $\bfV_{k + 1} = \bfV_k + \bfalpha_k\bfbeta_k$ and orthogonalize $\bfV_{k + 1}$\label{alg:ABAS:ApplyBasisUpdate}
\State Compute points $\bfp_{k + 1}$ by applying QDEIM to $\bfV_{k + 1}$ \label{alg:ABAS:AdaptP}
\EndFor\\
\Return Return trajectory $\bfQ$
\EndProcedure
\end{algorithmic}
\end{algorithm}

\subsection{Computational procedure and costs}
By combining the basis adaptation described in Section~\ref{sec:ADEIMAdaptBasis} and sampling-points adaptation of Section~\ref{sec:ADEIM:AdaptSamplingPoints}, we obtain the AADEIM reduced model 
\[
\tbfq_k = \tbff_k(\tbfq_{k-1}; \bfmu)\,,\qquad k = 1, \dots, K\,,
\]
where now the approximation $\tbff_k: \mathbb{R}^{\nh} \times \Dcal \to \mathbb{R}^{\nh}$ is
\[
\tbff_k(\tbfq_{k-1}) = (\bfP^T_k\bfV_k)^{\dagger}\bfP_k^T\bff(\bfV_k\tbfq_{k-1})\,,
\]
which depends on the time step $k = 1, \dots, K$ because the basis $\bfV_k$ and the points matrix $\bfP_k$ depend on the time step.
The AADEIM algorithm is summarized in Algorithm~\ref{alg:ABAS}. Inputs to the algorithm are the initial condition $\bfq_0 \in \Qcal_0$, full-model right-hand side function $\bff$, parameter $\bfmu$, reduced dimension $\nr$, initial window size $\winit$, window size $w$, and the frequency of updating the sampling points $\nz$. The algorithm returns the trajectory $\bfQ \in \mathbb{R}^{\nh \times K}$.

Lines~\ref{alg:ABAS:SolveFOM}--\ref{alg:ABAS:EndROMInit} initialize the AADEIM reduced model by first solving the full model for $\winit$ time steps to compute the snapshots and store them in the columns of the matrix $\bfQ$. From these $\winit$ snapshots, a POD basis matrix $\bfV_k$ for time step $k = \winit + 1$ is constructed. The time-integration loop starts in line~\ref{alg:ABAS:ROMLoop}. In each iteration $k = \winit + 1, \dots, K$, the reduced state $\tbfq_{k - 1}$ is propagated forward to obtain $\tbfq_k$ by solving the AADEIM reduced model for one time step. The first branch of the \texttt{if} clause in line~\ref{alg:ABAS:IfSampling} is entered if the sampling points are to be updated, which is the case every $\nz$-th time step. If the sampling points are updated, the full-model right-hand side function is evaluated at all $\nh$ components to compute the residual matrix $\bfR_k$. The new sampling points are selected based on the largest norm of the rows of the entry-wise squared residual matrix $\bfR_k$. If the sampling points are not updated, the full-model right-hand side function $\bff$ is evaluated only at the points corresponding to $\bfs_k$ and $\bfp_k$. All other components are approximated with empirical interpolation.
In lines~\ref{alg:ABAS:CompBasisUpdate} and~\ref{alg:ABAS:ApplyBasisUpdate}, the basis update $\bfalpha_k,\bfbeta_k$ is computed and then used to obtain the adapted basis matrix $\bfV_{k + 1} = \bfV_k + \bfalpha_k\bfbeta_k^T$. The points $\bfp_k$ are adapted to $\bfp_{k + 1}$ by applying QDEIM to the adapted basis matrix $\bfV_{k + 1}$ in line~\ref{alg:ABAS:AdaptP}. The method \texttt{solveFOM()} refers to the full-model solver and the method \texttt{qdeim()} to QDEIM \cite{QDEIM}

\section{Benchmarks of chemically reacting flow problems}\label{sec:PERFORM}
A collection of benchmarks for model reduction of transport-dominated problems is provided with PERFORM \cite{chris_wentland_2021_5517532}. Documentation of the code and benchmark problems is available online\footnote{\url{https://perform.readthedocs.io/}}. The benchmarks are motivated by combustion processes and modeled after the General Equations and Mesh Solver (GEMS), which provides a reacting flow solver in three spatial dimensions \cite{10.1007/3-540-31801-1_89}.

\begin{figure}[p]
\begin{tabular}{ccc}
\resizebox{0.33\columnwidth}{!}{\Huge\input{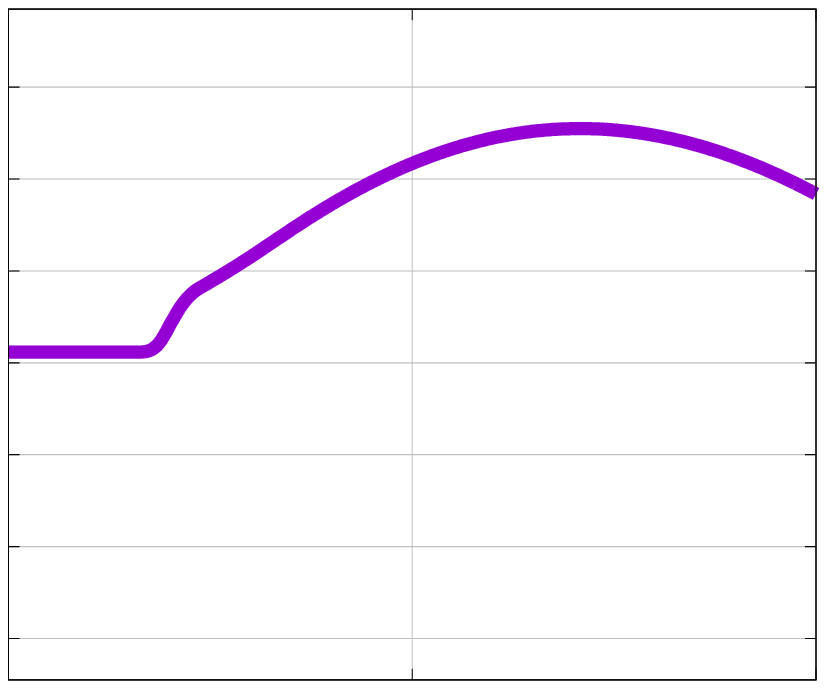}} & \resizebox{0.33\columnwidth}{!}{\Huge\input{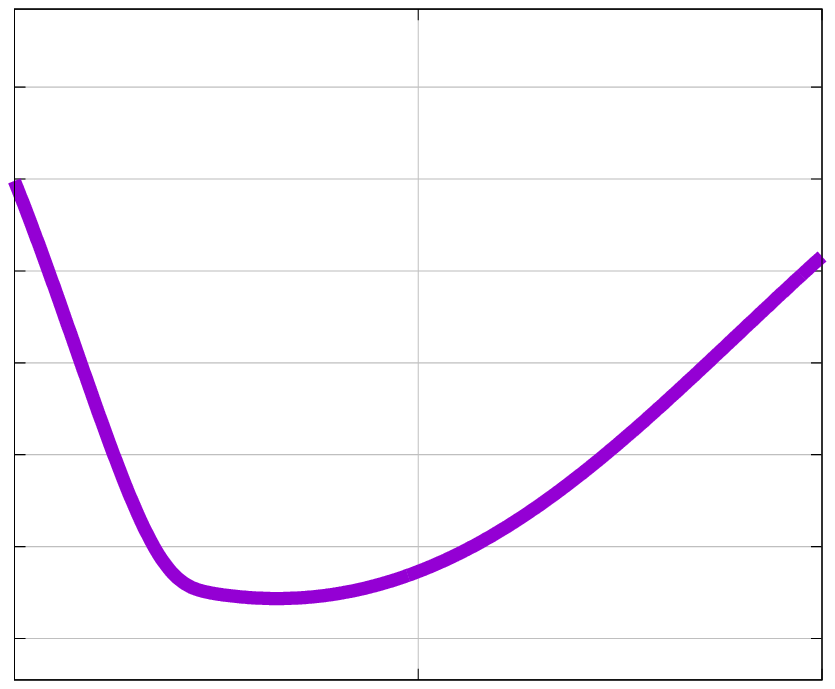}} & \resizebox{0.33\columnwidth}{!}{\Huge\input{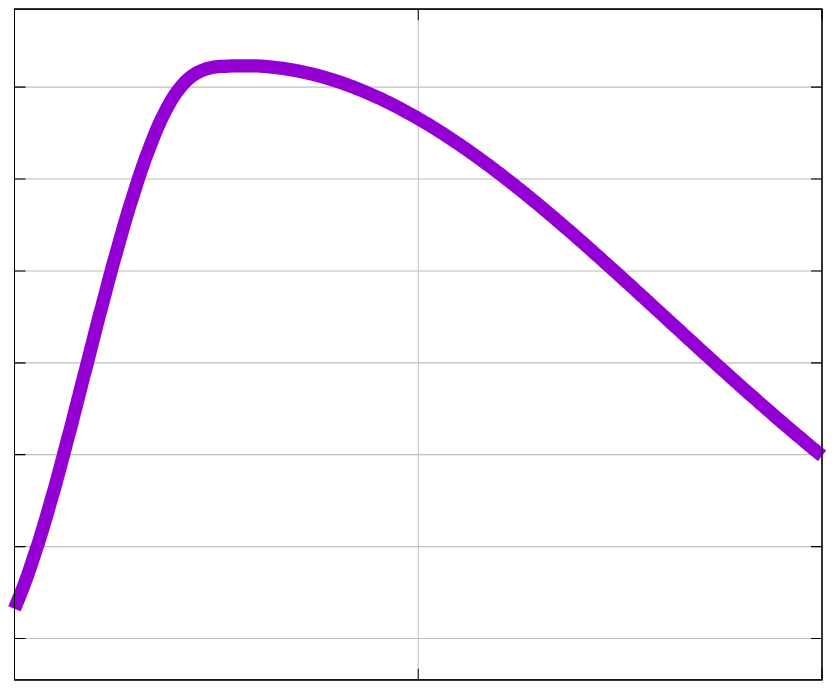}}\\
(a) pressure, $t = 7.5\times 10^{-6}$ & (b) pressure, $t = 2\times 10^{-5}$ & (c) pressure, $t = 3\times 10^{-5}$\\
\resizebox{0.33\columnwidth}{!}{\Huge\input{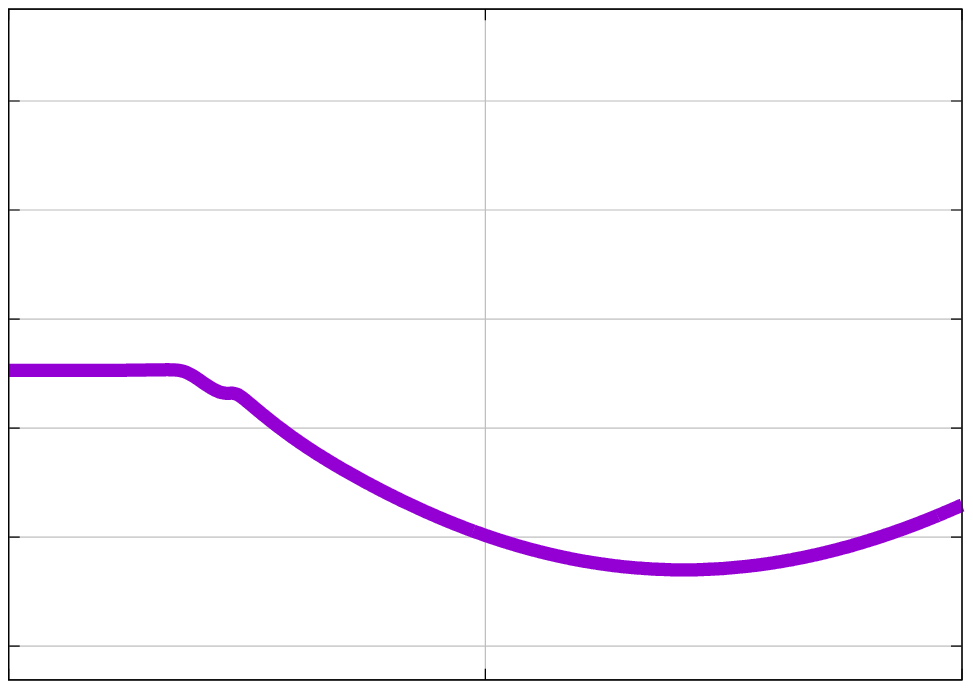}} & \resizebox{0.33\columnwidth}{!}{\Huge\input{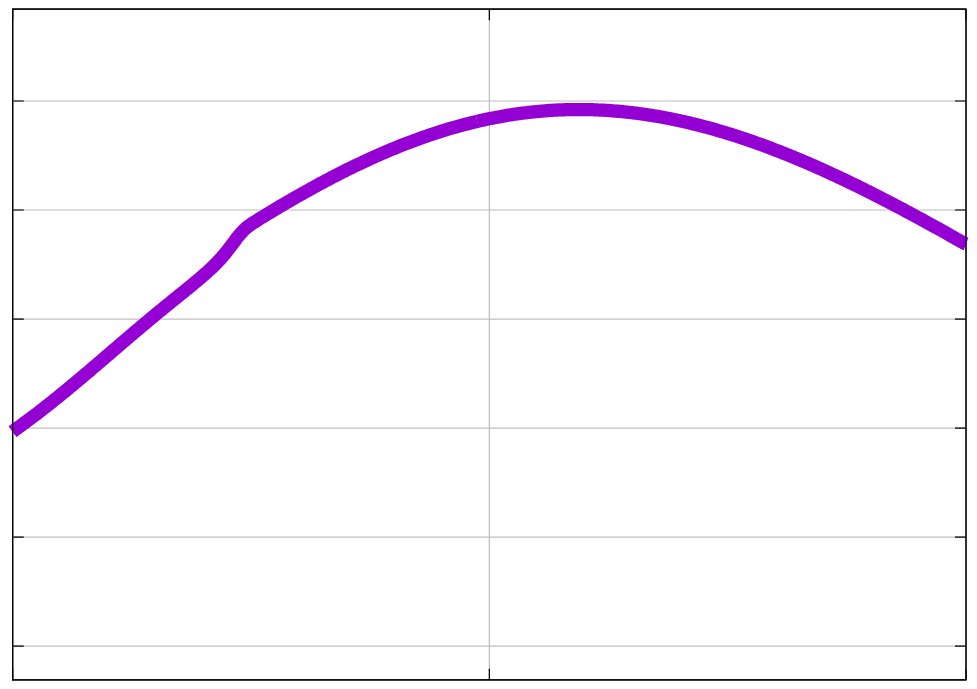}} & \resizebox{0.33\columnwidth}{!}{\Huge\input{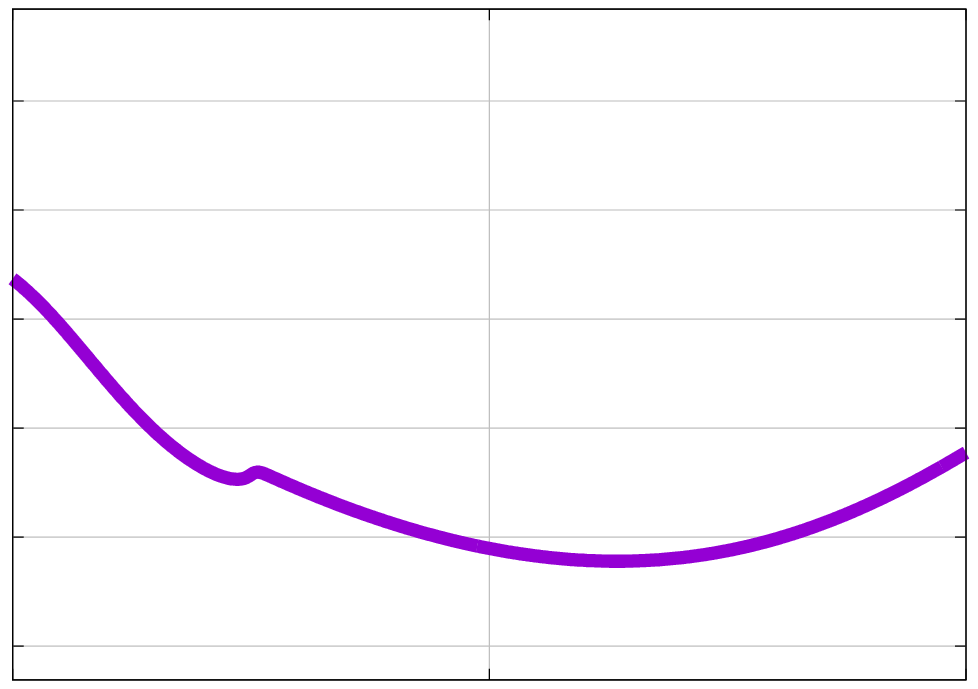}}\\
(d) velocity, $t = 7.5\times 10^{-6}$ & (e) velocity, $t = 2\times 10^{-5}$ & (f) velocity, $t = 3\times 10^{-5}$\\
\resizebox{0.33\columnwidth}{!}{\Huge\input{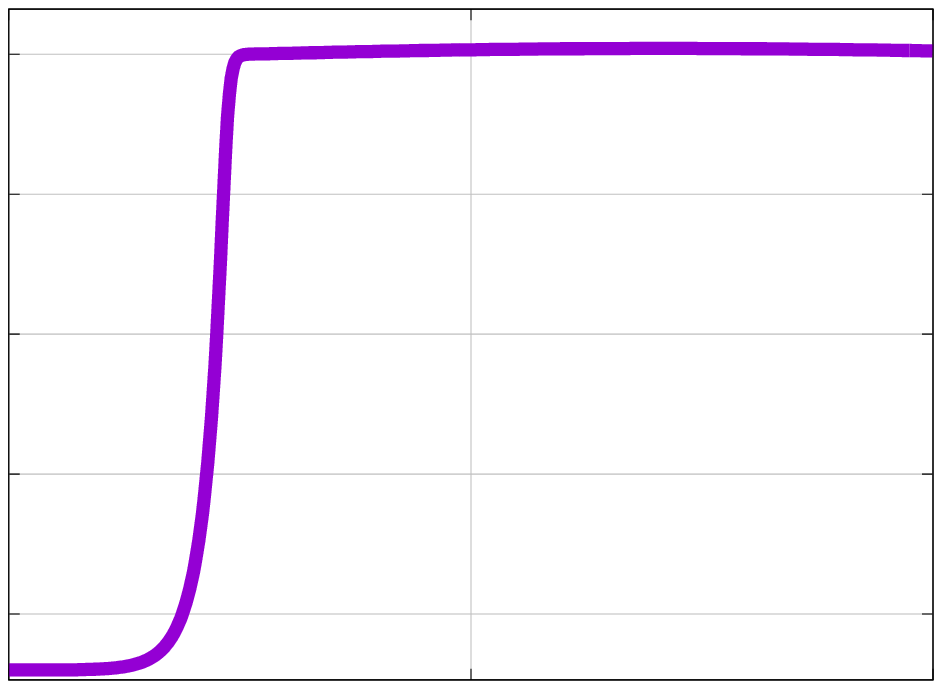}} & \resizebox{0.33\columnwidth}{!}{\Huge\input{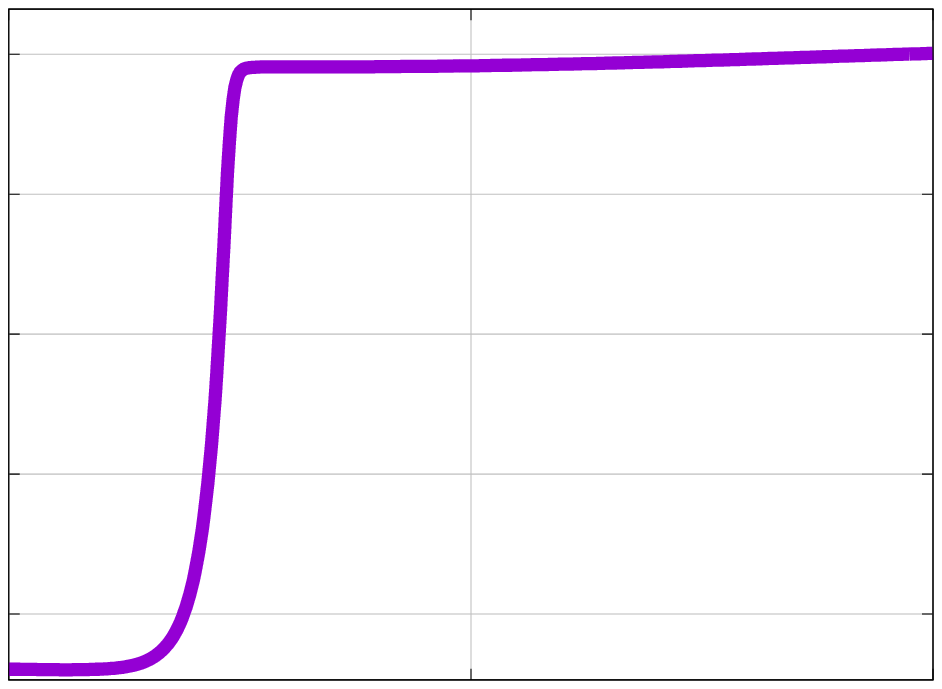}} & \resizebox{0.33\columnwidth}{!}{\Huge\input{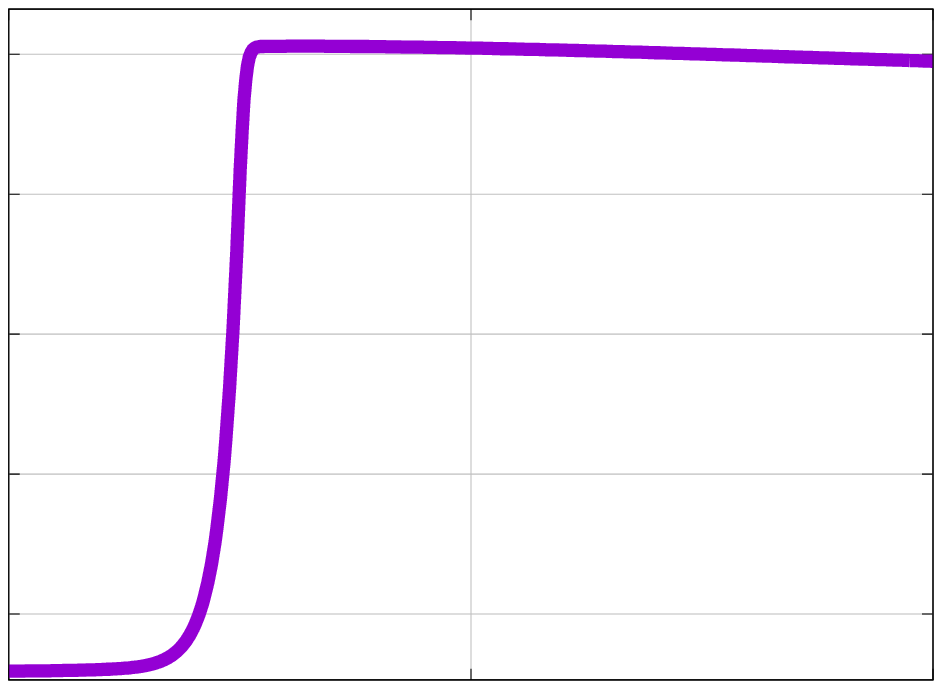}}\\
(g) temperature, $t = 7.5\times 10^{-6}$ & (h) temperature, $t = 2\times 10^{-5}$ & (i) temperature, $t = 3\times 10^{-5}$\\
\resizebox{0.33\columnwidth}{!}{\Huge\input{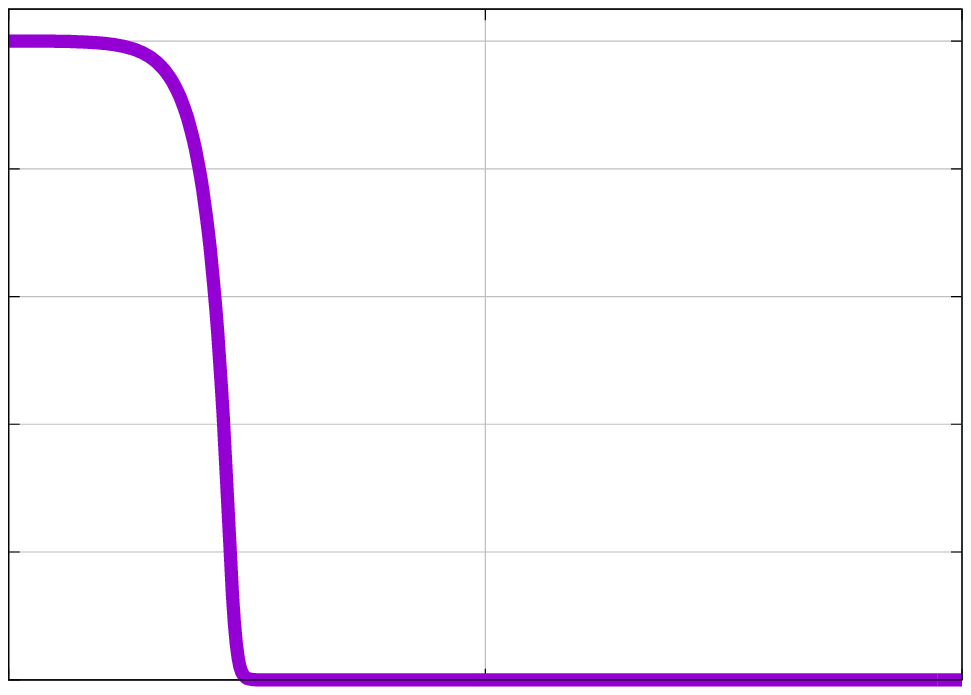}} & \resizebox{0.33\columnwidth}{!}{\Huge\input{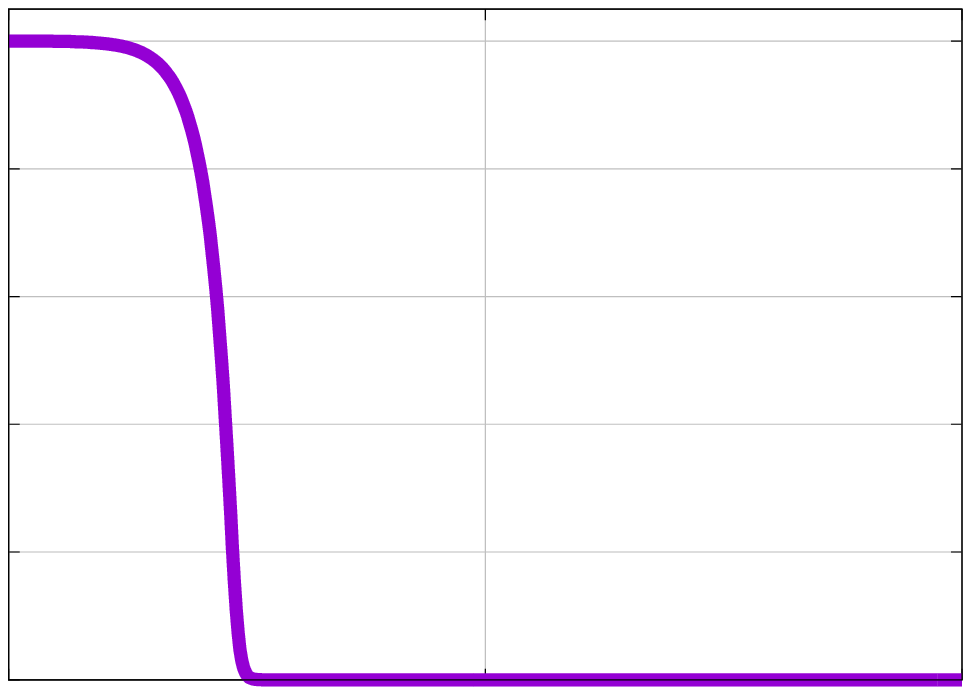}} & \resizebox{0.33\columnwidth}{!}{\Huge\input{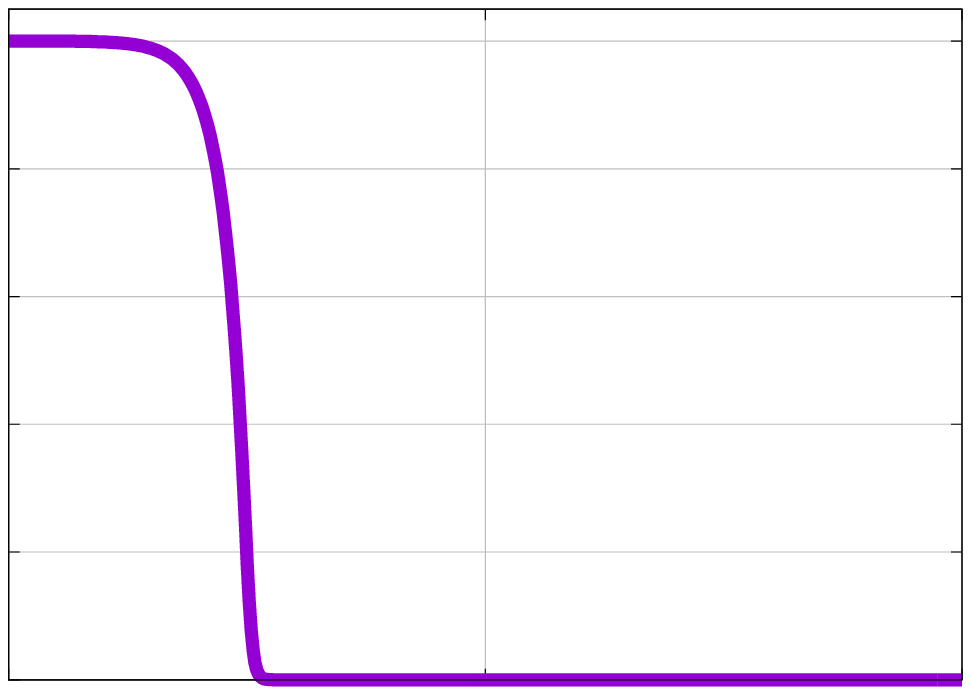}}\\
(j) mass fraction, $t = 7.5\times 10^{-6}$ & (k) mass fraction, $t = 2\times 10^{-5}$& (l) mass fraction, $t = 3\times 10^{-5}$
\end{tabular}
\caption{The states of the full model of the premixed flame problem. The sharp temperature and species gradients, and multiscale interactions between acoustics and flame, indicate that traditional static reduced models become inefficient.} 
\label{fig:FOMVis}
\end{figure}

\subsection{Numerical solver description}
PERFORM numerically solves the one-dimensional Navier-Stokes equations with chemical species transport and a chemical reaction source term:
\[
\frac{\partial}{\partial t}q(t, x) + \frac{\partial}{\partial x}\left(f(t, x) - f_v(t, x)\right) = f_s(t, x)\,,
\]
with
\begin{equation}
q = \begin{bmatrix}\rho\\\rho u\\\rho h^0 - p\\\rho Y_l\end{bmatrix}\,,\quad f = \begin{bmatrix}\rho u\\\rho u^2 + p\\\rho h^0 u\\\rho Y_l\end{bmatrix}\,,\quad f_v = \begin{bmatrix}0\\\tau \\u\tau - q\\ -\rho V_l Y_l\end{bmatrix}\,,\quad f_s = \begin{bmatrix}0\\0\\0\\\dot{\omega}_l\end{bmatrix}\,,
\label{eq:QuantitiesOfPDE}
\end{equation}
where $q$ is the conserved state at time $t$ and spatial coordinate $x$, $f$ is the inviscid flux vector, $f_v$ is the viscous flux vector, and $f_s$ is the source term. Additionally, $\rho$ is density, $u$ is velocity, $h^0$ is stagnation enthalpy, $p$ is static pressure, and $Y_l$ is mass fraction of the $l$th chemical species.
The reaction source $\dot{\omega}_l$ corresponds to the reaction model, which is described by an irreversible Arrhenius rate equation.
The problem is discretized in the spatial domain with a second-order accurate finite volume scheme. The inviscid flux is computed by the Roe scheme~\cite{Roe1981}. Gradients are limited by the Venkatakrishnan limiter~\cite{VENKATAKRISHNAN1993}. The time derivative is discretized with the first-order backwards differentiation formula scheme (i.e. backward Euler). Calculation of the viscous stress $\tau$, heat flux $q$, and diffusion velocity $V_l$, and any additional details about the implementation can be found in PERFORM's online documentation.

\subsection{Premixed flame with artificial forcing}\label{sec:PERFORM:Benchmark}
We consider a setup corresponding to a model premixed flame with artificial pressure forcing. There are two chemical species: ``reactant'' and ``product''. The reaction is a single-step irreversible mechanism that converts low-temperature reactant to high-temperature product, modeling a premixed combustion process. An artificial sinusoidal pressure forcing is applied at the outlet, which causes an acoustic wave to propagate upstream. The interaction between different length and time scales given by the system acoustics caused by the forcing and the flame leads to strongly nonlinear system dynamics with multiscale effects. The result are dynamics that evolve in high-dimensional spaces and thus are inefficient to reduce with static reduced models; see~Section~\ref{sec:ProblemFormulation}. The states of the full model and how they evolve over time are shown in Figure~\ref{fig:FOMVis}.

\section{Numerical results}\label{sec:NumExp}
We demonstrate nonlinear model reduction with online adaptive empirical interpolation on the model premixed flame problem introduced in Section~\ref{sec:PERFORM:Benchmark}. 

\begin{figure}[t]
\begin{tabular}{cc}
\includegraphics[width=0.48\linewidth]{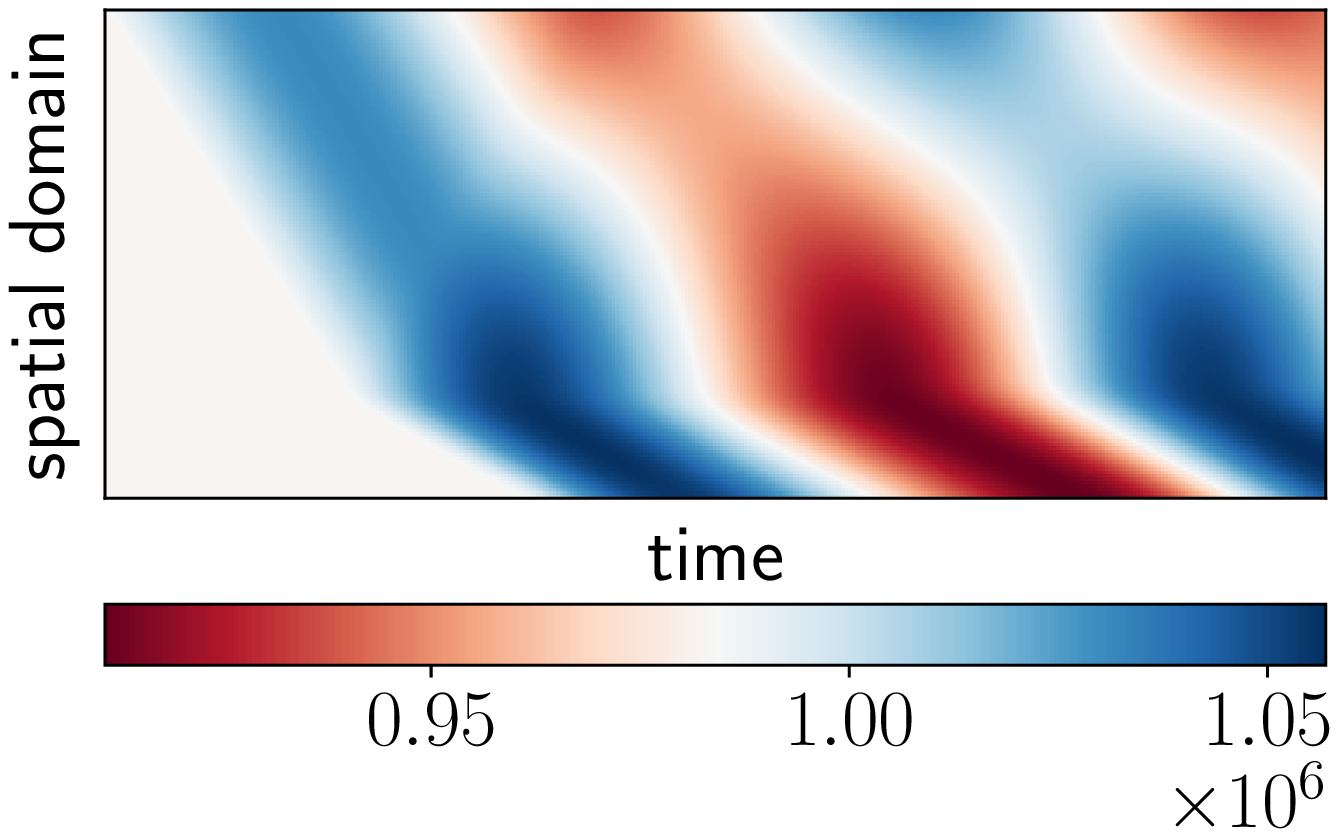} & \hspace*{-0.2cm}\includegraphics[width=0.48\linewidth]{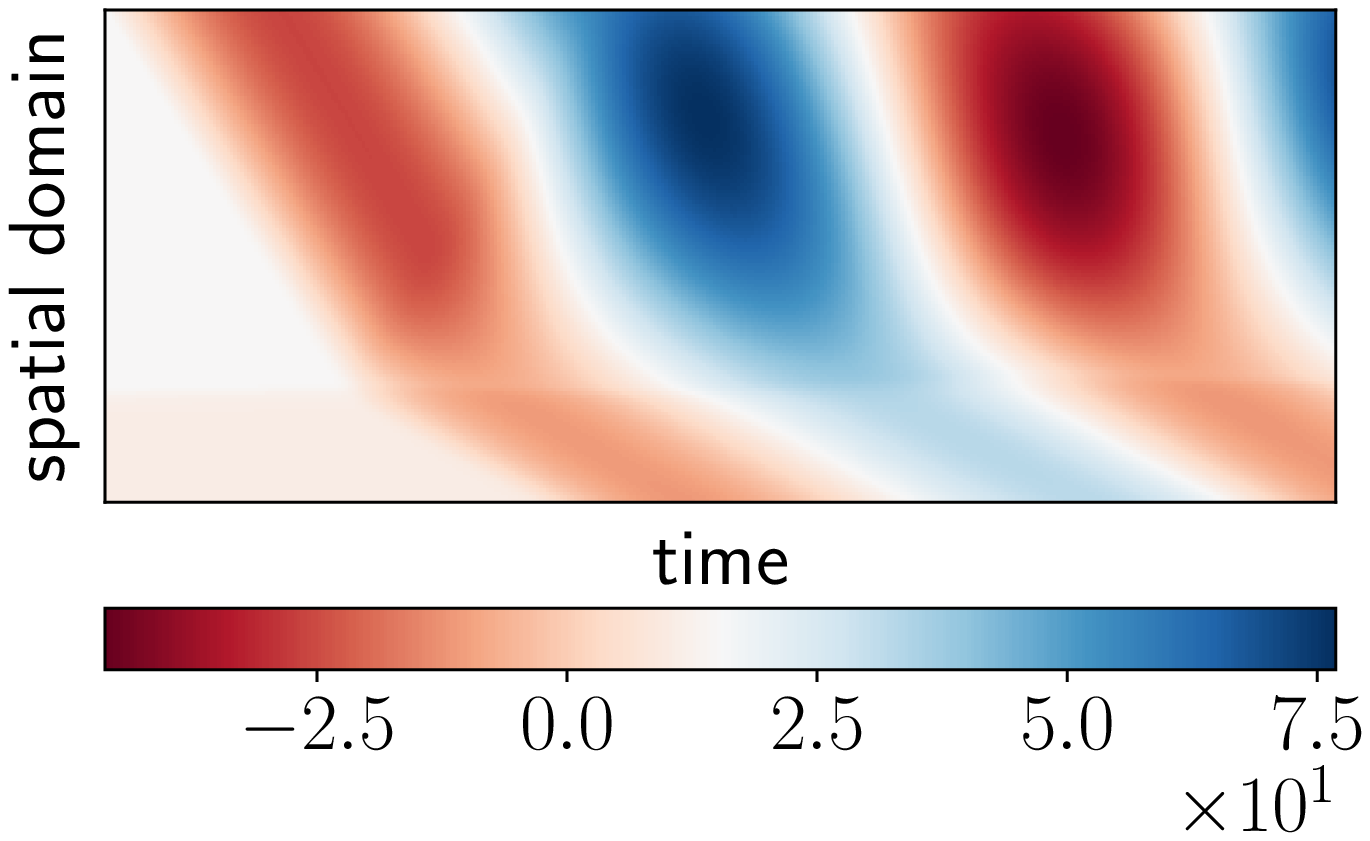}\\
(a) pressure & (b) velocity\\
\includegraphics[width=0.48\linewidth]{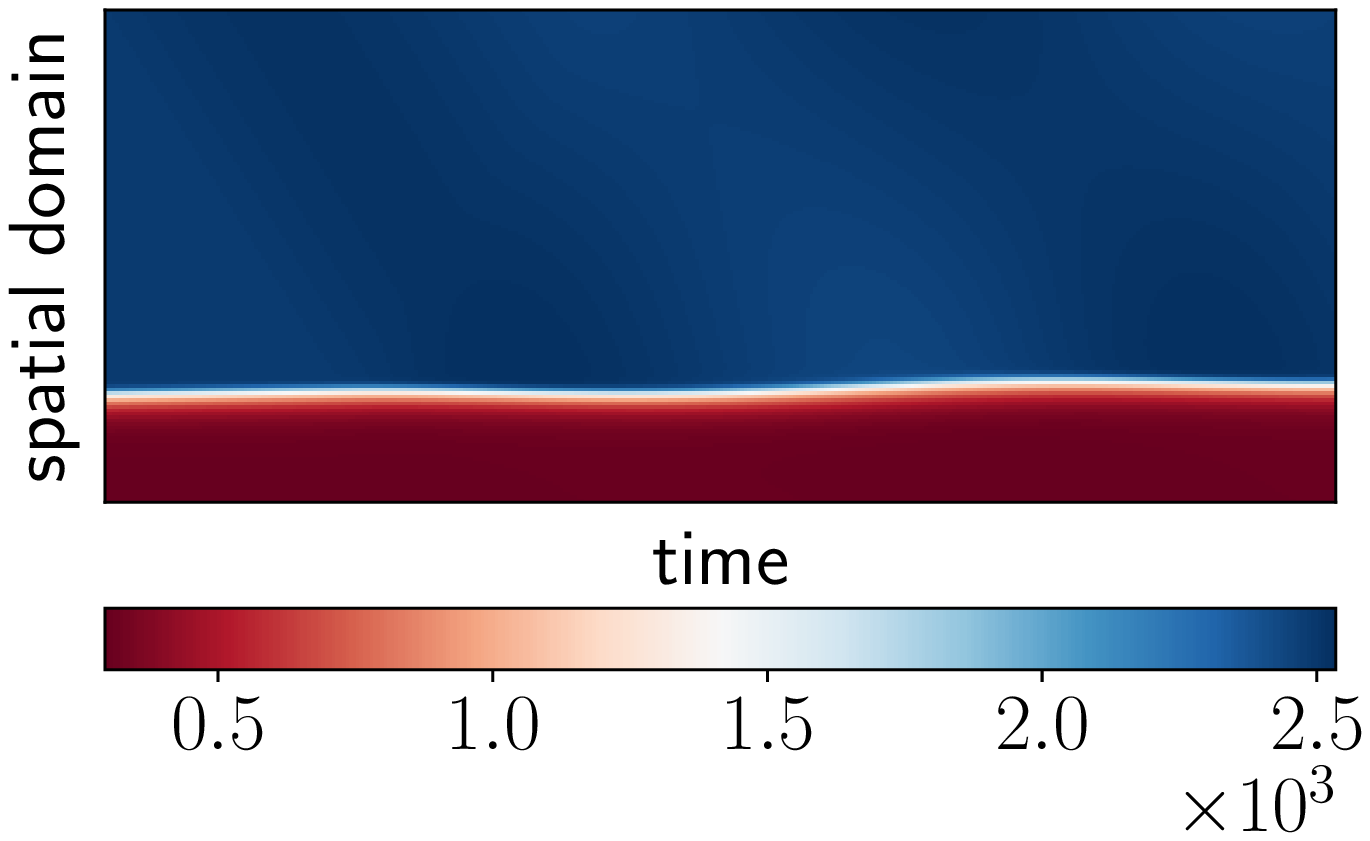} & \hspace*{-0.2cm}\includegraphics[width=0.48\linewidth]{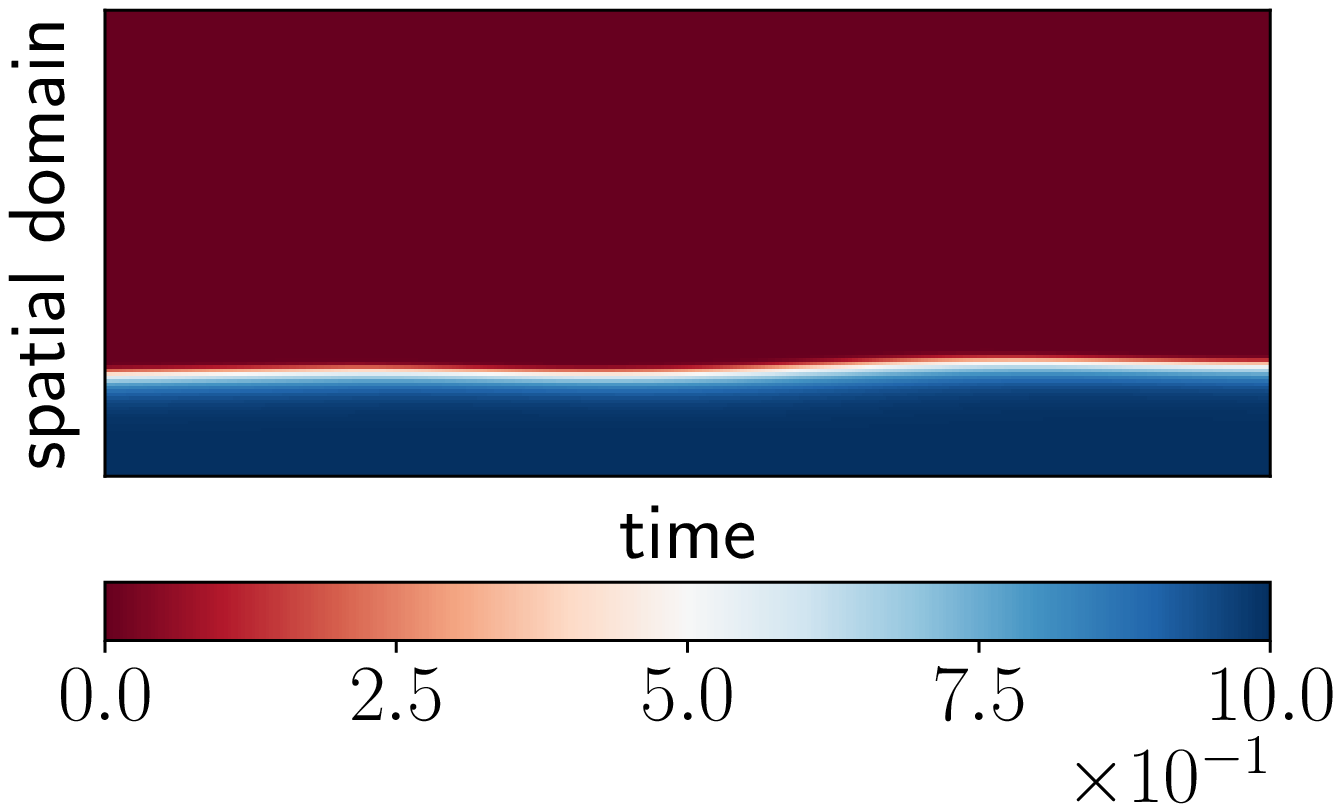}\\
(c) temperature & (d) species mass fraction
\end{tabular}
\caption{Full model: Plot (a) shows oscillations of pressure waves, which is evidence of the transport-dominated dynamics in this benchmark example.}
\label{fig:FOM2D}
\end{figure}

\begin{figure}[t]
\begin{tabular}{cc}
\resizebox{0.5\columnwidth}{!}{\LARGE\input{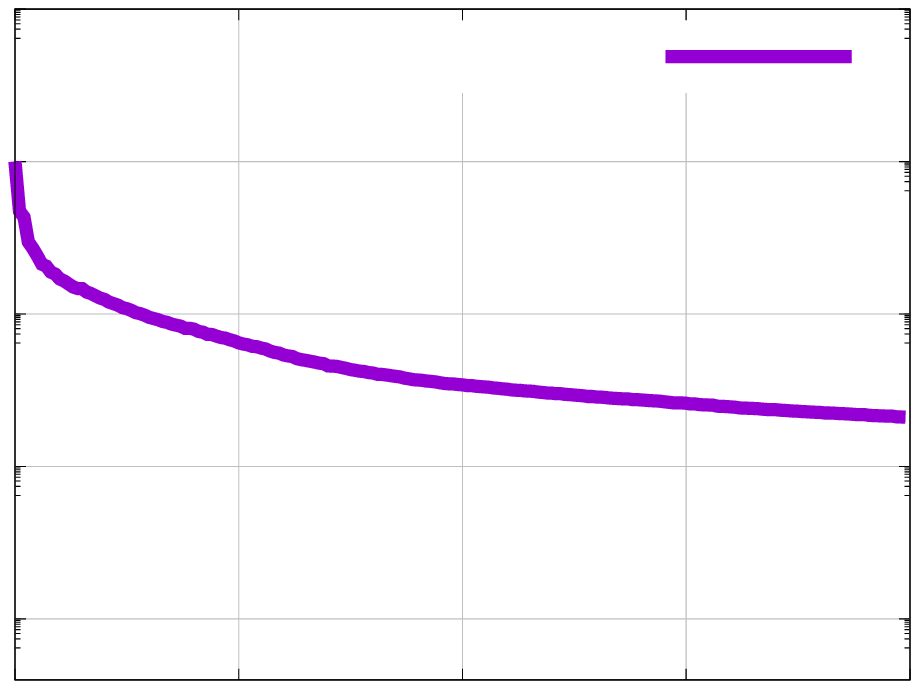}} &
\resizebox{0.5\columnwidth}{!}{\LARGE\input{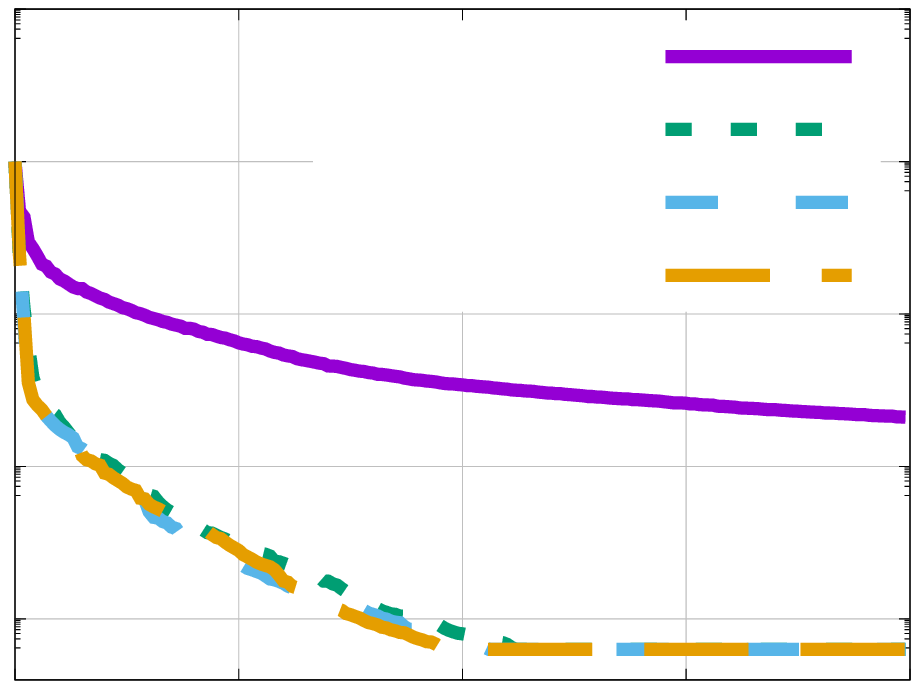}}\\
(a) global in time & (b) local in time
\end{tabular}
\caption{Singular values: Plot (a) shows the slow decay of the normalized singular values of the snapshots. Plot (b) shows that the normalized singular values computed from local trajectories in time decay faster.}
\label{fig:SingVal}
\end{figure}

\subsection{Numerical setup of full model}
Consider the problem described in Section~\ref{sec:PERFORM:Benchmark}. Each of the four conserved quantities $\rho$, $\rho u$, $\rho h^0 - p$, $\rho Y_1$ is discretized on $512$ equidistant grid points in the domain $\Omega = [0, 10^{-2}]$ m, which leads to a total of $\nh = 4 \times 512 = 2{,}048$ unknowns of the full model.
The time-step size is $\delta t =1 \times 10^{-9}$ s and the end time is $T = 3.5 \times 10^{-5}$ s, which is a total of $35{,}000$ time steps.
A 10\% sinusoidal pressure perturbation at a frequency of 50 kHz is applied at the outlet. 

Space-time plots of pressure, velocity, temperature, and species mass fraction are shown in Figure~\ref{fig:FOM2D}. Pressure and velocity shown in plots (a) and (b) of Figure~\ref{fig:FOM2D} show a transport-dominated behavior, which is in agreement with the pressure and velocity waves shown in Figure~\ref{fig:FOMVis}.

\begin{figure}[t]
\begin{tabular}{cc}
\includegraphics[width=0.48\linewidth]{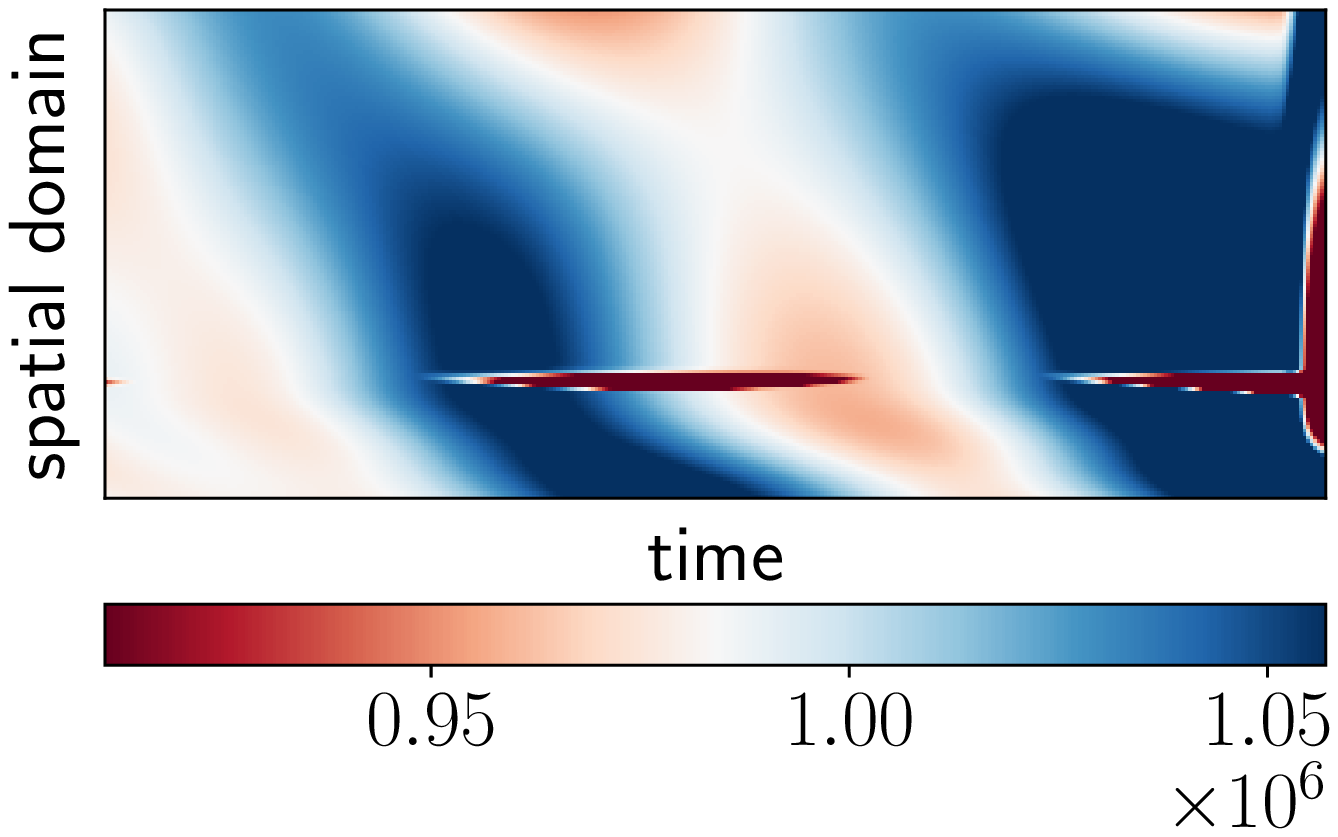} & \hspace*{-0.2cm}\includegraphics[width=0.48\linewidth]{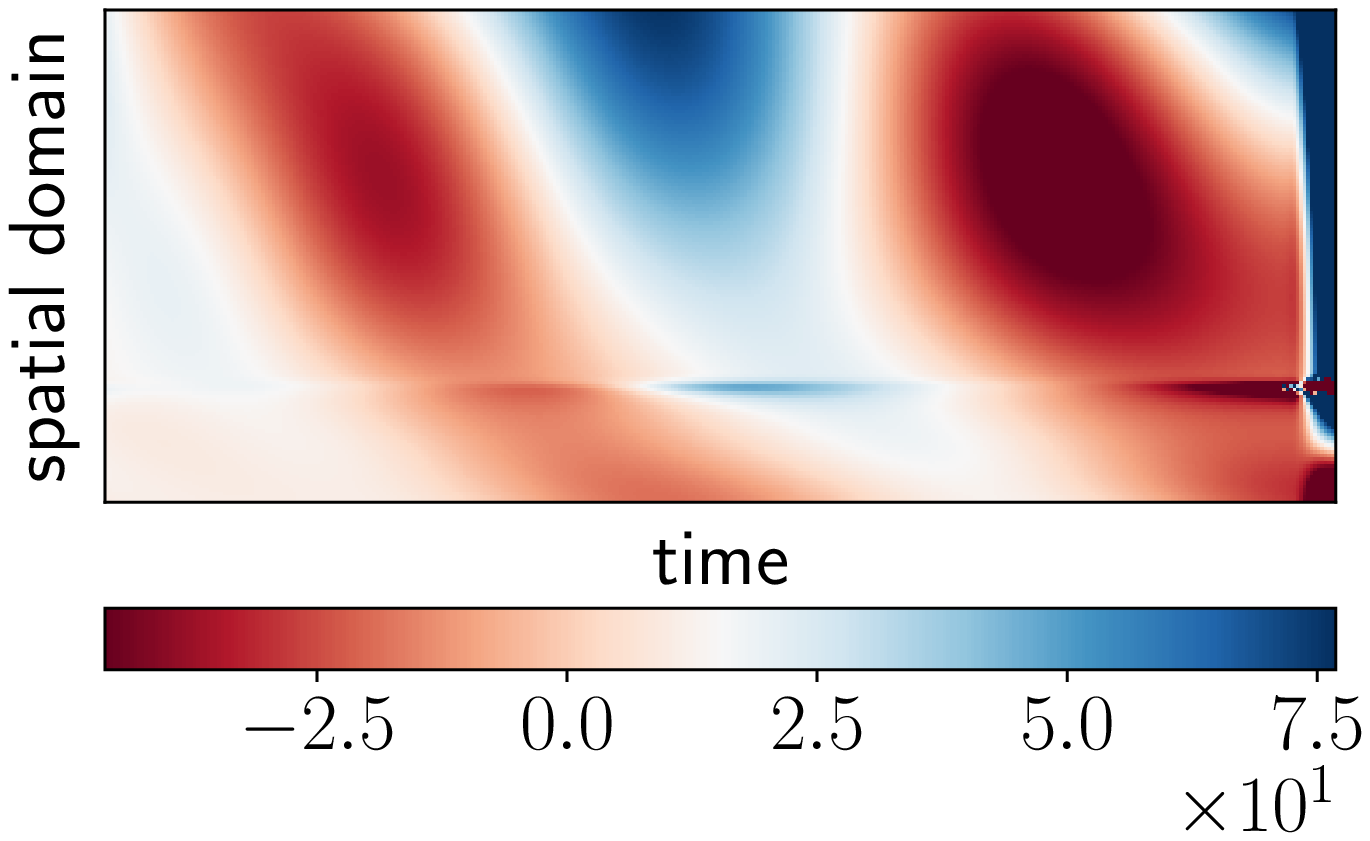}\\
(a) pressure & (b) velocity\\
\includegraphics[width=0.48\linewidth]{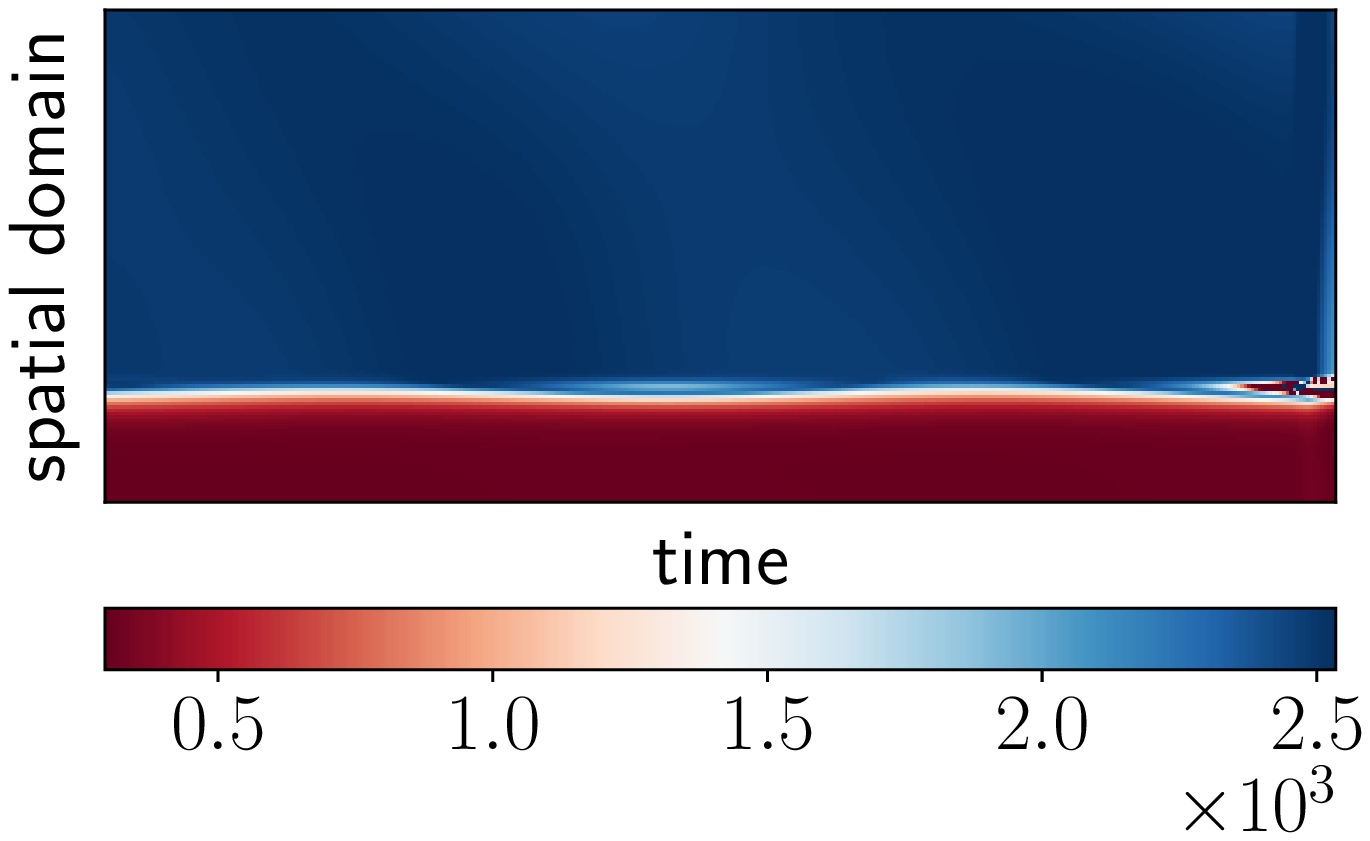} & \hspace*{-0.2cm}\includegraphics[width=0.48\linewidth]{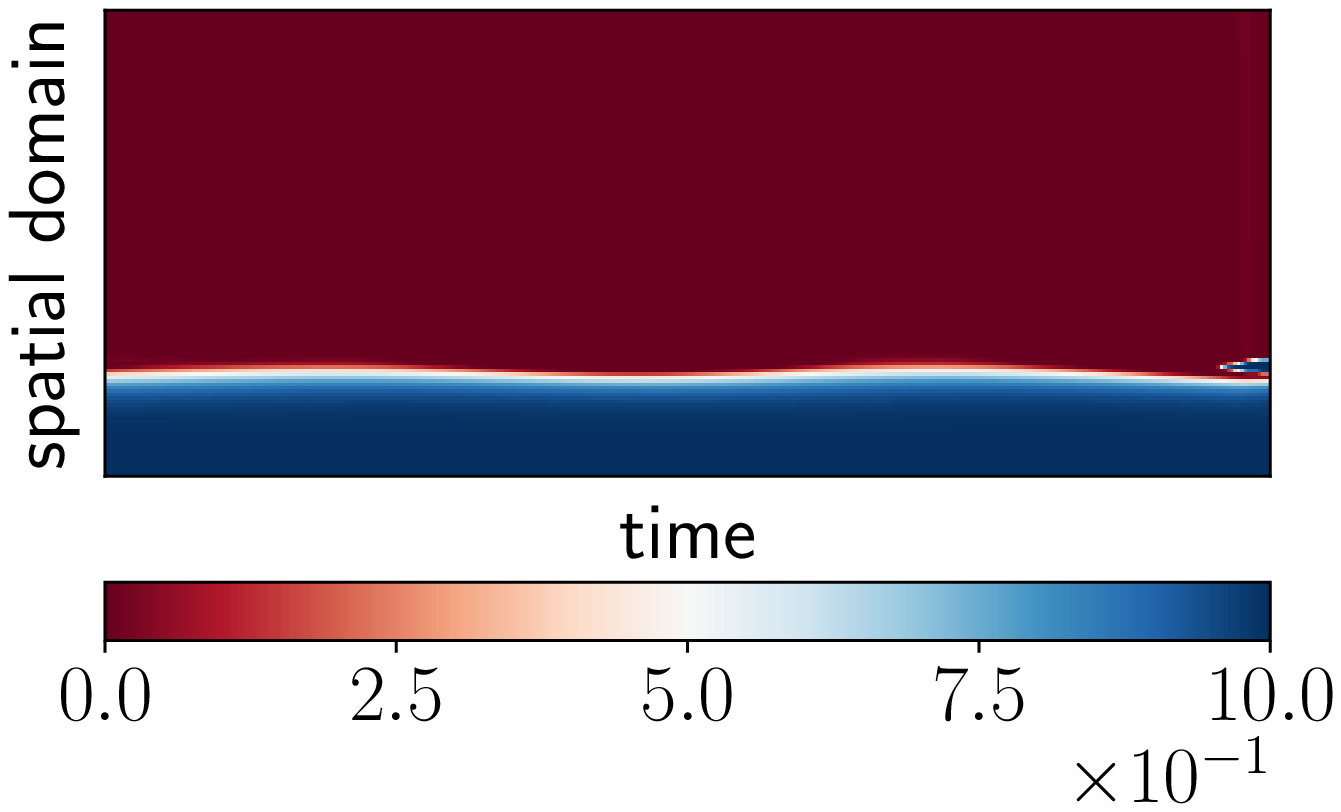}\\
(c) temperature & (d) species mass fraction
\end{tabular}
\caption{Static reduced model of dimension $\nr = 6$: Because of the coupling of dynamics over various length scales in the premixed flame example, a static reduced model is insufficient to provide accurate approximations of the full-model dynamics shown in Figure~\ref{fig:FOM2D}. The time-space plot corresponding to the static reduced model with $\nr = 7$ is not shown here because it provides a comparably poor approximation.}
\label{fig:StaticSpaceTimeD6}
\end{figure}

\subsection{Static reduced models}\label{sec:NumExp:StaticROM}
Snapshots are generated with the full model by querying \eqref{eq:FOM} for 35,000 time steps and storing the full model state every 50 time steps. The singular values of the snapshot matrix are shown in Figure~\ref{fig:SingVal}a, which decay slowly and indicate that static reduced models are inefficient. From the snapshots, we generate POD bases $\bfV$ of dimension $\nr = 6$ and $\nr = 7$ to derive reduced models with empirical interpolation. The interpolation points are selected with QDEIM \cite{QDEIM}. The empirical interpolation points are computed separately with QDEIM for each of the four variables in \eqref{eq:QuantitiesOfPDE}. Then, we select a component of the state $\bfq$ as an empirical interpolation point if it is selected for at least one variable.

The time-space plots of the static reduced approximation of the full-model dynamics are shown in Figure~\ref{fig:StaticSpaceTimeD6} for dimension $\nr = 6$. The approximation is poor, which is in agreement with the transport-dominated dynamics and the slow decay of the singular values shown in Figure~\ref{fig:SingVal}a. The time-space plot of the static reduced model of dimension $\nr = 7$ gives comparably poor approximations and is not shown here. It is important to note that the static reduced model is derived from snapshots over the whole time range from $t = 0$ to end time $T = 3.5 \times 10^{-5}$, which means that the static reduced model has to merely reconstruct the dynamics that were seen during training, rather than predicting unseen dynamics. This is in stark contrast to the  adaptive reduced model derived with AADEIM in the following subsection, where the reduced model will predict states far outside of the training window.

\begin{figure}[t]
\begin{tabular}{cc}
\includegraphics[width=0.48\linewidth]{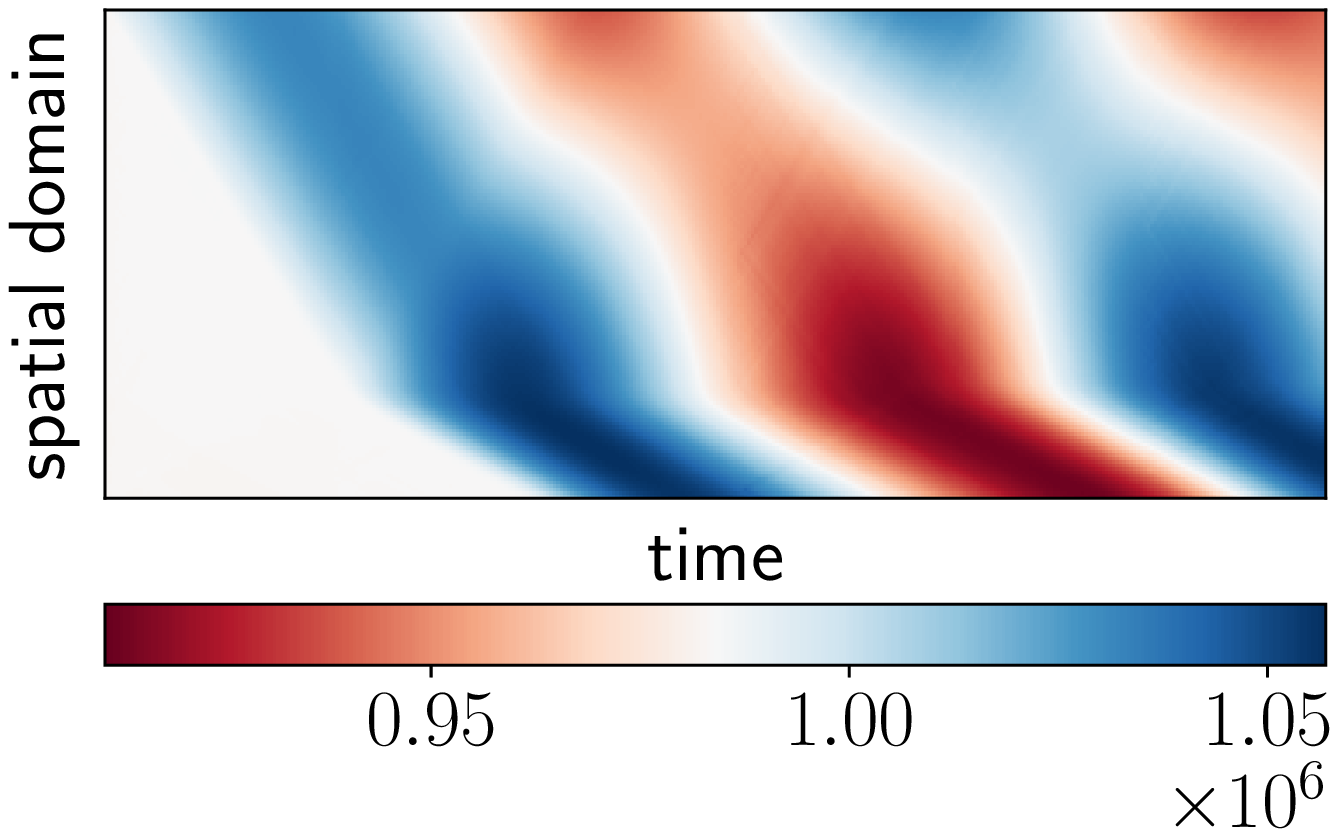} & \hspace*{-0.2cm}\includegraphics[width=0.48\linewidth]{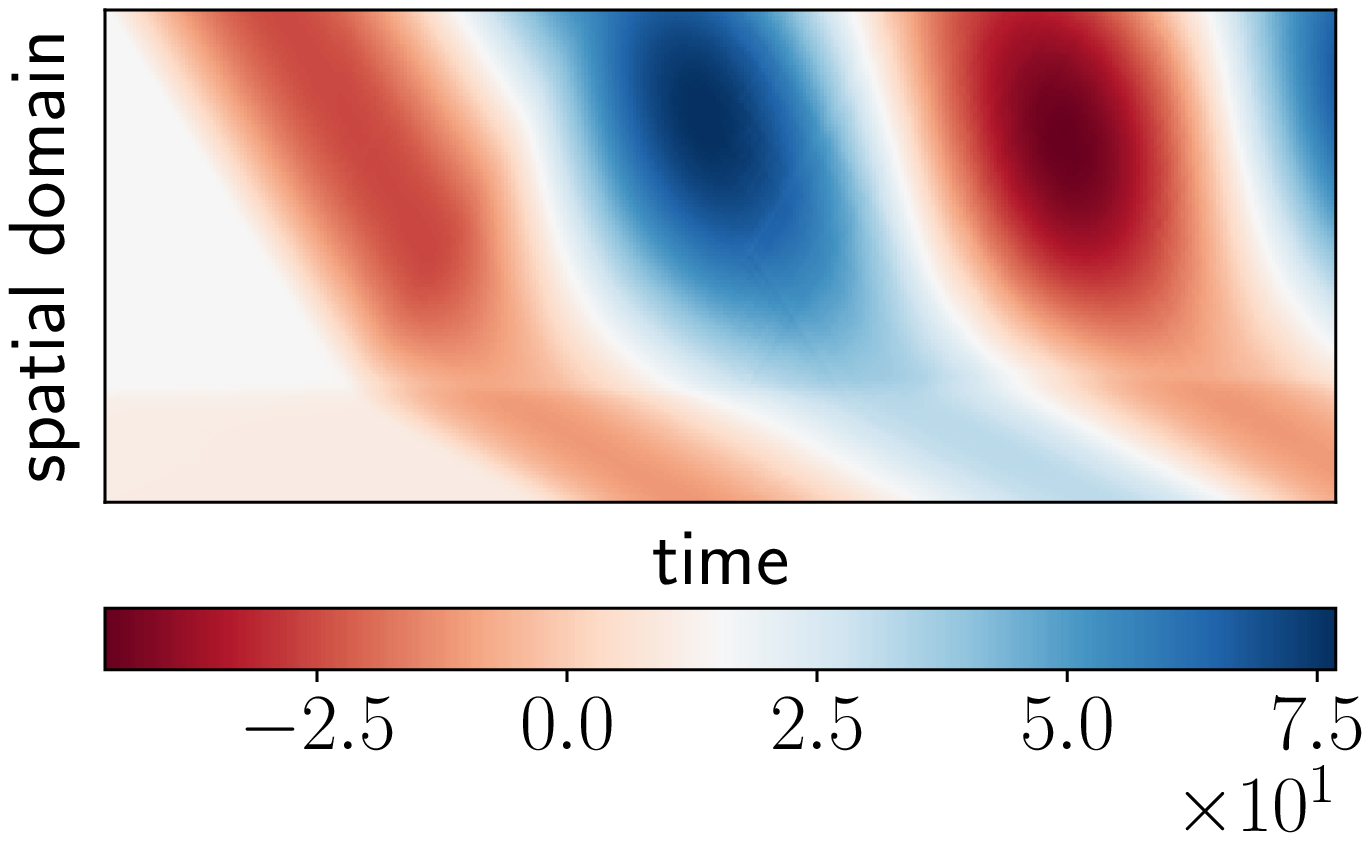}\\
(a) pressure & (b) velocity\\
\includegraphics[width=0.48\linewidth]{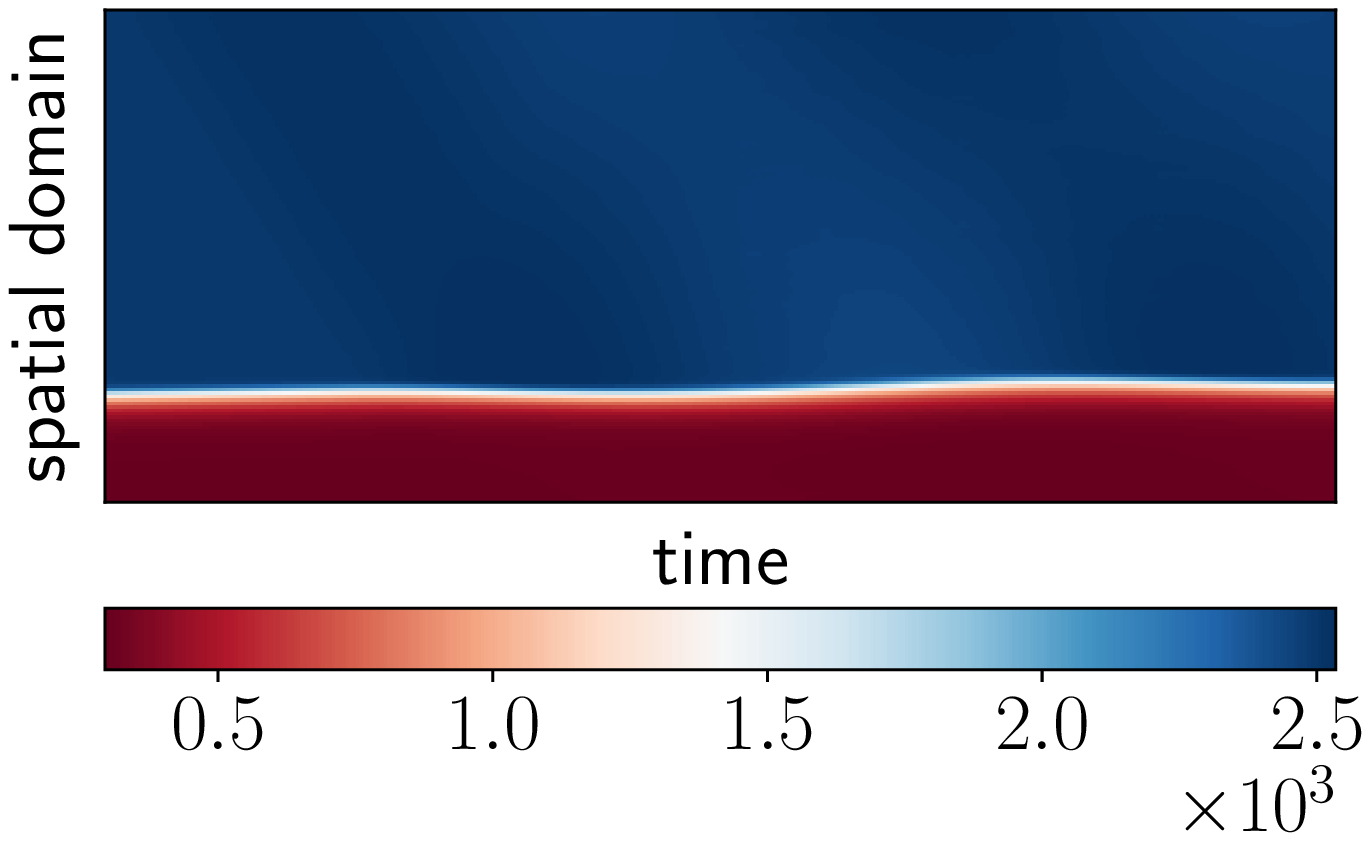} & \hspace*{-0.2cm}\includegraphics[width=0.48\linewidth]{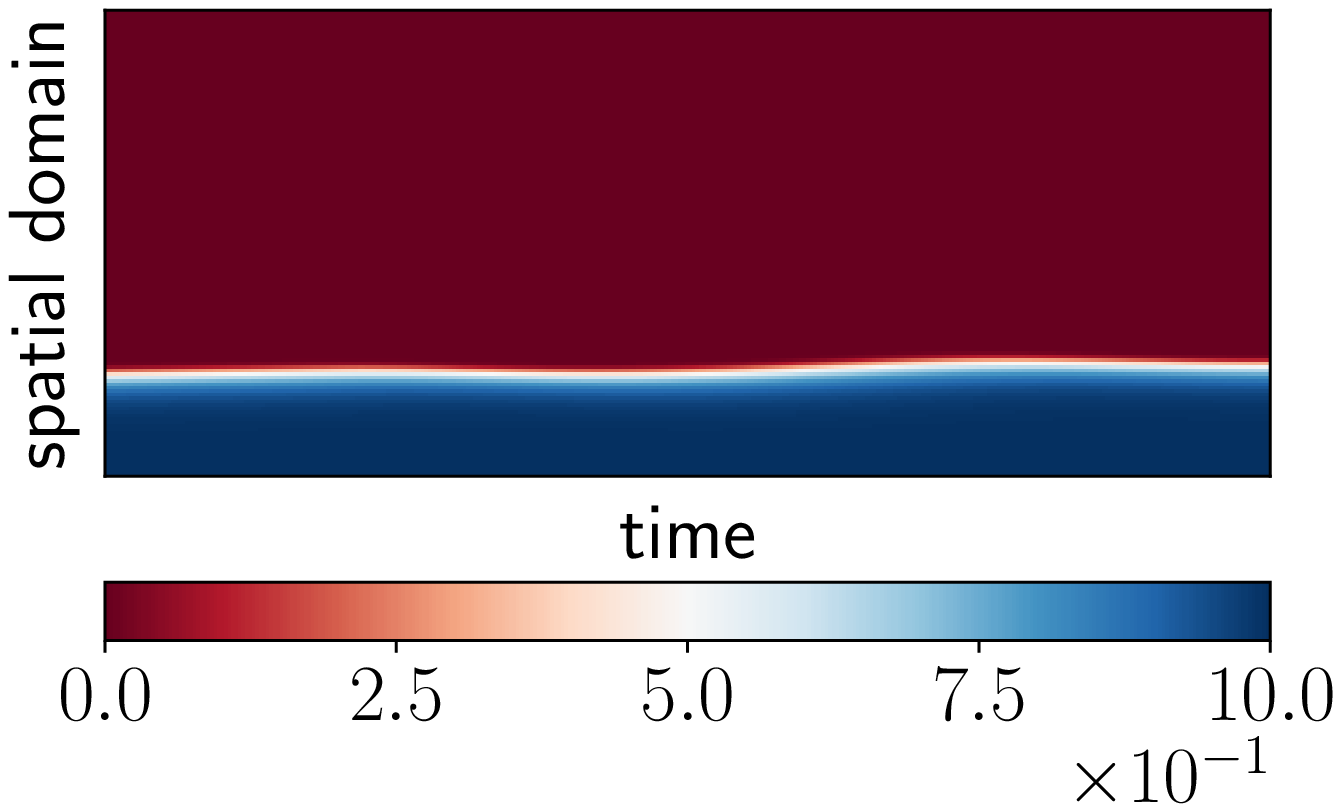}\\
(c) temperature & (d) species mass fraction
\end{tabular}
\caption{The online adaptive reduced model with $\nr = 6$ dimensions obtained with AADEIM provides accurate predictions of the full-model dynamics shown in Figure~\ref{fig:FOM2D}. The training snapshots used for initializing the AADEIM model cover dynamics up to time $t = 1.6 \times 10^{-8}$ and thus all states later in time up to end time $T = 3.5 \times 10^{-5}$ are predictions outside of the training data.}
\label{fig:ADEIM62DPlotHigh}
\end{figure}

\subsection{Reduced model with AADEIM}
We derive a reduced model with AADEIM of dimension $\nr = 6$. The initial window size is $\winit = 15$ and the window size is $w = \nr + 1 = 7$, as recommended in \cite{P18AADEIM}. Notice that an initial window $\winit = 15$ means that the AADEIM model predicts unseen dynamics (outside of training data) starting at time step $k = 16$, which corresponds to $t = 1.6 \times 10^{-8}$. This is in stark contrast to Section~\ref{sec:NumExp:StaticROM} where the static reduced model only has to reconstruct seen dynamics. The number of sampling points is $\nrs = 1{,}024$ and the frequency of adapting the sampling points is $\nz = 3$. This means that the sampling points are adapted every third time step; see~Algorithm~\ref{alg:ABAS}. The basis matrix $\bfV_k$ and the points matrix $\bfP_k$ are adapted every other time step.

\begin{figure}[t]
\begin{tabular}{cc}
\resizebox{0.5\columnwidth}{!}{\huge\input{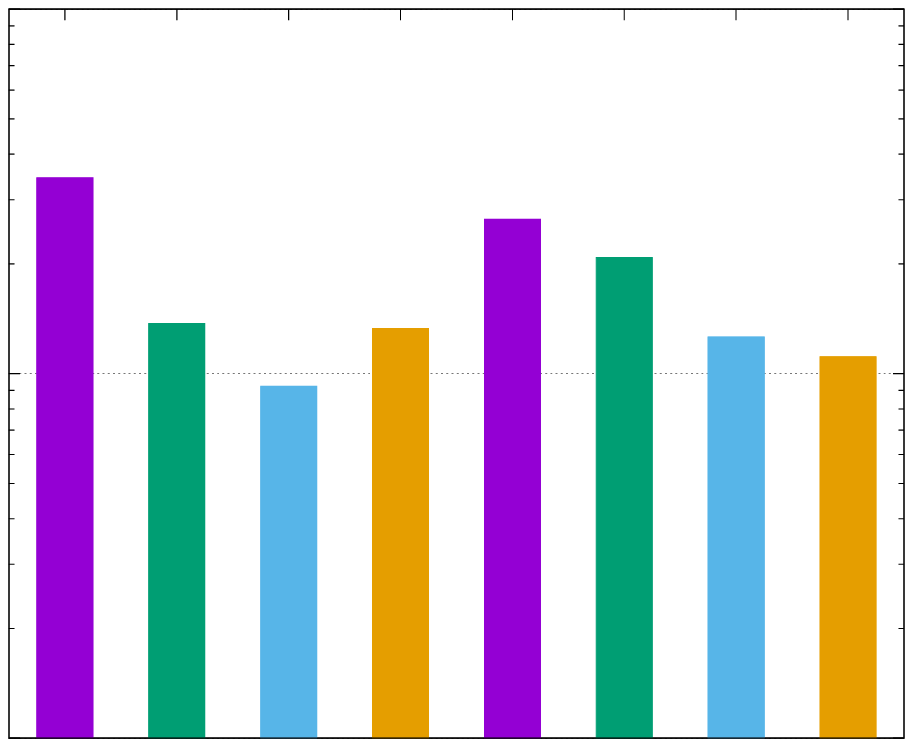}} &
\resizebox{0.5\columnwidth}{!}{\huge\input{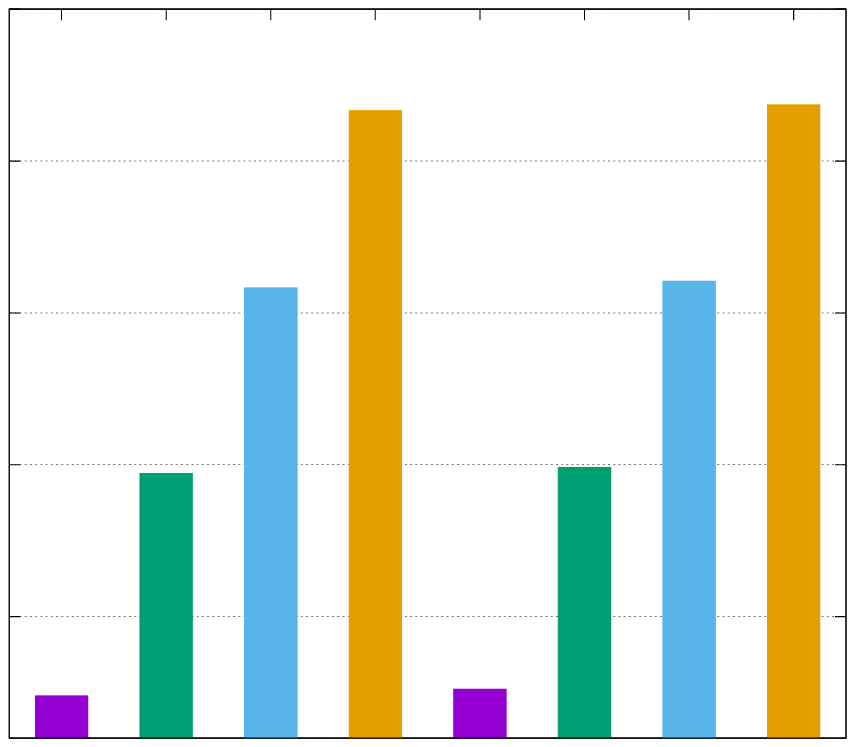}}\\
(a) error & (b) costs
\end{tabular}\\\centering\begin{minipage}{1.0\columnwidth}\fbox{\begin{tabular}{rl}
A: & dimension $n = 6$, update frequency $\nz = 4$, \#sampling points $\nrs = 768 $  \\
B: & dimension $n=6$, update frequency $\nz = 3$, \#sampling points $\nrs = 768$
\\
C: & dimension $n=6$, update frequency $\nz = 4$, \#sampling points $\nrs = 1024$
\\
D: & dimension $n=6$, update frequency $\nz = 3$, \#sampling points $\nrs = 1024$
\\
E: & dimension $n=7$, update frequency $\nz = 4$, \#sampling points $\nrs = 768$  \\
F: & dimension $n=7$, update frequency $\nz = 3$, \#sampling points $\nrs = 768$
\\
G: & dimension $n=7$, update frequency $\nz = 4$, \#sampling points $\nrs = 1024$
\\
H: & dimension $n=7$, update frequency $\nz = 3$, \#sampling points $\nrs = 1024$
\end{tabular}}~\\~\\~\\\end{minipage}
\caption{Performance of AADEIM reduced models for dimensions $\nr = 6$ and $\nr = 7$ with various combinations of number of sampling points $\nrs$ and update frequency $\nz$. One observation is that the AADEIM models outperform the static reduced models in terms of error \eqref{eq:ErrorPlot} over all parameters. Another observation is that a higher number of sampling points $\nrs$ tends to lead to lower errors in the AADEIM predictions in favor of higher costs.}
\label{fig:ErrorVsCosts}
\end{figure}

\begin{figure}[t]
\begin{tabular}{cc}
\resizebox{0.5\columnwidth}{!}{\huge\input{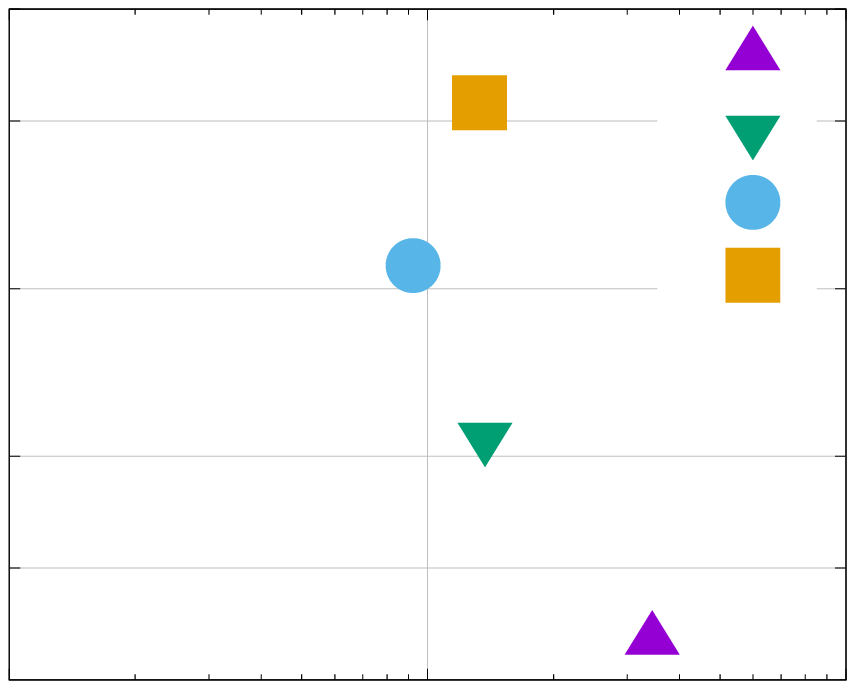}} &
\resizebox{0.5\columnwidth}{!}{\huge\input{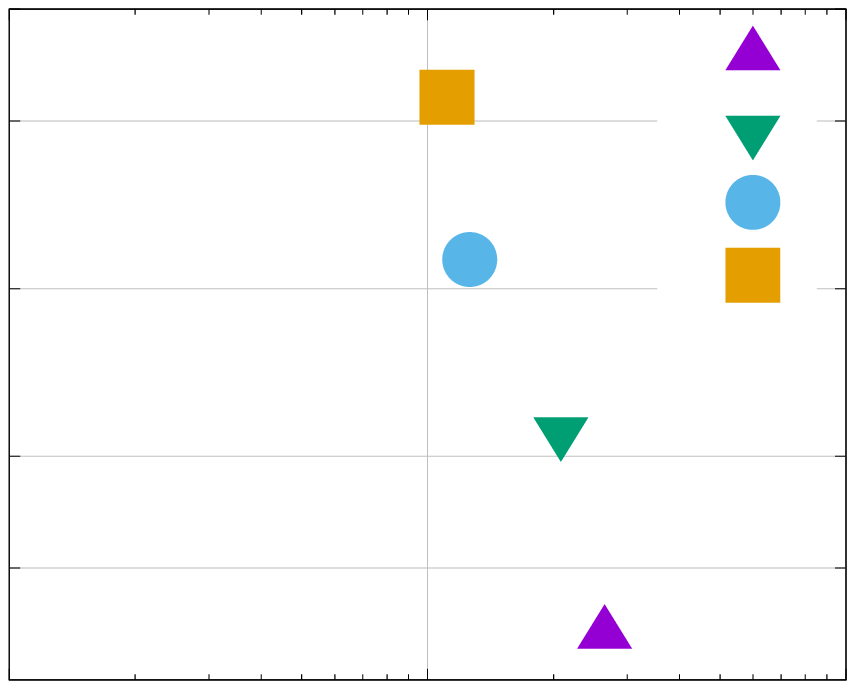}}\\
(a) dimension $\nr = 6$ & (b) dimension $\nr = 7$
\end{tabular}
\caption{Cost vs.~error of AADEIM models for dimension $\nr = 6$ and $\nr = 7$. All AADEIM models achieve a prediction error \eqref{eq:ErrorPlot} of roughly $10^{-3}$, which indicates an accurate prediction of future-state dynamics and which is in contrast to the approximations obtained with static reduced models of the same dimension. See Figure~\ref{fig:ErrorVsCosts} for legend.}
\label{fig:ErrorVsCostsScatter}
\end{figure}

The time-space plot of the prediction made with the AADEIM model is shown in Figure~\ref{fig:ADEIM62DPlotHigh}. The predicted states obtained with the AADEIM model are in close agreement with the full model (Figure~\ref{fig:FOM2D}), in contrast to the states derived with the static reduced model (Figure~\ref{fig:StaticSpaceTimeD6}).
This also in agreement with the fast decay of the singular values of snapshots in local time windows, as shown in Figure~\ref{fig:SingVal}b.

We further consider AADEIM models with dimension $\nr = 7$, initial window $\winit=12$, and $\nrs = 768$ sampling points and frequency of adapting the sampling points $\nz = 4$. We compute the average relative error as
\begin{equation}\label{eq:ErrorPlot}
e = \frac{\|\tilde{\bfQ} - \bfQ\|_F^2}{\|\bfQ\|_F^2}\,,
\end{equation}
where $\bfQ$ is the trajectory obtained with the full model and $\tilde{\bfQ}$ is the reduced trajectory obtained with AADEIM from Algorithm~\ref{alg:ABAS}. All combinations of models and their performance in terms of average relative error are shown in Figure~\ref{fig:ErrorVsCosts}a. As costs, we count the number of components of the full-model right-hand side function that need to be evaluated and report them in  Figure~\ref{fig:ErrorVsCosts}b; see also Figure~\ref{fig:ErrorVsCostsScatter}. All online adaptive reduced models achieve a comparable error of $10^{-3}$, where a higher number of sampling points $\nrs$ and a more frequent adaptation $\nz$ of the sampling points typically leads to lower errors in favor of higher costs. These numerical observations are in agreement with the principles of AADEIM and the results shown in \cite{Peherstorfer15aDEIM,P18AADEIM,CKMP19ADEIMQuasiOptimalPoints}. 

\begin{figure}[p]
\begin{tabular}{ccc}
\resizebox{0.33\columnwidth}{!}{\Huge\input{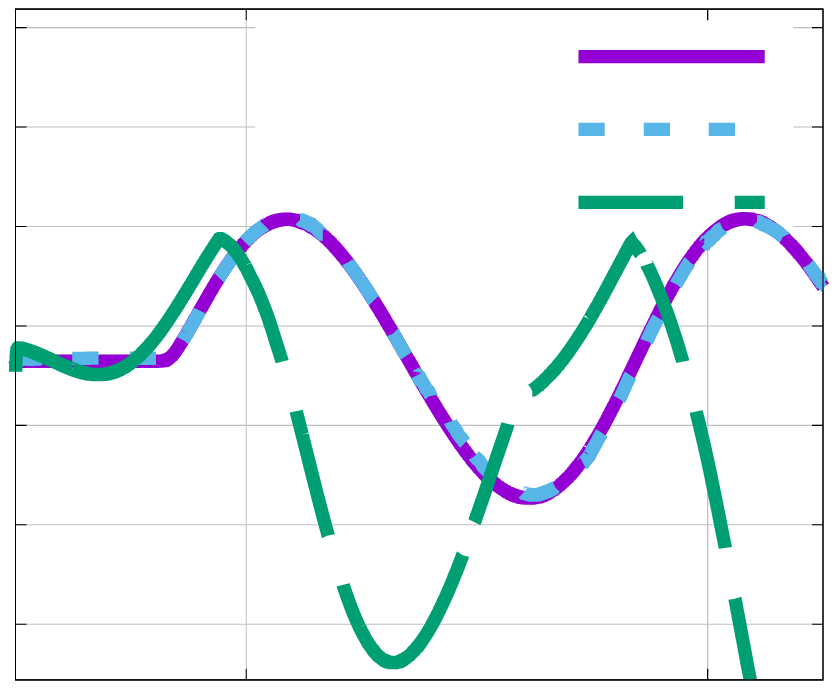}} & \resizebox{0.33\columnwidth}{!}{\Huge\input{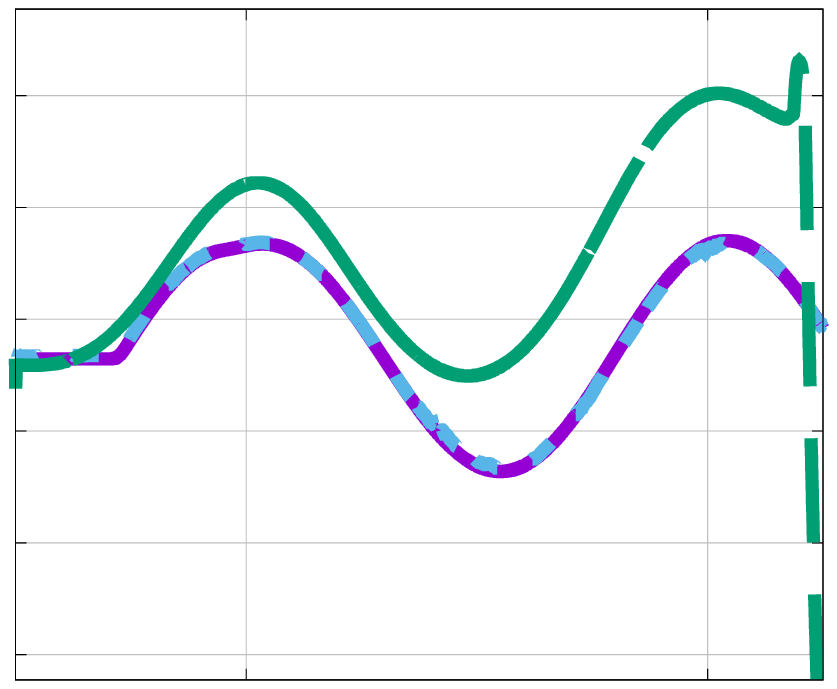}} & \resizebox{0.33\columnwidth}{!}{\Huge\input{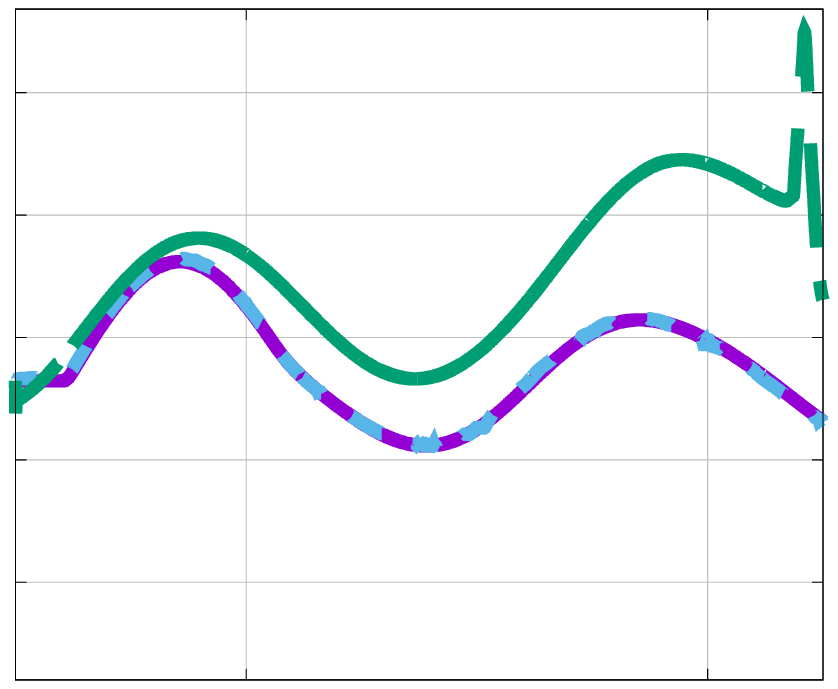}}\\
(a) pressure, probe 1 & (b) pressure, probe 2 & (c) pressure, probe 3\\
\resizebox{0.33\columnwidth}{!}{\Huge\input{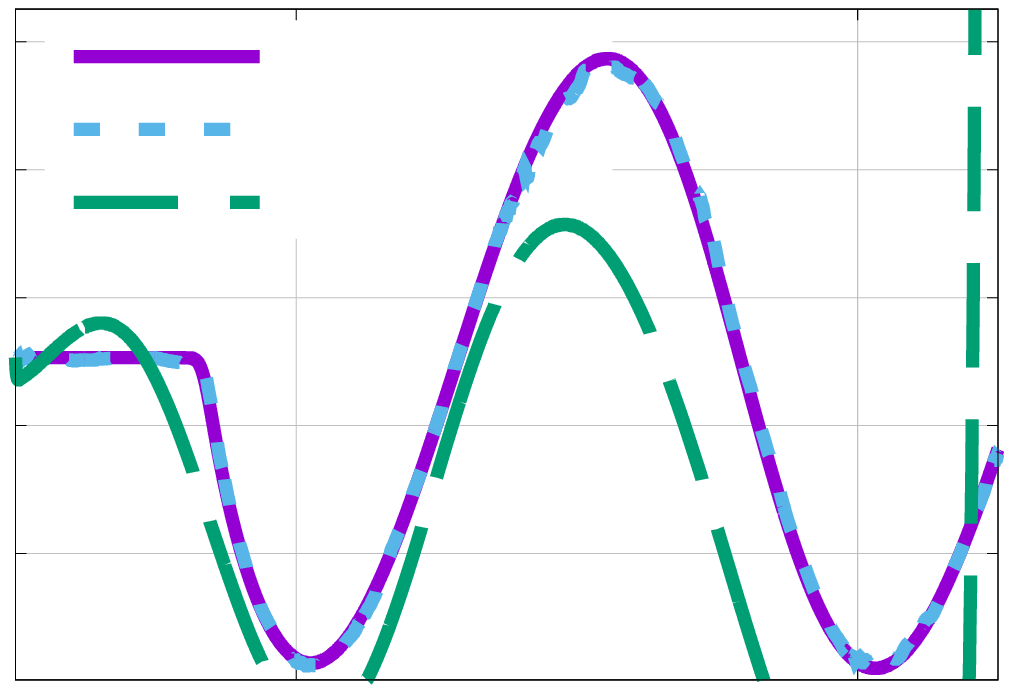}} & \resizebox{0.33\columnwidth}{!}{\Huge\input{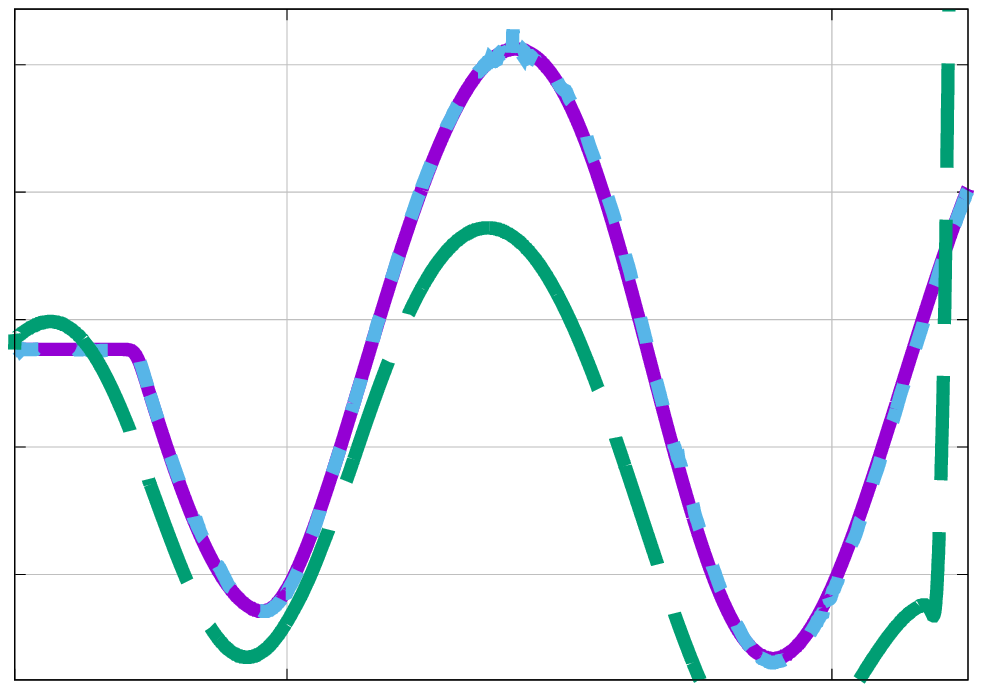}} & \resizebox{0.33\columnwidth}{!}{\Huge\input{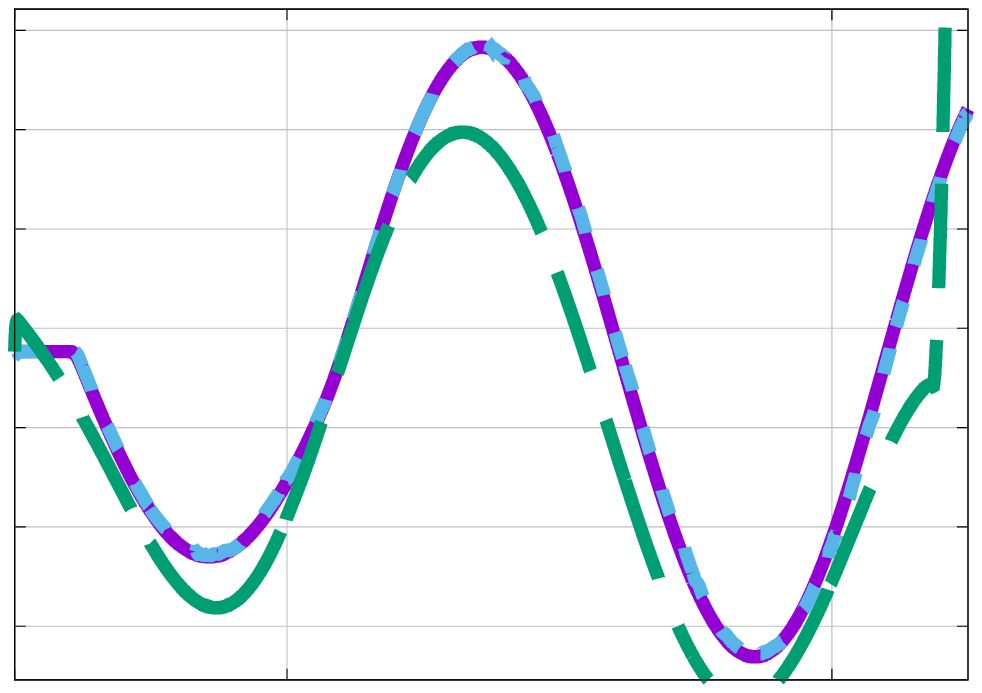}}\\
(d) velocity, probe 1 & (e) velocity, probe 2 & (f) velocity, probe 3\\
\resizebox{0.33\columnwidth}{!}{\Huge\input{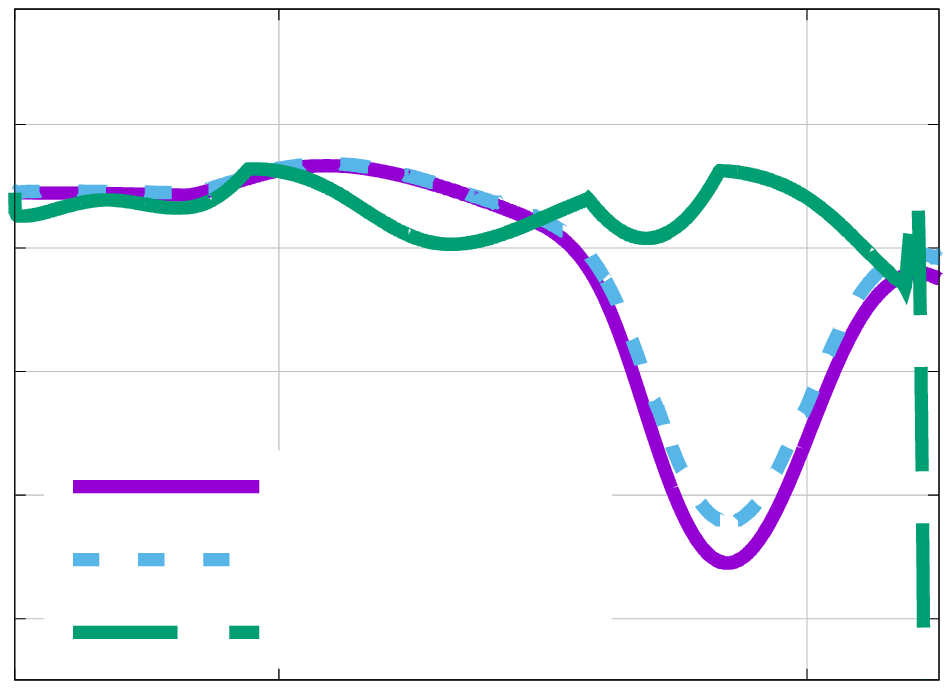}} & \resizebox{0.33\columnwidth}{!}{\Huge\input{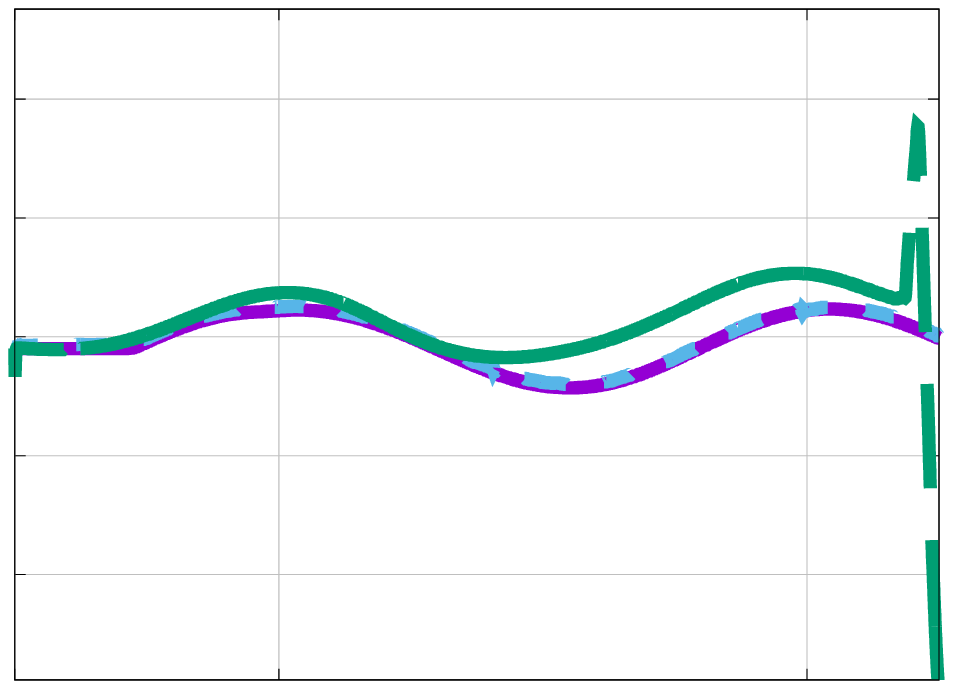}} & \resizebox{0.33\columnwidth}{!}{\Huge\input{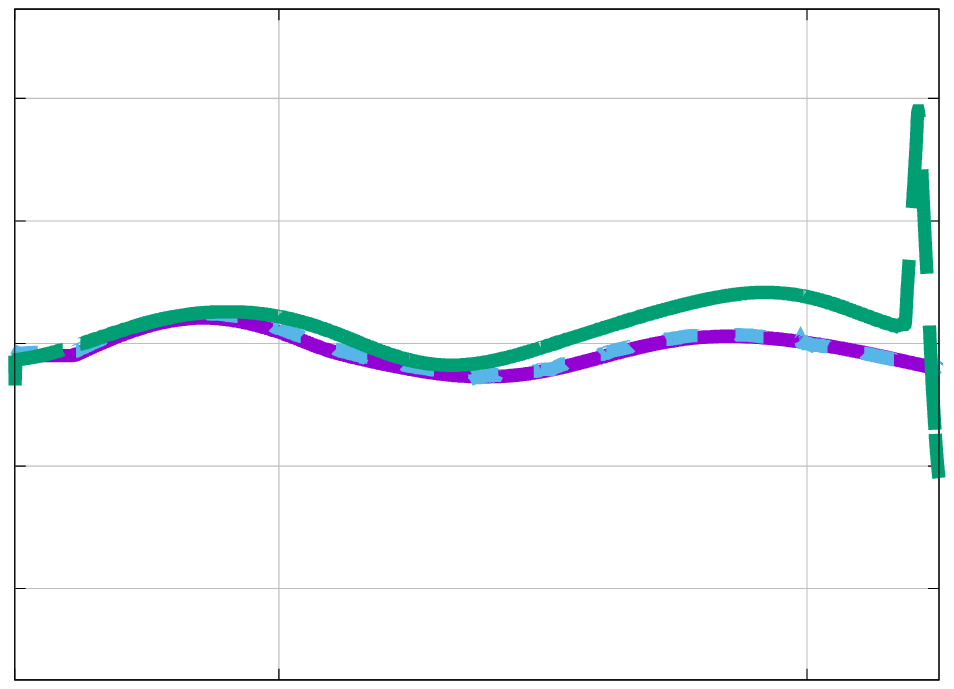}}\\
(g) temperature, probe 1 & (h) temperature, probe 2 & (i) temperature, probe 3\\
\resizebox{0.33\columnwidth}{!}{\Huge\input{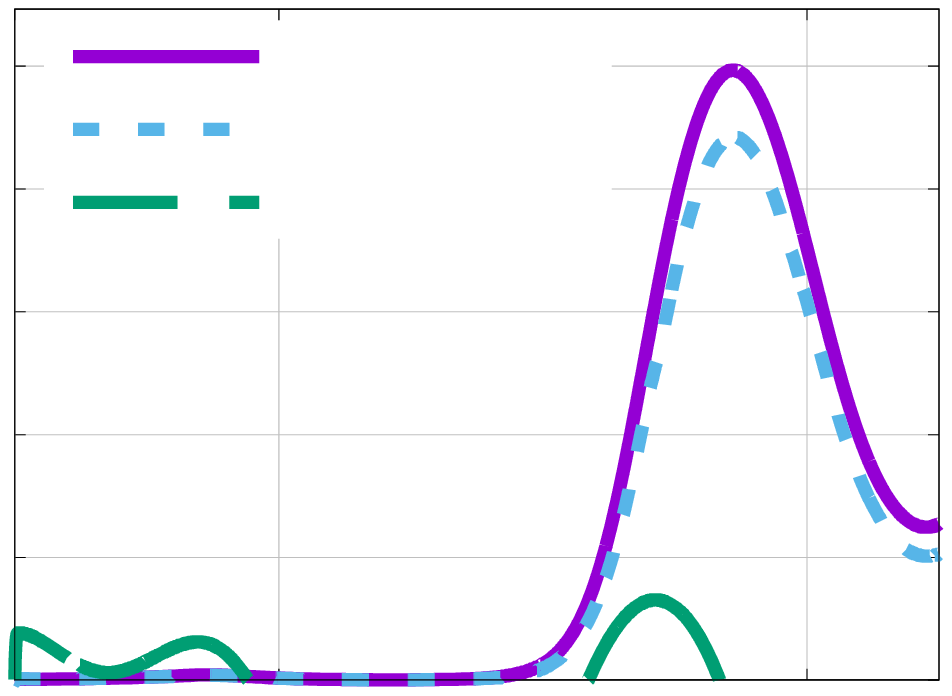}} & \resizebox{0.33\columnwidth}{!}{\Huge\input{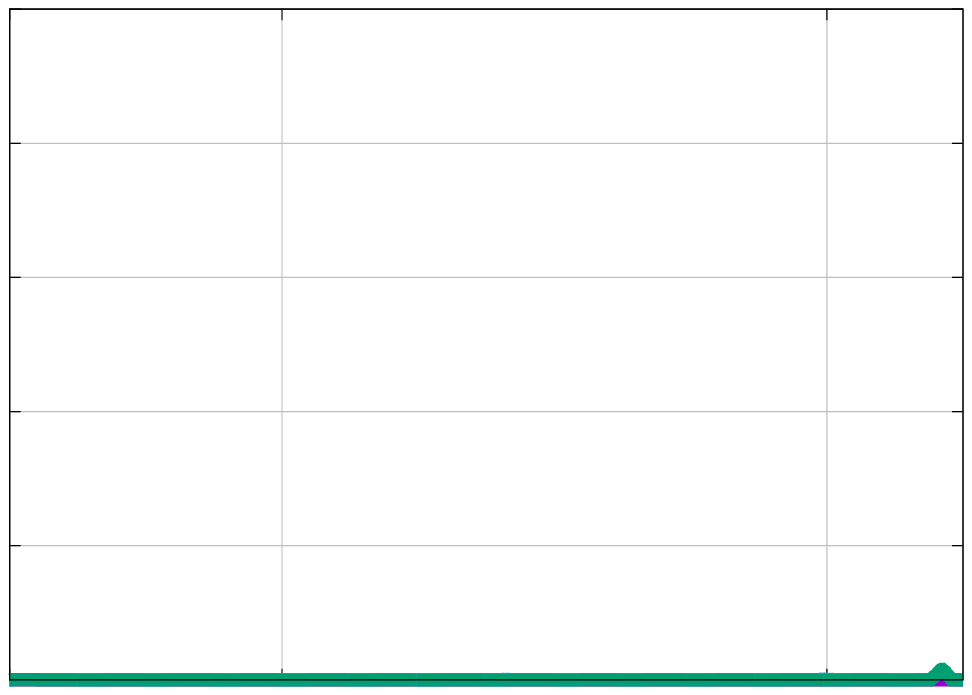}} & \resizebox{0.33\columnwidth}{!}{\Huge\input{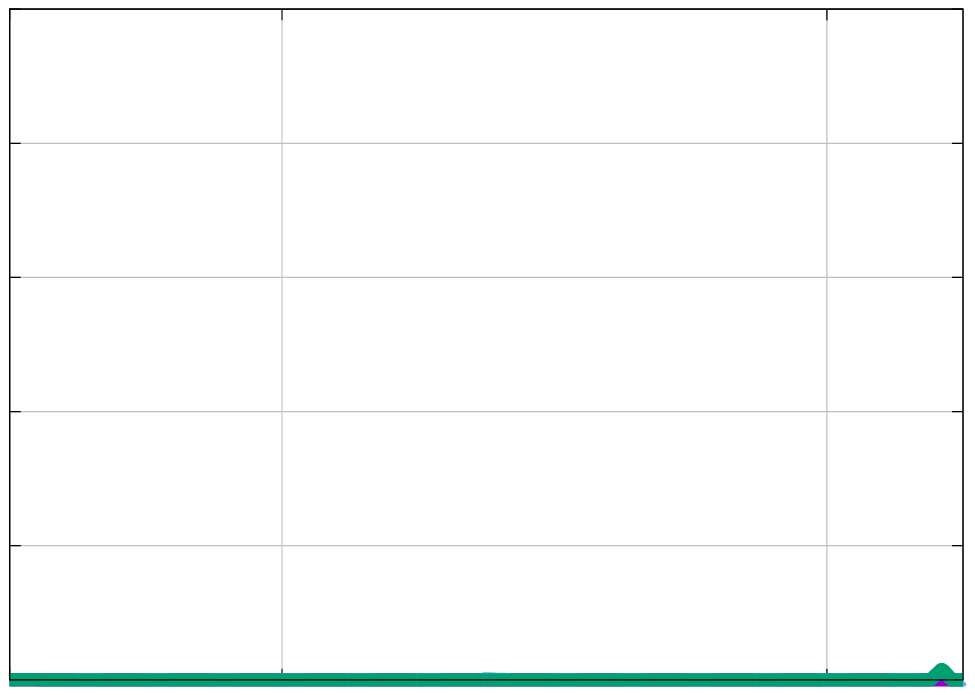}}\\
(j) mass fraction, probe 1 & (k) mass fraction, probe 2 & (l) mass fraction, probe 3\\
\end{tabular}
\caption{This figure compares the full-model states at three probe locations with the predictions obtained with the AADEIM and static reduced model of dimension $\nr = 6$. The AADEIM model provides accurate predictions for all quantities at all probe locations, whereas the static reduced model provides poor approximations. Note that the species mass fraction at probe location 2 and 3 is zero.} 
\label{fig:ProbeN6High}
\end{figure}

We now compare the AADEIM and static reduced model based on probes of the states at $x = 0.0025  $, $x = 0.005$, and $x = 0.0075$. The probes for dimension $\nr = 6$, update frequency $\nz = 3$ of the samples, and number of sampling points $\nrs = 1{,}024$ is shown in Figure~\ref{fig:ProbeN6High}. The probes obtained with the static reduced model of dimension $\nr = 6$ and the full model are plotted too. The AADEIM model provides an accurate predictions of the full model over all times and probe locations for all quantities, whereas the static reduced model fails to provide meaningful approximations. Note that the mass fraction at probe location 2 and 3 is zero.

\section{Conclusions}\label{sec:Conc}
The considered benchmark problem of a model premixed flame with artificial pressure forcing relies on strong simplifications of physics that are present in more realistic scenarios of chemically reacting flows. However, it preserves the transport-dominated and multiscale nature of the dynamics, which are major challenges for model reduction with linear approximations. We showed numerically that online adaptive model reduction with the AADEIM method provides accurate predictions of the flame dynamics with few degrees of freedom. The AADEIM method leverages two properties of the considered problem. First, the states of the considered problem have a local low-rank structure in the sense that the singular values decay quickly for snapshots in a local time window. Second, the residual of the AADEIM approximation is local in the spatial domain, which means that few sampling points are sufficient to inform the adaptation of the reduced basis. Reduced models based on AADEIM build on these two properties to derive nonlinear approximations of latent dynamics and so enable predictions of transport-dominated dynamics far outside of training regimes.

\bibliographystyle{splncs04}
\bibliography{main}

%% file: figures/FOM_snapshot_t7.5e-06_pressure.tex
\begingroup
  \makeatletter
  \providecommand\color[2][]{%
    \GenericError{(gnuplot) \space\space\space\@spaces}{%
      Package color not loaded in conjunction with
      terminal option `colourtext'%
    }{See the gnuplot documentation for explanation.%
    }{Either use 'blacktext' in gnuplot or load the package
      color.sty in LaTeX.}%
    \renewcommand\color[2][]{}%
  }%
  \providecommand\includegraphics[2][]{%
    \GenericError{(gnuplot) \space\space\space\@spaces}{%
      Package graphicx or graphics not loaded%
    }{See the gnuplot documentation for explanation.%
    }{The gnuplot epslatex terminal needs graphicx.sty or graphics.sty.}%
    \renewcommand\includegraphics[2][]{}%
  }%
  \providecommand\rotatebox[2]{#2}%
  \@ifundefined{ifGPcolor}{%
    \newif\ifGPcolor
    \GPcolortrue
  }{}%
  \@ifundefined{ifGPblacktext}{%
    \newif\ifGPblacktext
    \GPblacktexttrue
  }{}%
  \let\gplgaddtomacro\g@addto@macro
  \gdef\gplbacktext{}%
  \gdef\gplfronttext{}%
  \makeatother
  \ifGPblacktext
    \def\colorrgb#1{}%
    \def\colorgray#1{}%
  \else
    \ifGPcolor
      \def\colorrgb#1{\color[rgb]{#1}}%
      \def\colorgray#1{\color[gray]{#1}}%
      \expandafter\def\csname LTw\endcsname{\color{white}}%
      \expandafter\def\csname LTb\endcsname{\color{black}}%
      \expandafter\def\csname LTa\endcsname{\color{black}}%
      \expandafter\def\csname LT0\endcsname{\color[rgb]{1,0,0}}%
      \expandafter\def\csname LT1\endcsname{\color[rgb]{0,1,0}}%
      \expandafter\def\csname LT2\endcsname{\color[rgb]{0,0,1}}%
      \expandafter\def\csname LT3\endcsname{\color[rgb]{1,0,1}}%
      \expandafter\def\csname LT4\endcsname{\color[rgb]{0,1,1}}%
      \expandafter\def\csname LT5\endcsname{\color[rgb]{1,1,0}}%
      \expandafter\def\csname LT6\endcsname{\color[rgb]{0,0,0}}%
      \expandafter\def\csname LT7\endcsname{\color[rgb]{1,0.3,0}}%
      \expandafter\def\csname LT8\endcsname{\color[rgb]{0.5,0.5,0.5}}%
    \else
      \def\colorrgb#1{\color{black}}%
      \def\colorgray#1{\color[gray]{#1}}%
      \expandafter\def\csname LTw\endcsname{\color{white}}%
      \expandafter\def\csname LTb\endcsname{\color{black}}%
      \expandafter\def\csname LTa\endcsname{\color{black}}%
      \expandafter\def\csname LT0\endcsname{\color{black}}%
      \expandafter\def\csname LT1\endcsname{\color{black}}%
      \expandafter\def\csname LT2\endcsname{\color{black}}%
      \expandafter\def\csname LT3\endcsname{\color{black}}%
      \expandafter\def\csname LT4\endcsname{\color{black}}%
      \expandafter\def\csname LT5\endcsname{\color{black}}%
      \expandafter\def\csname LT6\endcsname{\color{black}}%
      \expandafter\def\csname LT7\endcsname{\color{black}}%
      \expandafter\def\csname LT8\endcsname{\color{black}}%
    \fi
  \fi
    \setlength{\unitlength}{0.0500bp}%
    \ifx\gptboxheight\undefined%
      \newlength{\gptboxheight}%
      \newlength{\gptboxwidth}%
      \newsavebox{\gptboxtext}%
    \fi%
    \setlength{\fboxrule}{0.5pt}%
    \setlength{\fboxsep}{1pt}%
\begin{picture}(7200.00,5040.00)%
    \gplgaddtomacro\gplbacktext{%
      \csname LTb\endcsname
      \put(1876,1134){\makebox(0,0)[r]{\strut{}9.20e+05}}%
      \csname LTb\endcsname
      \put(1876,1663){\makebox(0,0)[r]{\strut{}9.40e+05}}%
      \csname LTb\endcsname
      \put(1876,2193){\makebox(0,0)[r]{\strut{}9.60e+05}}%
      \csname LTb\endcsname
      \put(1876,2722){\makebox(0,0)[r]{\strut{}9.80e+05}}%
      \csname LTb\endcsname
      \put(1876,3251){\makebox(0,0)[r]{\strut{}1.00e+06}}%
      \csname LTb\endcsname
      \put(1876,3781){\makebox(0,0)[r]{\strut{}1.02e+06}}%
      \csname LTb\endcsname
      \put(1876,4310){\makebox(0,0)[r]{\strut{}1.04e+06}}%
      \csname LTb\endcsname
      \put(2044,616){\makebox(0,0){\strut{}0e+00}}%
      \csname LTb\endcsname
      \put(4370,616){\makebox(0,0){\strut{}5e-03}}%
      \csname LTb\endcsname
      \put(6695,616){\makebox(0,0){\strut{}1e-02}}%
    }%
    \gplgaddtomacro\gplfronttext{%
      \csname LTb\endcsname
      \put(-238,2827){\rotatebox{-270}{\makebox(0,0){\strut{}pressure}}}%
      \put(4369,196){\makebox(0,0){\strut{}space}}%
    }%
    \gplbacktext
    \put(0,0){\includegraphics[width={360.00bp},height={252.00bp}]{{FOM_snapshot_t7.5e-06_pressure}.eps}}%
    \gplfronttext
  \end{picture}%
\endgroup

%% file: figures/FOM_snapshot_t2e-05_pressure.tex
\begingroup
  \makeatletter
  \providecommand\color[2][]{%
    \GenericError{(gnuplot) \space\space\space\@spaces}{%
      Package color not loaded in conjunction with
      terminal option `colourtext'%
    }{See the gnuplot documentation for explanation.%
    }{Either use 'blacktext' in gnuplot or load the package
      color.sty in LaTeX.}%
    \renewcommand\color[2][]{}%
  }%
  \providecommand\includegraphics[2][]{%
    \GenericError{(gnuplot) \space\space\space\@spaces}{%
      Package graphicx or graphics not loaded%
    }{See the gnuplot documentation for explanation.%
    }{The gnuplot epslatex terminal needs graphicx.sty or graphics.sty.}%
    \renewcommand\includegraphics[2][]{}%
  }%
  \providecommand\rotatebox[2]{#2}%
  \@ifundefined{ifGPcolor}{%
    \newif\ifGPcolor
    \GPcolortrue
  }{}%
  \@ifundefined{ifGPblacktext}{%
    \newif\ifGPblacktext
    \GPblacktexttrue
  }{}%
  \let\gplgaddtomacro\g@addto@macro
  \gdef\gplbacktext{}%
  \gdef\gplfronttext{}%
  \makeatother
  \ifGPblacktext
    \def\colorrgb#1{}%
    \def\colorgray#1{}%
  \else
    \ifGPcolor
      \def\colorrgb#1{\color[rgb]{#1}}%
      \def\colorgray#1{\color[gray]{#1}}%
      \expandafter\def\csname LTw\endcsname{\color{white}}%
      \expandafter\def\csname LTb\endcsname{\color{black}}%
      \expandafter\def\csname LTa\endcsname{\color{black}}%
      \expandafter\def\csname LT0\endcsname{\color[rgb]{1,0,0}}%
      \expandafter\def\csname LT1\endcsname{\color[rgb]{0,1,0}}%
      \expandafter\def\csname LT2\endcsname{\color[rgb]{0,0,1}}%
      \expandafter\def\csname LT3\endcsname{\color[rgb]{1,0,1}}%
      \expandafter\def\csname LT4\endcsname{\color[rgb]{0,1,1}}%
      \expandafter\def\csname LT5\endcsname{\color[rgb]{1,1,0}}%
      \expandafter\def\csname LT6\endcsname{\color[rgb]{0,0,0}}%
      \expandafter\def\csname LT7\endcsname{\color[rgb]{1,0.3,0}}%
      \expandafter\def\csname LT8\endcsname{\color[rgb]{0.5,0.5,0.5}}%
    \else
      \def\colorrgb#1{\color{black}}%
      \def\colorgray#1{\color[gray]{#1}}%
      \expandafter\def\csname LTw\endcsname{\color{white}}%
      \expandafter\def\csname LTb\endcsname{\color{black}}%
      \expandafter\def\csname LTa\endcsname{\color{black}}%
      \expandafter\def\csname LT0\endcsname{\color{black}}%
      \expandafter\def\csname LT1\endcsname{\color{black}}%
      \expandafter\def\csname LT2\endcsname{\color{black}}%
      \expandafter\def\csname LT3\endcsname{\color{black}}%
      \expandafter\def\csname LT4\endcsname{\color{black}}%
      \expandafter\def\csname LT5\endcsname{\color{black}}%
      \expandafter\def\csname LT6\endcsname{\color{black}}%
      \expandafter\def\csname LT7\endcsname{\color{black}}%
      \expandafter\def\csname LT8\endcsname{\color{black}}%
    \fi
  \fi
    \setlength{\unitlength}{0.0500bp}%
    \ifx\gptboxheight\undefined%
      \newlength{\gptboxheight}%
      \newlength{\gptboxwidth}%
      \newsavebox{\gptboxtext}%
    \fi%
    \setlength{\fboxrule}{0.5pt}%
    \setlength{\fboxsep}{1pt}%
\begin{picture}(7200.00,5040.00)%
    \gplgaddtomacro\gplbacktext{%
      \csname LTb\endcsname
      \put(1876,1134){\makebox(0,0)[r]{\strut{}9.20e+05}}%
      \csname LTb\endcsname
      \put(1876,1663){\makebox(0,0)[r]{\strut{}9.40e+05}}%
      \csname LTb\endcsname
      \put(1876,2193){\makebox(0,0)[r]{\strut{}9.60e+05}}%
      \csname LTb\endcsname
      \put(1876,2722){\makebox(0,0)[r]{\strut{}9.80e+05}}%
      \csname LTb\endcsname
      \put(1876,3251){\makebox(0,0)[r]{\strut{}1.00e+06}}%
      \csname LTb\endcsname
      \put(1876,3781){\makebox(0,0)[r]{\strut{}1.02e+06}}%
      \csname LTb\endcsname
      \put(1876,4310){\makebox(0,0)[r]{\strut{}1.04e+06}}%
      \csname LTb\endcsname
      \put(2044,616){\makebox(0,0){\strut{}0e+00}}%
      \csname LTb\endcsname
      \put(4370,616){\makebox(0,0){\strut{}5e-03}}%
      \csname LTb\endcsname
      \put(6695,616){\makebox(0,0){\strut{}1e-02}}%
    }%
    \gplgaddtomacro\gplfronttext{%
      \csname LTb\endcsname
      \put(-238,2827){\rotatebox{-270}{\makebox(0,0){\strut{}pressure}}}%
      \put(4369,196){\makebox(0,0){\strut{}space}}%
    }%
    \gplbacktext
    \put(0,0){\includegraphics[width={360.00bp},height={252.00bp}]{FOM_snapshot_t2e-05_pressure}}%
    \gplfronttext
  \end{picture}%
\endgroup

%% file: figures/FOM_snapshot_t3e-05_pressure.tex
\begingroup
  \makeatletter
  \providecommand\color[2][]{%
    \GenericError{(gnuplot) \space\space\space\@spaces}{%
      Package color not loaded in conjunction with
      terminal option `colourtext'%
    }{See the gnuplot documentation for explanation.%
    }{Either use 'blacktext' in gnuplot or load the package
      color.sty in LaTeX.}%
    \renewcommand\color[2][]{}%
  }%
  \providecommand\includegraphics[2][]{%
    \GenericError{(gnuplot) \space\space\space\@spaces}{%
      Package graphicx or graphics not loaded%
    }{See the gnuplot documentation for explanation.%
    }{The gnuplot epslatex terminal needs graphicx.sty or graphics.sty.}%
    \renewcommand\includegraphics[2][]{}%
  }%
  \providecommand\rotatebox[2]{#2}%
  \@ifundefined{ifGPcolor}{%
    \newif\ifGPcolor
    \GPcolortrue
  }{}%
  \@ifundefined{ifGPblacktext}{%
    \newif\ifGPblacktext
    \GPblacktexttrue
  }{}%
  \let\gplgaddtomacro\g@addto@macro
  \gdef\gplbacktext{}%
  \gdef\gplfronttext{}%
  \makeatother
  \ifGPblacktext
    \def\colorrgb#1{}%
    \def\colorgray#1{}%
  \else
    \ifGPcolor
      \def\colorrgb#1{\color[rgb]{#1}}%
      \def\colorgray#1{\color[gray]{#1}}%
      \expandafter\def\csname LTw\endcsname{\color{white}}%
      \expandafter\def\csname LTb\endcsname{\color{black}}%
      \expandafter\def\csname LTa\endcsname{\color{black}}%
      \expandafter\def\csname LT0\endcsname{\color[rgb]{1,0,0}}%
      \expandafter\def\csname LT1\endcsname{\color[rgb]{0,1,0}}%
      \expandafter\def\csname LT2\endcsname{\color[rgb]{0,0,1}}%
      \expandafter\def\csname LT3\endcsname{\color[rgb]{1,0,1}}%
      \expandafter\def\csname LT4\endcsname{\color[rgb]{0,1,1}}%
      \expandafter\def\csname LT5\endcsname{\color[rgb]{1,1,0}}%
      \expandafter\def\csname LT6\endcsname{\color[rgb]{0,0,0}}%
      \expandafter\def\csname LT7\endcsname{\color[rgb]{1,0.3,0}}%
      \expandafter\def\csname LT8\endcsname{\color[rgb]{0.5,0.5,0.5}}%
    \else
      \def\colorrgb#1{\color{black}}%
      \def\colorgray#1{\color[gray]{#1}}%
      \expandafter\def\csname LTw\endcsname{\color{white}}%
      \expandafter\def\csname LTb\endcsname{\color{black}}%
      \expandafter\def\csname LTa\endcsname{\color{black}}%
      \expandafter\def\csname LT0\endcsname{\color{black}}%
      \expandafter\def\csname LT1\endcsname{\color{black}}%
      \expandafter\def\csname LT2\endcsname{\color{black}}%
      \expandafter\def\csname LT3\endcsname{\color{black}}%
      \expandafter\def\csname LT4\endcsname{\color{black}}%
      \expandafter\def\csname LT5\endcsname{\color{black}}%
      \expandafter\def\csname LT6\endcsname{\color{black}}%
      \expandafter\def\csname LT7\endcsname{\color{black}}%
      \expandafter\def\csname LT8\endcsname{\color{black}}%
    \fi
  \fi
    \setlength{\unitlength}{0.0500bp}%
    \ifx\gptboxheight\undefined%
      \newlength{\gptboxheight}%
      \newlength{\gptboxwidth}%
      \newsavebox{\gptboxtext}%
    \fi%
    \setlength{\fboxrule}{0.5pt}%
    \setlength{\fboxsep}{1pt}%
\begin{picture}(7200.00,5040.00)%
    \gplgaddtomacro\gplbacktext{%
      \csname LTb\endcsname
      \put(1876,1134){\makebox(0,0)[r]{\strut{}9.20e+05}}%
      \csname LTb\endcsname
      \put(1876,1663){\makebox(0,0)[r]{\strut{}9.40e+05}}%
      \csname LTb\endcsname
      \put(1876,2193){\makebox(0,0)[r]{\strut{}9.60e+05}}%
      \csname LTb\endcsname
      \put(1876,2722){\makebox(0,0)[r]{\strut{}9.80e+05}}%
      \csname LTb\endcsname
      \put(1876,3251){\makebox(0,0)[r]{\strut{}1.00e+06}}%
      \csname LTb\endcsname
      \put(1876,3781){\makebox(0,0)[r]{\strut{}1.02e+06}}%
      \csname LTb\endcsname
      \put(1876,4310){\makebox(0,0)[r]{\strut{}1.04e+06}}%
      \csname LTb\endcsname
      \put(2044,616){\makebox(0,0){\strut{}0e+00}}%
      \csname LTb\endcsname
      \put(4370,616){\makebox(0,0){\strut{}5e-03}}%
      \csname LTb\endcsname
      \put(6695,616){\makebox(0,0){\strut{}1e-02}}%
    }%
    \gplgaddtomacro\gplfronttext{%
      \csname LTb\endcsname
      \put(-238,2827){\rotatebox{-270}{\makebox(0,0){\strut{}pressure}}}%
      \put(4369,196){\makebox(0,0){\strut{}space}}%
    }%
    \gplbacktext
    \put(0,0){\includegraphics[width={360.00bp},height={252.00bp}]{FOM_snapshot_t3e-05_pressure}}%
    \gplfronttext
  \end{picture}%
\endgroup

%% file: figures/FOM_snapshot_t7.5e-06_velocity.tex
\begingroup
  \makeatletter
  \providecommand\color[2][]{%
    \GenericError{(gnuplot) \space\space\space\@spaces}{%
      Package color not loaded in conjunction with
      terminal option `colourtext'%
    }{See the gnuplot documentation for explanation.%
    }{Either use 'blacktext' in gnuplot or load the package
      color.sty in LaTeX.}%
    \renewcommand\color[2][]{}%
  }%
  \providecommand\includegraphics[2][]{%
    \GenericError{(gnuplot) \space\space\space\@spaces}{%
      Package graphicx or graphics not loaded%
    }{See the gnuplot documentation for explanation.%
    }{The gnuplot epslatex terminal needs graphicx.sty or graphics.sty.}%
    \renewcommand\includegraphics[2][]{}%
  }%
  \providecommand\rotatebox[2]{#2}%
  \@ifundefined{ifGPcolor}{%
    \newif\ifGPcolor
    \GPcolortrue
  }{}%
  \@ifundefined{ifGPblacktext}{%
    \newif\ifGPblacktext
    \GPblacktexttrue
  }{}%
  \let\gplgaddtomacro\g@addto@macro
  \gdef\gplbacktext{}%
  \gdef\gplfronttext{}%
  \makeatother
  \ifGPblacktext
    \def\colorrgb#1{}%
    \def\colorgray#1{}%
  \else
    \ifGPcolor
      \def\colorrgb#1{\color[rgb]{#1}}%
      \def\colorgray#1{\color[gray]{#1}}%
      \expandafter\def\csname LTw\endcsname{\color{white}}%
      \expandafter\def\csname LTb\endcsname{\color{black}}%
      \expandafter\def\csname LTa\endcsname{\color{black}}%
      \expandafter\def\csname LT0\endcsname{\color[rgb]{1,0,0}}%
      \expandafter\def\csname LT1\endcsname{\color[rgb]{0,1,0}}%
      \expandafter\def\csname LT2\endcsname{\color[rgb]{0,0,1}}%
      \expandafter\def\csname LT3\endcsname{\color[rgb]{1,0,1}}%
      \expandafter\def\csname LT4\endcsname{\color[rgb]{0,1,1}}%
      \expandafter\def\csname LT5\endcsname{\color[rgb]{1,1,0}}%
      \expandafter\def\csname LT6\endcsname{\color[rgb]{0,0,0}}%
      \expandafter\def\csname LT7\endcsname{\color[rgb]{1,0.3,0}}%
      \expandafter\def\csname LT8\endcsname{\color[rgb]{0.5,0.5,0.5}}%
    \else
      \def\colorrgb#1{\color{black}}%
      \def\colorgray#1{\color[gray]{#1}}%
      \expandafter\def\csname LTw\endcsname{\color{white}}%
      \expandafter\def\csname LTb\endcsname{\color{black}}%
      \expandafter\def\csname LTa\endcsname{\color{black}}%
      \expandafter\def\csname LT0\endcsname{\color{black}}%
      \expandafter\def\csname LT1\endcsname{\color{black}}%
      \expandafter\def\csname LT2\endcsname{\color{black}}%
      \expandafter\def\csname LT3\endcsname{\color{black}}%
      \expandafter\def\csname LT4\endcsname{\color{black}}%
      \expandafter\def\csname LT5\endcsname{\color{black}}%
      \expandafter\def\csname LT6\endcsname{\color{black}}%
      \expandafter\def\csname LT7\endcsname{\color{black}}%
      \expandafter\def\csname LT8\endcsname{\color{black}}%
    \fi
  \fi
    \setlength{\unitlength}{0.0500bp}%
    \ifx\gptboxheight\undefined%
      \newlength{\gptboxheight}%
      \newlength{\gptboxwidth}%
      \newsavebox{\gptboxtext}%
    \fi%
    \setlength{\fboxrule}{0.5pt}%
    \setlength{\fboxsep}{1pt}%
\begin{picture}(7200.00,5040.00)%
    \gplgaddtomacro\gplbacktext{%
      \csname LTb\endcsname
      \put(1036,1090){\makebox(0,0)[r]{\strut{}$-40$}}%
      \csname LTb\endcsname
      \put(1036,1718){\makebox(0,0)[r]{\strut{}$-20$}}%
      \csname LTb\endcsname
      \put(1036,2346){\makebox(0,0)[r]{\strut{}$0$}}%
      \csname LTb\endcsname
      \put(1036,2974){\makebox(0,0)[r]{\strut{}$20$}}%
      \csname LTb\endcsname
      \put(1036,3602){\makebox(0,0)[r]{\strut{}$40$}}%
      \csname LTb\endcsname
      \put(1036,4230){\makebox(0,0)[r]{\strut{}$60$}}%
      \csname LTb\endcsname
      \put(1204,616){\makebox(0,0){\strut{}0e+00}}%
      \csname LTb\endcsname
      \put(3950,616){\makebox(0,0){\strut{}5e-03}}%
      \csname LTb\endcsname
      \put(6695,616){\makebox(0,0){\strut{}1e-02}}%
    }%
    \gplgaddtomacro\gplfronttext{%
      \csname LTb\endcsname
      \put(98,2827){\rotatebox{-270}{\makebox(0,0){\strut{}velocity}}}%
      \put(3949,196){\makebox(0,0){\strut{}space}}%
    }%
    \gplbacktext
    \put(0,0){\includegraphics[width={360.00bp},height={252.00bp}]{{FOM_snapshot_t7.5e-06_velocity}.eps}}%
    \gplfronttext
  \end{picture}%
\endgroup

%% file: figures/FOM_snapshot_t2e-05_velocity.tex
\begingroup
  \makeatletter
  \providecommand\color[2][]{%
    \GenericError{(gnuplot) \space\space\space\@spaces}{%
      Package color not loaded in conjunction with
      terminal option `colourtext'%
    }{See the gnuplot documentation for explanation.%
    }{Either use 'blacktext' in gnuplot or load the package
      color.sty in LaTeX.}%
    \renewcommand\color[2][]{}%
  }%
  \providecommand\includegraphics[2][]{%
    \GenericError{(gnuplot) \space\space\space\@spaces}{%
      Package graphicx or graphics not loaded%
    }{See the gnuplot documentation for explanation.%
    }{The gnuplot epslatex terminal needs graphicx.sty or graphics.sty.}%
    \renewcommand\includegraphics[2][]{}%
  }%
  \providecommand\rotatebox[2]{#2}%
  \@ifundefined{ifGPcolor}{%
    \newif\ifGPcolor
    \GPcolortrue
  }{}%
  \@ifundefined{ifGPblacktext}{%
    \newif\ifGPblacktext
    \GPblacktexttrue
  }{}%
  \let\gplgaddtomacro\g@addto@macro
  \gdef\gplbacktext{}%
  \gdef\gplfronttext{}%
  \makeatother
  \ifGPblacktext
    \def\colorrgb#1{}%
    \def\colorgray#1{}%
  \else
    \ifGPcolor
      \def\colorrgb#1{\color[rgb]{#1}}%
      \def\colorgray#1{\color[gray]{#1}}%
      \expandafter\def\csname LTw\endcsname{\color{white}}%
      \expandafter\def\csname LTb\endcsname{\color{black}}%
      \expandafter\def\csname LTa\endcsname{\color{black}}%
      \expandafter\def\csname LT0\endcsname{\color[rgb]{1,0,0}}%
      \expandafter\def\csname LT1\endcsname{\color[rgb]{0,1,0}}%
      \expandafter\def\csname LT2\endcsname{\color[rgb]{0,0,1}}%
      \expandafter\def\csname LT3\endcsname{\color[rgb]{1,0,1}}%
      \expandafter\def\csname LT4\endcsname{\color[rgb]{0,1,1}}%
      \expandafter\def\csname LT5\endcsname{\color[rgb]{1,1,0}}%
      \expandafter\def\csname LT6\endcsname{\color[rgb]{0,0,0}}%
      \expandafter\def\csname LT7\endcsname{\color[rgb]{1,0.3,0}}%
      \expandafter\def\csname LT8\endcsname{\color[rgb]{0.5,0.5,0.5}}%
    \else
      \def\colorrgb#1{\color{black}}%
      \def\colorgray#1{\color[gray]{#1}}%
      \expandafter\def\csname LTw\endcsname{\color{white}}%
      \expandafter\def\csname LTb\endcsname{\color{black}}%
      \expandafter\def\csname LTa\endcsname{\color{black}}%
      \expandafter\def\csname LT0\endcsname{\color{black}}%
      \expandafter\def\csname LT1\endcsname{\color{black}}%
      \expandafter\def\csname LT2\endcsname{\color{black}}%
      \expandafter\def\csname LT3\endcsname{\color{black}}%
      \expandafter\def\csname LT4\endcsname{\color{black}}%
      \expandafter\def\csname LT5\endcsname{\color{black}}%
      \expandafter\def\csname LT6\endcsname{\color{black}}%
      \expandafter\def\csname LT7\endcsname{\color{black}}%
      \expandafter\def\csname LT8\endcsname{\color{black}}%
    \fi
  \fi
    \setlength{\unitlength}{0.0500bp}%
    \ifx\gptboxheight\undefined%
      \newlength{\gptboxheight}%
      \newlength{\gptboxwidth}%
      \newsavebox{\gptboxtext}%
    \fi%
    \setlength{\fboxrule}{0.5pt}%
    \setlength{\fboxsep}{1pt}%
\begin{picture}(7200.00,5040.00)%
    \gplgaddtomacro\gplbacktext{%
      \csname LTb\endcsname
      \put(1036,1090){\makebox(0,0)[r]{\strut{}$-40$}}%
      \csname LTb\endcsname
      \put(1036,1718){\makebox(0,0)[r]{\strut{}$-20$}}%
      \csname LTb\endcsname
      \put(1036,2346){\makebox(0,0)[r]{\strut{}$0$}}%
      \csname LTb\endcsname
      \put(1036,2974){\makebox(0,0)[r]{\strut{}$20$}}%
      \csname LTb\endcsname
      \put(1036,3602){\makebox(0,0)[r]{\strut{}$40$}}%
      \csname LTb\endcsname
      \put(1036,4230){\makebox(0,0)[r]{\strut{}$60$}}%
      \csname LTb\endcsname
      \put(1204,616){\makebox(0,0){\strut{}0e+00}}%
      \csname LTb\endcsname
      \put(3950,616){\makebox(0,0){\strut{}5e-03}}%
      \csname LTb\endcsname
      \put(6695,616){\makebox(0,0){\strut{}1e-02}}%
    }%
    \gplgaddtomacro\gplfronttext{%
      \csname LTb\endcsname
      \put(98,2827){\rotatebox{-270}{\makebox(0,0){\strut{}velocity}}}%
      \put(3949,196){\makebox(0,0){\strut{}space}}%
    }%
    \gplbacktext
    \put(0,0){\includegraphics[width={360.00bp},height={252.00bp}]{FOM_snapshot_t2e-05_velocity}}%
    \gplfronttext
  \end{picture}%
\endgroup

%% file: figures/FOM_snapshot_t3e-05_velocity.tex
\begingroup
  \makeatletter
  \providecommand\color[2][]{%
    \GenericError{(gnuplot) \space\space\space\@spaces}{%
      Package color not loaded in conjunction with
      terminal option `colourtext'%
    }{See the gnuplot documentation for explanation.%
    }{Either use 'blacktext' in gnuplot or load the package
      color.sty in LaTeX.}%
    \renewcommand\color[2][]{}%
  }%
  \providecommand\includegraphics[2][]{%
    \GenericError{(gnuplot) \space\space\space\@spaces}{%
      Package graphicx or graphics not loaded%
    }{See the gnuplot documentation for explanation.%
    }{The gnuplot epslatex terminal needs graphicx.sty or graphics.sty.}%
    \renewcommand\includegraphics[2][]{}%
  }%
  \providecommand\rotatebox[2]{#2}%
  \@ifundefined{ifGPcolor}{%
    \newif\ifGPcolor
    \GPcolortrue
  }{}%
  \@ifundefined{ifGPblacktext}{%
    \newif\ifGPblacktext
    \GPblacktexttrue
  }{}%
  \let\gplgaddtomacro\g@addto@macro
  \gdef\gplbacktext{}%
  \gdef\gplfronttext{}%
  \makeatother
  \ifGPblacktext
    \def\colorrgb#1{}%
    \def\colorgray#1{}%
  \else
    \ifGPcolor
      \def\colorrgb#1{\color[rgb]{#1}}%
      \def\colorgray#1{\color[gray]{#1}}%
      \expandafter\def\csname LTw\endcsname{\color{white}}%
      \expandafter\def\csname LTb\endcsname{\color{black}}%
      \expandafter\def\csname LTa\endcsname{\color{black}}%
      \expandafter\def\csname LT0\endcsname{\color[rgb]{1,0,0}}%
      \expandafter\def\csname LT1\endcsname{\color[rgb]{0,1,0}}%
      \expandafter\def\csname LT2\endcsname{\color[rgb]{0,0,1}}%
      \expandafter\def\csname LT3\endcsname{\color[rgb]{1,0,1}}%
      \expandafter\def\csname LT4\endcsname{\color[rgb]{0,1,1}}%
      \expandafter\def\csname LT5\endcsname{\color[rgb]{1,1,0}}%
      \expandafter\def\csname LT6\endcsname{\color[rgb]{0,0,0}}%
      \expandafter\def\csname LT7\endcsname{\color[rgb]{1,0.3,0}}%
      \expandafter\def\csname LT8\endcsname{\color[rgb]{0.5,0.5,0.5}}%
    \else
      \def\colorrgb#1{\color{black}}%
      \def\colorgray#1{\color[gray]{#1}}%
      \expandafter\def\csname LTw\endcsname{\color{white}}%
      \expandafter\def\csname LTb\endcsname{\color{black}}%
      \expandafter\def\csname LTa\endcsname{\color{black}}%
      \expandafter\def\csname LT0\endcsname{\color{black}}%
      \expandafter\def\csname LT1\endcsname{\color{black}}%
      \expandafter\def\csname LT2\endcsname{\color{black}}%
      \expandafter\def\csname LT3\endcsname{\color{black}}%
      \expandafter\def\csname LT4\endcsname{\color{black}}%
      \expandafter\def\csname LT5\endcsname{\color{black}}%
      \expandafter\def\csname LT6\endcsname{\color{black}}%
      \expandafter\def\csname LT7\endcsname{\color{black}}%
      \expandafter\def\csname LT8\endcsname{\color{black}}%
    \fi
  \fi
    \setlength{\unitlength}{0.0500bp}%
    \ifx\gptboxheight\undefined%
      \newlength{\gptboxheight}%
      \newlength{\gptboxwidth}%
      \newsavebox{\gptboxtext}%
    \fi%
    \setlength{\fboxrule}{0.5pt}%
    \setlength{\fboxsep}{1pt}%
\begin{picture}(7200.00,5040.00)%
    \gplgaddtomacro\gplbacktext{%
      \csname LTb\endcsname
      \put(1036,1090){\makebox(0,0)[r]{\strut{}$-40$}}%
      \csname LTb\endcsname
      \put(1036,1718){\makebox(0,0)[r]{\strut{}$-20$}}%
      \csname LTb\endcsname
      \put(1036,2346){\makebox(0,0)[r]{\strut{}$0$}}%
      \csname LTb\endcsname
      \put(1036,2974){\makebox(0,0)[r]{\strut{}$20$}}%
      \csname LTb\endcsname
      \put(1036,3602){\makebox(0,0)[r]{\strut{}$40$}}%
      \csname LTb\endcsname
      \put(1036,4230){\makebox(0,0)[r]{\strut{}$60$}}%
      \csname LTb\endcsname
      \put(1204,616){\makebox(0,0){\strut{}0e+00}}%
      \csname LTb\endcsname
      \put(3950,616){\makebox(0,0){\strut{}5e-03}}%
      \csname LTb\endcsname
      \put(6695,616){\makebox(0,0){\strut{}1e-02}}%
    }%
    \gplgaddtomacro\gplfronttext{%
      \csname LTb\endcsname
      \put(98,2827){\rotatebox{-270}{\makebox(0,0){\strut{}velocity}}}%
      \put(3949,196){\makebox(0,0){\strut{}space}}%
    }%
    \gplbacktext
    \put(0,0){\includegraphics[width={360.00bp},height={252.00bp}]{FOM_snapshot_t3e-05_velocity}}%
    \gplfronttext
  \end{picture}%
\endgroup

%% file: figures/FOM_snapshot_t7.5e-06_temperature.tex
\begingroup
  \makeatletter
  \providecommand\color[2][]{%
    \GenericError{(gnuplot) \space\space\space\@spaces}{%
      Package color not loaded in conjunction with
      terminal option `colourtext'%
    }{See the gnuplot documentation for explanation.%
    }{Either use 'blacktext' in gnuplot or load the package
      color.sty in LaTeX.}%
    \renewcommand\color[2][]{}%
  }%
  \providecommand\includegraphics[2][]{%
    \GenericError{(gnuplot) \space\space\space\@spaces}{%
      Package graphicx or graphics not loaded%
    }{See the gnuplot documentation for explanation.%
    }{The gnuplot epslatex terminal needs graphicx.sty or graphics.sty.}%
    \renewcommand\includegraphics[2][]{}%
  }%
  \providecommand\rotatebox[2]{#2}%
  \@ifundefined{ifGPcolor}{%
    \newif\ifGPcolor
    \GPcolortrue
  }{}%
  \@ifundefined{ifGPblacktext}{%
    \newif\ifGPblacktext
    \GPblacktexttrue
  }{}%
  \let\gplgaddtomacro\g@addto@macro
  \gdef\gplbacktext{}%
  \gdef\gplfronttext{}%
  \makeatother
  \ifGPblacktext
    \def\colorrgb#1{}%
    \def\colorgray#1{}%
  \else
    \ifGPcolor
      \def\colorrgb#1{\color[rgb]{#1}}%
      \def\colorgray#1{\color[gray]{#1}}%
      \expandafter\def\csname LTw\endcsname{\color{white}}%
      \expandafter\def\csname LTb\endcsname{\color{black}}%
      \expandafter\def\csname LTa\endcsname{\color{black}}%
      \expandafter\def\csname LT0\endcsname{\color[rgb]{1,0,0}}%
      \expandafter\def\csname LT1\endcsname{\color[rgb]{0,1,0}}%
      \expandafter\def\csname LT2\endcsname{\color[rgb]{0,0,1}}%
      \expandafter\def\csname LT3\endcsname{\color[rgb]{1,0,1}}%
      \expandafter\def\csname LT4\endcsname{\color[rgb]{0,1,1}}%
      \expandafter\def\csname LT5\endcsname{\color[rgb]{1,1,0}}%
      \expandafter\def\csname LT6\endcsname{\color[rgb]{0,0,0}}%
      \expandafter\def\csname LT7\endcsname{\color[rgb]{1,0.3,0}}%
      \expandafter\def\csname LT8\endcsname{\color[rgb]{0.5,0.5,0.5}}%
    \else
      \def\colorrgb#1{\color{black}}%
      \def\colorgray#1{\color[gray]{#1}}%
      \expandafter\def\csname LTw\endcsname{\color{white}}%
      \expandafter\def\csname LTb\endcsname{\color{black}}%
      \expandafter\def\csname LTa\endcsname{\color{black}}%
      \expandafter\def\csname LT0\endcsname{\color{black}}%
      \expandafter\def\csname LT1\endcsname{\color{black}}%
      \expandafter\def\csname LT2\endcsname{\color{black}}%
      \expandafter\def\csname LT3\endcsname{\color{black}}%
      \expandafter\def\csname LT4\endcsname{\color{black}}%
      \expandafter\def\csname LT5\endcsname{\color{black}}%
      \expandafter\def\csname LT6\endcsname{\color{black}}%
      \expandafter\def\csname LT7\endcsname{\color{black}}%
      \expandafter\def\csname LT8\endcsname{\color{black}}%
    \fi
  \fi
    \setlength{\unitlength}{0.0500bp}%
    \ifx\gptboxheight\undefined%
      \newlength{\gptboxheight}%
      \newlength{\gptboxwidth}%
      \newsavebox{\gptboxtext}%
    \fi%
    \setlength{\fboxrule}{0.5pt}%
    \setlength{\fboxsep}{1pt}%
\begin{picture}(7200.00,5040.00)%
    \gplgaddtomacro\gplbacktext{%
      \csname LTb\endcsname
      \put(1204,1275){\makebox(0,0)[r]{\strut{}$500$}}%
      \csname LTb\endcsname
      \put(1204,2081){\makebox(0,0)[r]{\strut{}$1000$}}%
      \csname LTb\endcsname
      \put(1204,2887){\makebox(0,0)[r]{\strut{}$1500$}}%
      \csname LTb\endcsname
      \put(1204,3693){\makebox(0,0)[r]{\strut{}$2000$}}%
      \csname LTb\endcsname
      \put(1204,4499){\makebox(0,0)[r]{\strut{}$2500$}}%
      \csname LTb\endcsname
      \put(1372,616){\makebox(0,0){\strut{}0e+00}}%
      \csname LTb\endcsname
      \put(4034,616){\makebox(0,0){\strut{}5e-03}}%
      \csname LTb\endcsname
      \put(6695,616){\makebox(0,0){\strut{}1e-02}}%
    }%
    \gplgaddtomacro\gplfronttext{%
      \csname LTb\endcsname
      \put(-70,2827){\rotatebox{-270}{\makebox(0,0){\strut{}temperature}}}%
      \put(4033,196){\makebox(0,0){\strut{}space}}%
    }%
    \gplbacktext
    \put(0,0){\includegraphics[width={360.00bp},height={252.00bp}]{{FOM_snapshot_t7.5e-06_temperature}.eps}}%
    \gplfronttext
  \end{picture}%
\endgroup

%% file: figures/FOM_snapshot_t2e-05_temperature.tex
\begingroup
  \makeatletter
  \providecommand\color[2][]{%
    \GenericError{(gnuplot) \space\space\space\@spaces}{%
      Package color not loaded in conjunction with
      terminal option `colourtext'%
    }{See the gnuplot documentation for explanation.%
    }{Either use 'blacktext' in gnuplot or load the package
      color.sty in LaTeX.}%
    \renewcommand\color[2][]{}%
  }%
  \providecommand\includegraphics[2][]{%
    \GenericError{(gnuplot) \space\space\space\@spaces}{%
      Package graphicx or graphics not loaded%
    }{See the gnuplot documentation for explanation.%
    }{The gnuplot epslatex terminal needs graphicx.sty or graphics.sty.}%
    \renewcommand\includegraphics[2][]{}%
  }%
  \providecommand\rotatebox[2]{#2}%
  \@ifundefined{ifGPcolor}{%
    \newif\ifGPcolor
    \GPcolortrue
  }{}%
  \@ifundefined{ifGPblacktext}{%
    \newif\ifGPblacktext
    \GPblacktexttrue
  }{}%
  \let\gplgaddtomacro\g@addto@macro
  \gdef\gplbacktext{}%
  \gdef\gplfronttext{}%
  \makeatother
  \ifGPblacktext
    \def\colorrgb#1{}%
    \def\colorgray#1{}%
  \else
    \ifGPcolor
      \def\colorrgb#1{\color[rgb]{#1}}%
      \def\colorgray#1{\color[gray]{#1}}%
      \expandafter\def\csname LTw\endcsname{\color{white}}%
      \expandafter\def\csname LTb\endcsname{\color{black}}%
      \expandafter\def\csname LTa\endcsname{\color{black}}%
      \expandafter\def\csname LT0\endcsname{\color[rgb]{1,0,0}}%
      \expandafter\def\csname LT1\endcsname{\color[rgb]{0,1,0}}%
      \expandafter\def\csname LT2\endcsname{\color[rgb]{0,0,1}}%
      \expandafter\def\csname LT3\endcsname{\color[rgb]{1,0,1}}%
      \expandafter\def\csname LT4\endcsname{\color[rgb]{0,1,1}}%
      \expandafter\def\csname LT5\endcsname{\color[rgb]{1,1,0}}%
      \expandafter\def\csname LT6\endcsname{\color[rgb]{0,0,0}}%
      \expandafter\def\csname LT7\endcsname{\color[rgb]{1,0.3,0}}%
      \expandafter\def\csname LT8\endcsname{\color[rgb]{0.5,0.5,0.5}}%
    \else
      \def\colorrgb#1{\color{black}}%
      \def\colorgray#1{\color[gray]{#1}}%
      \expandafter\def\csname LTw\endcsname{\color{white}}%
      \expandafter\def\csname LTb\endcsname{\color{black}}%
      \expandafter\def\csname LTa\endcsname{\color{black}}%
      \expandafter\def\csname LT0\endcsname{\color{black}}%
      \expandafter\def\csname LT1\endcsname{\color{black}}%
      \expandafter\def\csname LT2\endcsname{\color{black}}%
      \expandafter\def\csname LT3\endcsname{\color{black}}%
      \expandafter\def\csname LT4\endcsname{\color{black}}%
      \expandafter\def\csname LT5\endcsname{\color{black}}%
      \expandafter\def\csname LT6\endcsname{\color{black}}%
      \expandafter\def\csname LT7\endcsname{\color{black}}%
      \expandafter\def\csname LT8\endcsname{\color{black}}%
    \fi
  \fi
    \setlength{\unitlength}{0.0500bp}%
    \ifx\gptboxheight\undefined%
      \newlength{\gptboxheight}%
      \newlength{\gptboxwidth}%
      \newsavebox{\gptboxtext}%
    \fi%
    \setlength{\fboxrule}{0.5pt}%
    \setlength{\fboxsep}{1pt}%
\begin{picture}(7200.00,5040.00)%
    \gplgaddtomacro\gplbacktext{%
      \csname LTb\endcsname
      \put(1204,1275){\makebox(0,0)[r]{\strut{}$500$}}%
      \csname LTb\endcsname
      \put(1204,2081){\makebox(0,0)[r]{\strut{}$1000$}}%
      \csname LTb\endcsname
      \put(1204,2887){\makebox(0,0)[r]{\strut{}$1500$}}%
      \csname LTb\endcsname
      \put(1204,3693){\makebox(0,0)[r]{\strut{}$2000$}}%
      \csname LTb\endcsname
      \put(1204,4499){\makebox(0,0)[r]{\strut{}$2500$}}%
      \csname LTb\endcsname
      \put(1372,616){\makebox(0,0){\strut{}0e+00}}%
      \csname LTb\endcsname
      \put(4034,616){\makebox(0,0){\strut{}5e-03}}%
      \csname LTb\endcsname
      \put(6695,616){\makebox(0,0){\strut{}1e-02}}%
    }%
    \gplgaddtomacro\gplfronttext{%
      \csname LTb\endcsname
      \put(-70,2827){\rotatebox{-270}{\makebox(0,0){\strut{}temperature}}}%
      \put(4033,196){\makebox(0,0){\strut{}space}}%
    }%
    \gplbacktext
    \put(0,0){\includegraphics[width={360.00bp},height={252.00bp}]{FOM_snapshot_t2e-05_temperature}}%
    \gplfronttext
  \end{picture}%
\endgroup

%% file: figures/FOM_snapshot_t3e-05_temperature.tex
\begingroup
  \makeatletter
  \providecommand\color[2][]{%
    \GenericError{(gnuplot) \space\space\space\@spaces}{%
      Package color not loaded in conjunction with
      terminal option `colourtext'%
    }{See the gnuplot documentation for explanation.%
    }{Either use 'blacktext' in gnuplot or load the package
      color.sty in LaTeX.}%
    \renewcommand\color[2][]{}%
  }%
  \providecommand\includegraphics[2][]{%
    \GenericError{(gnuplot) \space\space\space\@spaces}{%
      Package graphicx or graphics not loaded%
    }{See the gnuplot documentation for explanation.%
    }{The gnuplot epslatex terminal needs graphicx.sty or graphics.sty.}%
    \renewcommand\includegraphics[2][]{}%
  }%
  \providecommand\rotatebox[2]{#2}%
  \@ifundefined{ifGPcolor}{%
    \newif\ifGPcolor
    \GPcolortrue
  }{}%
  \@ifundefined{ifGPblacktext}{%
    \newif\ifGPblacktext
    \GPblacktexttrue
  }{}%
  \let\gplgaddtomacro\g@addto@macro
  \gdef\gplbacktext{}%
  \gdef\gplfronttext{}%
  \makeatother
  \ifGPblacktext
    \def\colorrgb#1{}%
    \def\colorgray#1{}%
  \else
    \ifGPcolor
      \def\colorrgb#1{\color[rgb]{#1}}%
      \def\colorgray#1{\color[gray]{#1}}%
      \expandafter\def\csname LTw\endcsname{\color{white}}%
      \expandafter\def\csname LTb\endcsname{\color{black}}%
      \expandafter\def\csname LTa\endcsname{\color{black}}%
      \expandafter\def\csname LT0\endcsname{\color[rgb]{1,0,0}}%
      \expandafter\def\csname LT1\endcsname{\color[rgb]{0,1,0}}%
      \expandafter\def\csname LT2\endcsname{\color[rgb]{0,0,1}}%
      \expandafter\def\csname LT3\endcsname{\color[rgb]{1,0,1}}%
      \expandafter\def\csname LT4\endcsname{\color[rgb]{0,1,1}}%
      \expandafter\def\csname LT5\endcsname{\color[rgb]{1,1,0}}%
      \expandafter\def\csname LT6\endcsname{\color[rgb]{0,0,0}}%
      \expandafter\def\csname LT7\endcsname{\color[rgb]{1,0.3,0}}%
      \expandafter\def\csname LT8\endcsname{\color[rgb]{0.5,0.5,0.5}}%
    \else
      \def\colorrgb#1{\color{black}}%
      \def\colorgray#1{\color[gray]{#1}}%
      \expandafter\def\csname LTw\endcsname{\color{white}}%
      \expandafter\def\csname LTb\endcsname{\color{black}}%
      \expandafter\def\csname LTa\endcsname{\color{black}}%
      \expandafter\def\csname LT0\endcsname{\color{black}}%
      \expandafter\def\csname LT1\endcsname{\color{black}}%
      \expandafter\def\csname LT2\endcsname{\color{black}}%
      \expandafter\def\csname LT3\endcsname{\color{black}}%
      \expandafter\def\csname LT4\endcsname{\color{black}}%
      \expandafter\def\csname LT5\endcsname{\color{black}}%
      \expandafter\def\csname LT6\endcsname{\color{black}}%
      \expandafter\def\csname LT7\endcsname{\color{black}}%
      \expandafter\def\csname LT8\endcsname{\color{black}}%
    \fi
  \fi
    \setlength{\unitlength}{0.0500bp}%
    \ifx\gptboxheight\undefined%
      \newlength{\gptboxheight}%
      \newlength{\gptboxwidth}%
      \newsavebox{\gptboxtext}%
    \fi%
    \setlength{\fboxrule}{0.5pt}%
    \setlength{\fboxsep}{1pt}%
\begin{picture}(7200.00,5040.00)%
    \gplgaddtomacro\gplbacktext{%
      \csname LTb\endcsname
      \put(1204,1275){\makebox(0,0)[r]{\strut{}$500$}}%
      \csname LTb\endcsname
      \put(1204,2081){\makebox(0,0)[r]{\strut{}$1000$}}%
      \csname LTb\endcsname
      \put(1204,2887){\makebox(0,0)[r]{\strut{}$1500$}}%
      \csname LTb\endcsname
      \put(1204,3693){\makebox(0,0)[r]{\strut{}$2000$}}%
      \csname LTb\endcsname
      \put(1204,4499){\makebox(0,0)[r]{\strut{}$2500$}}%
      \csname LTb\endcsname
      \put(1372,616){\makebox(0,0){\strut{}0e+00}}%
      \csname LTb\endcsname
      \put(4034,616){\makebox(0,0){\strut{}5e-03}}%
      \csname LTb\endcsname
      \put(6695,616){\makebox(0,0){\strut{}1e-02}}%
    }%
    \gplgaddtomacro\gplfronttext{%
      \csname LTb\endcsname
      \put(-70,2827){\rotatebox{-270}{\makebox(0,0){\strut{}temperature}}}%
      \put(4033,196){\makebox(0,0){\strut{}space}}%
    }%
    \gplbacktext
    \put(0,0){\includegraphics[width={360.00bp},height={252.00bp}]{FOM_snapshot_t3e-05_temperature}}%
    \gplfronttext
  \end{picture}%
\endgroup

%% file: figures/FOM_snapshot_t7.5e-06_massfrac.tex
\begingroup
  \makeatletter
  \providecommand\color[2][]{%
    \GenericError{(gnuplot) \space\space\space\@spaces}{%
      Package color not loaded in conjunction with
      terminal option `colourtext'%
    }{See the gnuplot documentation for explanation.%
    }{Either use 'blacktext' in gnuplot or load the package
      color.sty in LaTeX.}%
    \renewcommand\color[2][]{}%
  }%
  \providecommand\includegraphics[2][]{%
    \GenericError{(gnuplot) \space\space\space\@spaces}{%
      Package graphicx or graphics not loaded%
    }{See the gnuplot documentation for explanation.%
    }{The gnuplot epslatex terminal needs graphicx.sty or graphics.sty.}%
    \renewcommand\includegraphics[2][]{}%
  }%
  \providecommand\rotatebox[2]{#2}%
  \@ifundefined{ifGPcolor}{%
    \newif\ifGPcolor
    \GPcolortrue
  }{}%
  \@ifundefined{ifGPblacktext}{%
    \newif\ifGPblacktext
    \GPblacktexttrue
  }{}%
  \let\gplgaddtomacro\g@addto@macro
  \gdef\gplbacktext{}%
  \gdef\gplfronttext{}%
  \makeatother
  \ifGPblacktext
    \def\colorrgb#1{}%
    \def\colorgray#1{}%
  \else
    \ifGPcolor
      \def\colorrgb#1{\color[rgb]{#1}}%
      \def\colorgray#1{\color[gray]{#1}}%
      \expandafter\def\csname LTw\endcsname{\color{white}}%
      \expandafter\def\csname LTb\endcsname{\color{black}}%
      \expandafter\def\csname LTa\endcsname{\color{black}}%
      \expandafter\def\csname LT0\endcsname{\color[rgb]{1,0,0}}%
      \expandafter\def\csname LT1\endcsname{\color[rgb]{0,1,0}}%
      \expandafter\def\csname LT2\endcsname{\color[rgb]{0,0,1}}%
      \expandafter\def\csname LT3\endcsname{\color[rgb]{1,0,1}}%
      \expandafter\def\csname LT4\endcsname{\color[rgb]{0,1,1}}%
      \expandafter\def\csname LT5\endcsname{\color[rgb]{1,1,0}}%
      \expandafter\def\csname LT6\endcsname{\color[rgb]{0,0,0}}%
      \expandafter\def\csname LT7\endcsname{\color[rgb]{1,0.3,0}}%
      \expandafter\def\csname LT8\endcsname{\color[rgb]{0.5,0.5,0.5}}%
    \else
      \def\colorrgb#1{\color{black}}%
      \def\colorgray#1{\color[gray]{#1}}%
      \expandafter\def\csname LTw\endcsname{\color{white}}%
      \expandafter\def\csname LTb\endcsname{\color{black}}%
      \expandafter\def\csname LTa\endcsname{\color{black}}%
      \expandafter\def\csname LT0\endcsname{\color{black}}%
      \expandafter\def\csname LT1\endcsname{\color{black}}%
      \expandafter\def\csname LT2\endcsname{\color{black}}%
      \expandafter\def\csname LT3\endcsname{\color{black}}%
      \expandafter\def\csname LT4\endcsname{\color{black}}%
      \expandafter\def\csname LT5\endcsname{\color{black}}%
      \expandafter\def\csname LT6\endcsname{\color{black}}%
      \expandafter\def\csname LT7\endcsname{\color{black}}%
      \expandafter\def\csname LT8\endcsname{\color{black}}%
    \fi
  \fi
    \setlength{\unitlength}{0.0500bp}%
    \ifx\gptboxheight\undefined%
      \newlength{\gptboxheight}%
      \newlength{\gptboxwidth}%
      \newsavebox{\gptboxtext}%
    \fi%
    \setlength{\fboxrule}{0.5pt}%
    \setlength{\fboxsep}{1pt}%
\begin{picture}(7200.00,5040.00)%
    \gplgaddtomacro\gplbacktext{%
      \csname LTb\endcsname
      \put(1036,896){\makebox(0,0)[r]{\strut{}$0$}}%
      \csname LTb\endcsname
      \put(1036,1632){\makebox(0,0)[r]{\strut{}$0.2$}}%
      \csname LTb\endcsname
      \put(1036,2368){\makebox(0,0)[r]{\strut{}$0.4$}}%
      \csname LTb\endcsname
      \put(1036,3103){\makebox(0,0)[r]{\strut{}$0.6$}}%
      \csname LTb\endcsname
      \put(1036,3839){\makebox(0,0)[r]{\strut{}$0.8$}}%
      \csname LTb\endcsname
      \put(1036,4575){\makebox(0,0)[r]{\strut{}$1$}}%
      \csname LTb\endcsname
      \put(1204,616){\makebox(0,0){\strut{}0e+00}}%
      \csname LTb\endcsname
      \put(3950,616){\makebox(0,0){\strut{}5e-03}}%
      \csname LTb\endcsname
      \put(6695,616){\makebox(0,0){\strut{}1e-02}}%
    }%
    \gplgaddtomacro\gplfronttext{%
      \csname LTb\endcsname
      \put(98,2827){\rotatebox{-270}{\makebox(0,0){\strut{}species mass fraction}}}%
      \put(3949,196){\makebox(0,0){\strut{}space}}%
    }%
    \gplbacktext
    \put(0,0){\includegraphics[width={360.00bp},height={252.00bp}]{{FOM_snapshot_t7.5e-06_massfrac}.eps}}%
    \gplfronttext
  \end{picture}%
\endgroup

%% file: figures/FOM_snapshot_t2e-05_massfrac.tex
\begingroup
  \makeatletter
  \providecommand\color[2][]{%
    \GenericError{(gnuplot) \space\space\space\@spaces}{%
      Package color not loaded in conjunction with
      terminal option `colourtext'%
    }{See the gnuplot documentation for explanation.%
    }{Either use 'blacktext' in gnuplot or load the package
      color.sty in LaTeX.}%
    \renewcommand\color[2][]{}%
  }%
  \providecommand\includegraphics[2][]{%
    \GenericError{(gnuplot) \space\space\space\@spaces}{%
      Package graphicx or graphics not loaded%
    }{See the gnuplot documentation for explanation.%
    }{The gnuplot epslatex terminal needs graphicx.sty or graphics.sty.}%
    \renewcommand\includegraphics[2][]{}%
  }%
  \providecommand\rotatebox[2]{#2}%
  \@ifundefined{ifGPcolor}{%
    \newif\ifGPcolor
    \GPcolortrue
  }{}%
  \@ifundefined{ifGPblacktext}{%
    \newif\ifGPblacktext
    \GPblacktexttrue
  }{}%
  \let\gplgaddtomacro\g@addto@macro
  \gdef\gplbacktext{}%
  \gdef\gplfronttext{}%
  \makeatother
  \ifGPblacktext
    \def\colorrgb#1{}%
    \def\colorgray#1{}%
  \else
    \ifGPcolor
      \def\colorrgb#1{\color[rgb]{#1}}%
      \def\colorgray#1{\color[gray]{#1}}%
      \expandafter\def\csname LTw\endcsname{\color{white}}%
      \expandafter\def\csname LTb\endcsname{\color{black}}%
      \expandafter\def\csname LTa\endcsname{\color{black}}%
      \expandafter\def\csname LT0\endcsname{\color[rgb]{1,0,0}}%
      \expandafter\def\csname LT1\endcsname{\color[rgb]{0,1,0}}%
      \expandafter\def\csname LT2\endcsname{\color[rgb]{0,0,1}}%
      \expandafter\def\csname LT3\endcsname{\color[rgb]{1,0,1}}%
      \expandafter\def\csname LT4\endcsname{\color[rgb]{0,1,1}}%
      \expandafter\def\csname LT5\endcsname{\color[rgb]{1,1,0}}%
      \expandafter\def\csname LT6\endcsname{\color[rgb]{0,0,0}}%
      \expandafter\def\csname LT7\endcsname{\color[rgb]{1,0.3,0}}%
      \expandafter\def\csname LT8\endcsname{\color[rgb]{0.5,0.5,0.5}}%
    \else
      \def\colorrgb#1{\color{black}}%
      \def\colorgray#1{\color[gray]{#1}}%
      \expandafter\def\csname LTw\endcsname{\color{white}}%
      \expandafter\def\csname LTb\endcsname{\color{black}}%
      \expandafter\def\csname LTa\endcsname{\color{black}}%
      \expandafter\def\csname LT0\endcsname{\color{black}}%
      \expandafter\def\csname LT1\endcsname{\color{black}}%
      \expandafter\def\csname LT2\endcsname{\color{black}}%
      \expandafter\def\csname LT3\endcsname{\color{black}}%
      \expandafter\def\csname LT4\endcsname{\color{black}}%
      \expandafter\def\csname LT5\endcsname{\color{black}}%
      \expandafter\def\csname LT6\endcsname{\color{black}}%
      \expandafter\def\csname LT7\endcsname{\color{black}}%
      \expandafter\def\csname LT8\endcsname{\color{black}}%
    \fi
  \fi
    \setlength{\unitlength}{0.0500bp}%
    \ifx\gptboxheight\undefined%
      \newlength{\gptboxheight}%
      \newlength{\gptboxwidth}%
      \newsavebox{\gptboxtext}%
    \fi%
    \setlength{\fboxrule}{0.5pt}%
    \setlength{\fboxsep}{1pt}%
\begin{picture}(7200.00,5040.00)%
    \gplgaddtomacro\gplbacktext{%
      \csname LTb\endcsname
      \put(1036,896){\makebox(0,0)[r]{\strut{}$0$}}%
      \csname LTb\endcsname
      \put(1036,1632){\makebox(0,0)[r]{\strut{}$0.2$}}%
      \csname LTb\endcsname
      \put(1036,2368){\makebox(0,0)[r]{\strut{}$0.4$}}%
      \csname LTb\endcsname
      \put(1036,3103){\makebox(0,0)[r]{\strut{}$0.6$}}%
      \csname LTb\endcsname
      \put(1036,3839){\makebox(0,0)[r]{\strut{}$0.8$}}%
      \csname LTb\endcsname
      \put(1036,4575){\makebox(0,0)[r]{\strut{}$1$}}%
      \csname LTb\endcsname
      \put(1204,616){\makebox(0,0){\strut{}0e+00}}%
      \csname LTb\endcsname
      \put(3950,616){\makebox(0,0){\strut{}5e-03}}%
      \csname LTb\endcsname
      \put(6695,616){\makebox(0,0){\strut{}1e-02}}%
    }%
    \gplgaddtomacro\gplfronttext{%
      \csname LTb\endcsname
      \put(98,2827){\rotatebox{-270}{\makebox(0,0){\strut{}species mass fraction}}}%
      \put(3949,196){\makebox(0,0){\strut{}space}}%
    }%
    \gplbacktext
    \put(0,0){\includegraphics[width={360.00bp},height={252.00bp}]{FOM_snapshot_t2e-05_massfrac}}%
    \gplfronttext
  \end{picture}%
\endgroup

%% file: figures/FOM_snapshot_t3e-05_massfrac.tex
\begingroup
  \makeatletter
  \providecommand\color[2][]{%
    \GenericError{(gnuplot) \space\space\space\@spaces}{%
      Package color not loaded in conjunction with
      terminal option `colourtext'%
    }{See the gnuplot documentation for explanation.%
    }{Either use 'blacktext' in gnuplot or load the package
      color.sty in LaTeX.}%
    \renewcommand\color[2][]{}%
  }%
  \providecommand\includegraphics[2][]{%
    \GenericError{(gnuplot) \space\space\space\@spaces}{%
      Package graphicx or graphics not loaded%
    }{See the gnuplot documentation for explanation.%
    }{The gnuplot epslatex terminal needs graphicx.sty or graphics.sty.}%
    \renewcommand\includegraphics[2][]{}%
  }%
  \providecommand\rotatebox[2]{#2}%
  \@ifundefined{ifGPcolor}{%
    \newif\ifGPcolor
    \GPcolortrue
  }{}%
  \@ifundefined{ifGPblacktext}{%
    \newif\ifGPblacktext
    \GPblacktexttrue
  }{}%
  \let\gplgaddtomacro\g@addto@macro
  \gdef\gplbacktext{}%
  \gdef\gplfronttext{}%
  \makeatother
  \ifGPblacktext
    \def\colorrgb#1{}%
    \def\colorgray#1{}%
  \else
    \ifGPcolor
      \def\colorrgb#1{\color[rgb]{#1}}%
      \def\colorgray#1{\color[gray]{#1}}%
      \expandafter\def\csname LTw\endcsname{\color{white}}%
      \expandafter\def\csname LTb\endcsname{\color{black}}%
      \expandafter\def\csname LTa\endcsname{\color{black}}%
      \expandafter\def\csname LT0\endcsname{\color[rgb]{1,0,0}}%
      \expandafter\def\csname LT1\endcsname{\color[rgb]{0,1,0}}%
      \expandafter\def\csname LT2\endcsname{\color[rgb]{0,0,1}}%
      \expandafter\def\csname LT3\endcsname{\color[rgb]{1,0,1}}%
      \expandafter\def\csname LT4\endcsname{\color[rgb]{0,1,1}}%
      \expandafter\def\csname LT5\endcsname{\color[rgb]{1,1,0}}%
      \expandafter\def\csname LT6\endcsname{\color[rgb]{0,0,0}}%
      \expandafter\def\csname LT7\endcsname{\color[rgb]{1,0.3,0}}%
      \expandafter\def\csname LT8\endcsname{\color[rgb]{0.5,0.5,0.5}}%
    \else
      \def\colorrgb#1{\color{black}}%
      \def\colorgray#1{\color[gray]{#1}}%
      \expandafter\def\csname LTw\endcsname{\color{white}}%
      \expandafter\def\csname LTb\endcsname{\color{black}}%
      \expandafter\def\csname LTa\endcsname{\color{black}}%
      \expandafter\def\csname LT0\endcsname{\color{black}}%
      \expandafter\def\csname LT1\endcsname{\color{black}}%
      \expandafter\def\csname LT2\endcsname{\color{black}}%
      \expandafter\def\csname LT3\endcsname{\color{black}}%
      \expandafter\def\csname LT4\endcsname{\color{black}}%
      \expandafter\def\csname LT5\endcsname{\color{black}}%
      \expandafter\def\csname LT6\endcsname{\color{black}}%
      \expandafter\def\csname LT7\endcsname{\color{black}}%
      \expandafter\def\csname LT8\endcsname{\color{black}}%
    \fi
  \fi
    \setlength{\unitlength}{0.0500bp}%
    \ifx\gptboxheight\undefined%
      \newlength{\gptboxheight}%
      \newlength{\gptboxwidth}%
      \newsavebox{\gptboxtext}%
    \fi%
    \setlength{\fboxrule}{0.5pt}%
    \setlength{\fboxsep}{1pt}%
\begin{picture}(7200.00,5040.00)%
    \gplgaddtomacro\gplbacktext{%
      \csname LTb\endcsname
      \put(1036,896){\makebox(0,0)[r]{\strut{}$0$}}%
      \csname LTb\endcsname
      \put(1036,1632){\makebox(0,0)[r]{\strut{}$0.2$}}%
      \csname LTb\endcsname
      \put(1036,2368){\makebox(0,0)[r]{\strut{}$0.4$}}%
      \csname LTb\endcsname
      \put(1036,3103){\makebox(0,0)[r]{\strut{}$0.6$}}%
      \csname LTb\endcsname
      \put(1036,3839){\makebox(0,0)[r]{\strut{}$0.8$}}%
      \csname LTb\endcsname
      \put(1036,4575){\makebox(0,0)[r]{\strut{}$1$}}%
      \csname LTb\endcsname
      \put(1204,616){\makebox(0,0){\strut{}0e+00}}%
      \csname LTb\endcsname
      \put(3950,616){\makebox(0,0){\strut{}5e-03}}%
      \csname LTb\endcsname
      \put(6695,616){\makebox(0,0){\strut{}1e-02}}%
    }%
    \gplgaddtomacro\gplfronttext{%
      \csname LTb\endcsname
      \put(98,2827){\rotatebox{-270}{\makebox(0,0){\strut{}species mass fraction}}}%
      \put(3949,196){\makebox(0,0){\strut{}space}}%
    }%
    \gplbacktext
    \put(0,0){\includegraphics[width={360.00bp},height={252.00bp}]{FOM_snapshot_t3e-05_massfrac}}%
    \gplfronttext
  \end{picture}%
\endgroup

%% file: figures/singVal.tex
\begingroup
  \makeatletter
  \providecommand\color[2][]{%
    \GenericError{(gnuplot) \space\space\space\@spaces}{%
      Package color not loaded in conjunction with
      terminal option `colourtext'%
    }{See the gnuplot documentation for explanation.%
    }{Either use 'blacktext' in gnuplot or load the package
      color.sty in LaTeX.}%
    \renewcommand\color[2][]{}%
  }%
  \providecommand\includegraphics[2][]{%
    \GenericError{(gnuplot) \space\space\space\@spaces}{%
      Package graphicx or graphics not loaded%
    }{See the gnuplot documentation for explanation.%
    }{The gnuplot epslatex terminal needs graphicx.sty or graphics.sty.}%
    \renewcommand\includegraphics[2][]{}%
  }%
  \providecommand\rotatebox[2]{#2}%
  \@ifundefined{ifGPcolor}{%
    \newif\ifGPcolor
    \GPcolortrue
  }{}%
  \@ifundefined{ifGPblacktext}{%
    \newif\ifGPblacktext
    \GPblacktexttrue
  }{}%
  \let\gplgaddtomacro\g@addto@macro
  \gdef\gplbacktext{}%
  \gdef\gplfronttext{}%
  \makeatother
  \ifGPblacktext
    \def\colorrgb#1{}%
    \def\colorgray#1{}%
  \else
    \ifGPcolor
      \def\colorrgb#1{\color[rgb]{#1}}%
      \def\colorgray#1{\color[gray]{#1}}%
      \expandafter\def\csname LTw\endcsname{\color{white}}%
      \expandafter\def\csname LTb\endcsname{\color{black}}%
      \expandafter\def\csname LTa\endcsname{\color{black}}%
      \expandafter\def\csname LT0\endcsname{\color[rgb]{1,0,0}}%
      \expandafter\def\csname LT1\endcsname{\color[rgb]{0,1,0}}%
      \expandafter\def\csname LT2\endcsname{\color[rgb]{0,0,1}}%
      \expandafter\def\csname LT3\endcsname{\color[rgb]{1,0,1}}%
      \expandafter\def\csname LT4\endcsname{\color[rgb]{0,1,1}}%
      \expandafter\def\csname LT5\endcsname{\color[rgb]{1,1,0}}%
      \expandafter\def\csname LT6\endcsname{\color[rgb]{0,0,0}}%
      \expandafter\def\csname LT7\endcsname{\color[rgb]{1,0.3,0}}%
      \expandafter\def\csname LT8\endcsname{\color[rgb]{0.5,0.5,0.5}}%
    \else
      \def\colorrgb#1{\color{black}}%
      \def\colorgray#1{\color[gray]{#1}}%
      \expandafter\def\csname LTw\endcsname{\color{white}}%
      \expandafter\def\csname LTb\endcsname{\color{black}}%
      \expandafter\def\csname LTa\endcsname{\color{black}}%
      \expandafter\def\csname LT0\endcsname{\color{black}}%
      \expandafter\def\csname LT1\endcsname{\color{black}}%
      \expandafter\def\csname LT2\endcsname{\color{black}}%
      \expandafter\def\csname LT3\endcsname{\color{black}}%
      \expandafter\def\csname LT4\endcsname{\color{black}}%
      \expandafter\def\csname LT5\endcsname{\color{black}}%
      \expandafter\def\csname LT6\endcsname{\color{black}}%
      \expandafter\def\csname LT7\endcsname{\color{black}}%
      \expandafter\def\csname LT8\endcsname{\color{black}}%
    \fi
  \fi
    \setlength{\unitlength}{0.0500bp}%
    \ifx\gptboxheight\undefined%
      \newlength{\gptboxheight}%
      \newlength{\gptboxwidth}%
      \newsavebox{\gptboxtext}%
    \fi%
    \setlength{\fboxrule}{0.5pt}%
    \setlength{\fboxsep}{1pt}%
\begin{picture}(7200.00,5040.00)%
    \gplgaddtomacro\gplbacktext{%
      \csname LTb\endcsname
      \put(1372,1247){\makebox(0,0)[r]{\strut{}1e-15}}%
      \csname LTb\endcsname
      \put(1372,2125){\makebox(0,0)[r]{\strut{}1e-10}}%
      \csname LTb\endcsname
      \put(1372,3003){\makebox(0,0)[r]{\strut{}1e-05}}%
      \csname LTb\endcsname
      \put(1372,3881){\makebox(0,0)[r]{\strut{}1e+00}}%
      \csname LTb\endcsname
      \put(1372,4759){\makebox(0,0)[r]{\strut{}1e+05}}%
      \csname LTb\endcsname
      \put(1540,616){\makebox(0,0){\strut{}$0$}}%
      \csname LTb\endcsname
      \put(2829,616){\makebox(0,0){\strut{}$50$}}%
      \csname LTb\endcsname
      \put(4118,616){\makebox(0,0){\strut{}$100$}}%
      \csname LTb\endcsname
      \put(5406,616){\makebox(0,0){\strut{}$150$}}%
      \csname LTb\endcsname
      \put(6695,616){\makebox(0,0){\strut{}$200$}}%
    }%
    \gplgaddtomacro\gplfronttext{%
      \csname LTb\endcsname
      \put(98,2827){\rotatebox{-270}{\makebox(0,0){\strut{}normalized singular value}}}%
      \put(4117,196){\makebox(0,0){\strut{}index}}%
      \csname LTb\endcsname
      \put(5120,4486){\makebox(0,0)[r]{\strut{}global}}%
    }%
    \gplbacktext
    \put(0,0){\includegraphics[width={360.00bp},height={252.00bp}]{singVal}}%
    \gplfronttext
  \end{picture}%
\endgroup

%% file: figures/singVal_local.tex
\begingroup
  \makeatletter
  \providecommand\color[2][]{%
    \GenericError{(gnuplot) \space\space\space\@spaces}{%
      Package color not loaded in conjunction with
      terminal option `colourtext'%
    }{See the gnuplot documentation for explanation.%
    }{Either use 'blacktext' in gnuplot or load the package
      color.sty in LaTeX.}%
    \renewcommand\color[2][]{}%
  }%
  \providecommand\includegraphics[2][]{%
    \GenericError{(gnuplot) \space\space\space\@spaces}{%
      Package graphicx or graphics not loaded%
    }{See the gnuplot documentation for explanation.%
    }{The gnuplot epslatex terminal needs graphicx.sty or graphics.sty.}%
    \renewcommand\includegraphics[2][]{}%
  }%
  \providecommand\rotatebox[2]{#2}%
  \@ifundefined{ifGPcolor}{%
    \newif\ifGPcolor
    \GPcolortrue
  }{}%
  \@ifundefined{ifGPblacktext}{%
    \newif\ifGPblacktext
    \GPblacktexttrue
  }{}%
  \let\gplgaddtomacro\g@addto@macro
  \gdef\gplbacktext{}%
  \gdef\gplfronttext{}%
  \makeatother
  \ifGPblacktext
    \def\colorrgb#1{}%
    \def\colorgray#1{}%
  \else
    \ifGPcolor
      \def\colorrgb#1{\color[rgb]{#1}}%
      \def\colorgray#1{\color[gray]{#1}}%
      \expandafter\def\csname LTw\endcsname{\color{white}}%
      \expandafter\def\csname LTb\endcsname{\color{black}}%
      \expandafter\def\csname LTa\endcsname{\color{black}}%
      \expandafter\def\csname LT0\endcsname{\color[rgb]{1,0,0}}%
      \expandafter\def\csname LT1\endcsname{\color[rgb]{0,1,0}}%
      \expandafter\def\csname LT2\endcsname{\color[rgb]{0,0,1}}%
      \expandafter\def\csname LT3\endcsname{\color[rgb]{1,0,1}}%
      \expandafter\def\csname LT4\endcsname{\color[rgb]{0,1,1}}%
      \expandafter\def\csname LT5\endcsname{\color[rgb]{1,1,0}}%
      \expandafter\def\csname LT6\endcsname{\color[rgb]{0,0,0}}%
      \expandafter\def\csname LT7\endcsname{\color[rgb]{1,0.3,0}}%
      \expandafter\def\csname LT8\endcsname{\color[rgb]{0.5,0.5,0.5}}%
    \else
      \def\colorrgb#1{\color{black}}%
      \def\colorgray#1{\color[gray]{#1}}%
      \expandafter\def\csname LTw\endcsname{\color{white}}%
      \expandafter\def\csname LTb\endcsname{\color{black}}%
      \expandafter\def\csname LTa\endcsname{\color{black}}%
      \expandafter\def\csname LT0\endcsname{\color{black}}%
      \expandafter\def\csname LT1\endcsname{\color{black}}%
      \expandafter\def\csname LT2\endcsname{\color{black}}%
      \expandafter\def\csname LT3\endcsname{\color{black}}%
      \expandafter\def\csname LT4\endcsname{\color{black}}%
      \expandafter\def\csname LT5\endcsname{\color{black}}%
      \expandafter\def\csname LT6\endcsname{\color{black}}%
      \expandafter\def\csname LT7\endcsname{\color{black}}%
      \expandafter\def\csname LT8\endcsname{\color{black}}%
    \fi
  \fi
    \setlength{\unitlength}{0.0500bp}%
    \ifx\gptboxheight\undefined%
      \newlength{\gptboxheight}%
      \newlength{\gptboxwidth}%
      \newsavebox{\gptboxtext}%
    \fi%
    \setlength{\fboxrule}{0.5pt}%
    \setlength{\fboxsep}{1pt}%
\begin{picture}(7200.00,5040.00)%
    \gplgaddtomacro\gplbacktext{%
      \csname LTb\endcsname
      \put(1372,1247){\makebox(0,0)[r]{\strut{}1e-15}}%
      \csname LTb\endcsname
      \put(1372,2125){\makebox(0,0)[r]{\strut{}1e-10}}%
      \csname LTb\endcsname
      \put(1372,3003){\makebox(0,0)[r]{\strut{}1e-05}}%
      \csname LTb\endcsname
      \put(1372,3881){\makebox(0,0)[r]{\strut{}1e+00}}%
      \csname LTb\endcsname
      \put(1372,4759){\makebox(0,0)[r]{\strut{}1e+05}}%
      \csname LTb\endcsname
      \put(1540,616){\makebox(0,0){\strut{}$0$}}%
      \csname LTb\endcsname
      \put(2829,616){\makebox(0,0){\strut{}$50$}}%
      \csname LTb\endcsname
      \put(4118,616){\makebox(0,0){\strut{}$100$}}%
      \csname LTb\endcsname
      \put(5406,616){\makebox(0,0){\strut{}$150$}}%
      \csname LTb\endcsname
      \put(6695,616){\makebox(0,0){\strut{}$200$}}%
    }%
    \gplgaddtomacro\gplfronttext{%
      \csname LTb\endcsname
      \put(98,2827){\rotatebox{-270}{\makebox(0,0){\strut{}normalized singular value}}}%
      \put(4117,196){\makebox(0,0){\strut{}index}}%
      \csname LTb\endcsname
      \put(5120,4486){\makebox(0,0)[r]{\strut{}global}}%
      \csname LTb\endcsname
      \put(5120,4066){\makebox(0,0)[r]{\strut{}$t \in [7.5,7.75] \times 10^{-6}$}}%
      \csname LTb\endcsname
      \put(5120,3646){\makebox(0,0)[r]{\strut{}$t \in [2,2.025] \times 10^{-5}$}}%
      \csname LTb\endcsname
      \put(5120,3226){\makebox(0,0)[r]{\strut{}$t \in [3,3.025] \times 10^{-5}$}}%
    }%
    \gplbacktext
    \put(0,0){\includegraphics[width={360.00bp},height={252.00bp}]{singVal_local}}%
    \gplfronttext
  \end{picture}%
\endgroup

%% file: figures/RelErrBar.tex
\begingroup
  \makeatletter
  \providecommand\color[2][]{%
    \GenericError{(gnuplot) \space\space\space\@spaces}{%
      Package color not loaded in conjunction with
      terminal option `colourtext'%
    }{See the gnuplot documentation for explanation.%
    }{Either use 'blacktext' in gnuplot or load the package
      color.sty in LaTeX.}%
    \renewcommand\color[2][]{}%
  }%
  \providecommand\includegraphics[2][]{%
    \GenericError{(gnuplot) \space\space\space\@spaces}{%
      Package graphicx or graphics not loaded%
    }{See the gnuplot documentation for explanation.%
    }{The gnuplot epslatex terminal needs graphicx.sty or graphics.sty.}%
    \renewcommand\includegraphics[2][]{}%
  }%
  \providecommand\rotatebox[2]{#2}%
  \@ifundefined{ifGPcolor}{%
    \newif\ifGPcolor
    \GPcolortrue
  }{}%
  \@ifundefined{ifGPblacktext}{%
    \newif\ifGPblacktext
    \GPblacktexttrue
  }{}%
  \let\gplgaddtomacro\g@addto@macro
  \gdef\gplbacktext{}%
  \gdef\gplfronttext{}%
  \makeatother
  \ifGPblacktext
    \def\colorrgb#1{}%
    \def\colorgray#1{}%
  \else
    \ifGPcolor
      \def\colorrgb#1{\color[rgb]{#1}}%
      \def\colorgray#1{\color[gray]{#1}}%
      \expandafter\def\csname LTw\endcsname{\color{white}}%
      \expandafter\def\csname LTb\endcsname{\color{black}}%
      \expandafter\def\csname LTa\endcsname{\color{black}}%
      \expandafter\def\csname LT0\endcsname{\color[rgb]{1,0,0}}%
      \expandafter\def\csname LT1\endcsname{\color[rgb]{0,1,0}}%
      \expandafter\def\csname LT2\endcsname{\color[rgb]{0,0,1}}%
      \expandafter\def\csname LT3\endcsname{\color[rgb]{1,0,1}}%
      \expandafter\def\csname LT4\endcsname{\color[rgb]{0,1,1}}%
      \expandafter\def\csname LT5\endcsname{\color[rgb]{1,1,0}}%
      \expandafter\def\csname LT6\endcsname{\color[rgb]{0,0,0}}%
      \expandafter\def\csname LT7\endcsname{\color[rgb]{1,0.3,0}}%
      \expandafter\def\csname LT8\endcsname{\color[rgb]{0.5,0.5,0.5}}%
    \else
      \def\colorrgb#1{\color{black}}%
      \def\colorgray#1{\color[gray]{#1}}%
      \expandafter\def\csname LTw\endcsname{\color{white}}%
      \expandafter\def\csname LTb\endcsname{\color{black}}%
      \expandafter\def\csname LTa\endcsname{\color{black}}%
      \expandafter\def\csname LT0\endcsname{\color{black}}%
      \expandafter\def\csname LT1\endcsname{\color{black}}%
      \expandafter\def\csname LT2\endcsname{\color{black}}%
      \expandafter\def\csname LT3\endcsname{\color{black}}%
      \expandafter\def\csname LT4\endcsname{\color{black}}%
      \expandafter\def\csname LT5\endcsname{\color{black}}%
      \expandafter\def\csname LT6\endcsname{\color{black}}%
      \expandafter\def\csname LT7\endcsname{\color{black}}%
      \expandafter\def\csname LT8\endcsname{\color{black}}%
    \fi
  \fi
    \setlength{\unitlength}{0.0500bp}%
    \ifx\gptboxheight\undefined%
      \newlength{\gptboxheight}%
      \newlength{\gptboxwidth}%
      \newsavebox{\gptboxtext}%
    \fi%
    \setlength{\fboxrule}{0.5pt}%
    \setlength{\fboxsep}{1pt}%
\begin{picture}(7200.00,5040.00)%
    \gplgaddtomacro\gplbacktext{%
      \csname LTb\endcsname
      \put(1372,560){\makebox(0,0)[r]{\strut{}1e-04}}%
      \csname LTb\endcsname
      \put(1372,2660){\makebox(0,0)[r]{\strut{}1e-03}}%
      \csname LTb\endcsname
      \put(1372,4759){\makebox(0,0)[r]{\strut{}1e-02}}%
      \put(1862,280){\makebox(0,0){\strut{}A}}%
      \put(2507,280){\makebox(0,0){\strut{}B}}%
      \put(3151,280){\makebox(0,0){\strut{}C}}%
      \put(3795,280){\makebox(0,0){\strut{}D}}%
      \put(4440,280){\makebox(0,0){\strut{}E}}%
      \put(5084,280){\makebox(0,0){\strut{}F}}%
      \put(5728,280){\makebox(0,0){\strut{}G}}%
      \put(6373,280){\makebox(0,0){\strut{}H}}%
    }%
    \gplgaddtomacro\gplfronttext{%
      \csname LTb\endcsname
      \put(-238,2659){\rotatebox{-270}{\makebox(0,0){\strut{}avg rel error}}}%
    }%
    \gplbacktext
    \put(0,0){\includegraphics[width={360.00bp},height={252.00bp}]{RelErrBar}}%
    \gplfronttext
  \end{picture}%
\endgroup

%% file: figures/NEvalsBar.tex
\begingroup
  \makeatletter
  \providecommand\color[2][]{%
    \GenericError{(gnuplot) \space\space\space\@spaces}{%
      Package color not loaded in conjunction with
      terminal option `colourtext'%
    }{See the gnuplot documentation for explanation.%
    }{Either use 'blacktext' in gnuplot or load the package
      color.sty in LaTeX.}%
    \renewcommand\color[2][]{}%
  }%
  \providecommand\includegraphics[2][]{%
    \GenericError{(gnuplot) \space\space\space\@spaces}{%
      Package graphicx or graphics not loaded%
    }{See the gnuplot documentation for explanation.%
    }{The gnuplot epslatex terminal needs graphicx.sty or graphics.sty.}%
    \renewcommand\includegraphics[2][]{}%
  }%
  \providecommand\rotatebox[2]{#2}%
  \@ifundefined{ifGPcolor}{%
    \newif\ifGPcolor
    \GPcolortrue
  }{}%
  \@ifundefined{ifGPblacktext}{%
    \newif\ifGPblacktext
    \GPblacktexttrue
  }{}%
  \let\gplgaddtomacro\g@addto@macro
  \gdef\gplbacktext{}%
  \gdef\gplfronttext{}%
  \makeatother
  \ifGPblacktext
    \def\colorrgb#1{}%
    \def\colorgray#1{}%
  \else
    \ifGPcolor
      \def\colorrgb#1{\color[rgb]{#1}}%
      \def\colorgray#1{\color[gray]{#1}}%
      \expandafter\def\csname LTw\endcsname{\color{white}}%
      \expandafter\def\csname LTb\endcsname{\color{black}}%
      \expandafter\def\csname LTa\endcsname{\color{black}}%
      \expandafter\def\csname LT0\endcsname{\color[rgb]{1,0,0}}%
      \expandafter\def\csname LT1\endcsname{\color[rgb]{0,1,0}}%
      \expandafter\def\csname LT2\endcsname{\color[rgb]{0,0,1}}%
      \expandafter\def\csname LT3\endcsname{\color[rgb]{1,0,1}}%
      \expandafter\def\csname LT4\endcsname{\color[rgb]{0,1,1}}%
      \expandafter\def\csname LT5\endcsname{\color[rgb]{1,1,0}}%
      \expandafter\def\csname LT6\endcsname{\color[rgb]{0,0,0}}%
      \expandafter\def\csname LT7\endcsname{\color[rgb]{1,0.3,0}}%
      \expandafter\def\csname LT8\endcsname{\color[rgb]{0.5,0.5,0.5}}%
    \else
      \def\colorrgb#1{\color{black}}%
      \def\colorgray#1{\color[gray]{#1}}%
      \expandafter\def\csname LTw\endcsname{\color{white}}%
      \expandafter\def\csname LTb\endcsname{\color{black}}%
      \expandafter\def\csname LTa\endcsname{\color{black}}%
      \expandafter\def\csname LT0\endcsname{\color{black}}%
      \expandafter\def\csname LT1\endcsname{\color{black}}%
      \expandafter\def\csname LT2\endcsname{\color{black}}%
      \expandafter\def\csname LT3\endcsname{\color{black}}%
      \expandafter\def\csname LT4\endcsname{\color{black}}%
      \expandafter\def\csname LT5\endcsname{\color{black}}%
      \expandafter\def\csname LT6\endcsname{\color{black}}%
      \expandafter\def\csname LT7\endcsname{\color{black}}%
      \expandafter\def\csname LT8\endcsname{\color{black}}%
    \fi
  \fi
    \setlength{\unitlength}{0.0500bp}%
    \ifx\gptboxheight\undefined%
      \newlength{\gptboxheight}%
      \newlength{\gptboxwidth}%
      \newsavebox{\gptboxtext}%
    \fi%
    \setlength{\fboxrule}{0.5pt}%
    \setlength{\fboxsep}{1pt}%
\begin{picture}(7200.00,5040.00)%
    \gplgaddtomacro\gplbacktext{%
      \csname LTb\endcsname
      \put(1708,1260){\makebox(0,0)[r]{\strut{}4.0e+07}}%
      \csname LTb\endcsname
      \put(1708,2135){\makebox(0,0)[r]{\strut{}4.2e+07}}%
      \csname LTb\endcsname
      \put(1708,3009){\makebox(0,0)[r]{\strut{}4.5e+07}}%
      \csname LTb\endcsname
      \put(1708,3884){\makebox(0,0)[r]{\strut{}4.8e+07}}%
      \csname LTb\endcsname
      \put(1708,4759){\makebox(0,0)[r]{\strut{}5.0e+07}}%
      \put(2177,280){\makebox(0,0){\strut{}A}}%
      \put(2780,280){\makebox(0,0){\strut{}B}}%
      \put(3382,280){\makebox(0,0){\strut{}C}}%
      \put(3984,280){\makebox(0,0){\strut{}D}}%
      \put(4587,280){\makebox(0,0){\strut{}E}}%
      \put(5189,280){\makebox(0,0){\strut{}F}}%
      \put(5791,280){\makebox(0,0){\strut{}G}}%
      \put(6394,280){\makebox(0,0){\strut{}H}}%
    }%
    \gplgaddtomacro\gplfronttext{%
      \csname LTb\endcsname
      \put(-238,2659){\rotatebox{-270}{\makebox(0,0){\strut{}\# evaluations}}}%
    }%
    \gplbacktext
    \put(0,0){\includegraphics[width={360.00bp},height={252.00bp}]{NEvalsBar}}%
    \gplfronttext
  \end{picture}%
\endgroup

%% file: figures/scatter_dim6.tex
\begingroup
  \makeatletter
  \providecommand\color[2][]{%
    \GenericError{(gnuplot) \space\space\space\@spaces}{%
      Package color not loaded in conjunction with
      terminal option `colourtext'%
    }{See the gnuplot documentation for explanation.%
    }{Either use 'blacktext' in gnuplot or load the package
      color.sty in LaTeX.}%
    \renewcommand\color[2][]{}%
  }%
  \providecommand\includegraphics[2][]{%
    \GenericError{(gnuplot) \space\space\space\@spaces}{%
      Package graphicx or graphics not loaded%
    }{See the gnuplot documentation for explanation.%
    }{The gnuplot epslatex terminal needs graphicx.sty or graphics.sty.}%
    \renewcommand\includegraphics[2][]{}%
  }%
  \providecommand\rotatebox[2]{#2}%
  \@ifundefined{ifGPcolor}{%
    \newif\ifGPcolor
    \GPcolortrue
  }{}%
  \@ifundefined{ifGPblacktext}{%
    \newif\ifGPblacktext
    \GPblacktexttrue
  }{}%
  \let\gplgaddtomacro\g@addto@macro
  \gdef\gplbacktext{}%
  \gdef\gplfronttext{}%
  \makeatother
  \ifGPblacktext
    \def\colorrgb#1{}%
    \def\colorgray#1{}%
  \else
    \ifGPcolor
      \def\colorrgb#1{\color[rgb]{#1}}%
      \def\colorgray#1{\color[gray]{#1}}%
      \expandafter\def\csname LTw\endcsname{\color{white}}%
      \expandafter\def\csname LTb\endcsname{\color{black}}%
      \expandafter\def\csname LTa\endcsname{\color{black}}%
      \expandafter\def\csname LT0\endcsname{\color[rgb]{1,0,0}}%
      \expandafter\def\csname LT1\endcsname{\color[rgb]{0,1,0}}%
      \expandafter\def\csname LT2\endcsname{\color[rgb]{0,0,1}}%
      \expandafter\def\csname LT3\endcsname{\color[rgb]{1,0,1}}%
      \expandafter\def\csname LT4\endcsname{\color[rgb]{0,1,1}}%
      \expandafter\def\csname LT5\endcsname{\color[rgb]{1,1,0}}%
      \expandafter\def\csname LT6\endcsname{\color[rgb]{0,0,0}}%
      \expandafter\def\csname LT7\endcsname{\color[rgb]{1,0.3,0}}%
      \expandafter\def\csname LT8\endcsname{\color[rgb]{0.5,0.5,0.5}}%
    \else
      \def\colorrgb#1{\color{black}}%
      \def\colorgray#1{\color[gray]{#1}}%
      \expandafter\def\csname LTw\endcsname{\color{white}}%
      \expandafter\def\csname LTb\endcsname{\color{black}}%
      \expandafter\def\csname LTa\endcsname{\color{black}}%
      \expandafter\def\csname LT0\endcsname{\color{black}}%
      \expandafter\def\csname LT1\endcsname{\color{black}}%
      \expandafter\def\csname LT2\endcsname{\color{black}}%
      \expandafter\def\csname LT3\endcsname{\color{black}}%
      \expandafter\def\csname LT4\endcsname{\color{black}}%
      \expandafter\def\csname LT5\endcsname{\color{black}}%
      \expandafter\def\csname LT6\endcsname{\color{black}}%
      \expandafter\def\csname LT7\endcsname{\color{black}}%
      \expandafter\def\csname LT8\endcsname{\color{black}}%
    \fi
  \fi
    \setlength{\unitlength}{0.0500bp}%
    \ifx\gptboxheight\undefined%
      \newlength{\gptboxheight}%
      \newlength{\gptboxwidth}%
      \newsavebox{\gptboxtext}%
    \fi%
    \setlength{\fboxrule}{0.5pt}%
    \setlength{\fboxsep}{1pt}%
\begin{picture}(7200.00,5040.00)%
    \gplgaddtomacro\gplbacktext{%
      \csname LTb\endcsname
      \put(1708,1540){\makebox(0,0)[r]{\strut{}4.0e+07}}%
      \csname LTb\endcsname
      \put(1708,2184){\makebox(0,0)[r]{\strut{}4.2e+07}}%
      \csname LTb\endcsname
      \put(1708,3149){\makebox(0,0)[r]{\strut{}4.5e+07}}%
      \csname LTb\endcsname
      \put(1708,4115){\makebox(0,0)[r]{\strut{}4.8e+07}}%
      \csname LTb\endcsname
      \put(1708,4759){\makebox(0,0)[r]{\strut{}5.0e+07}}%
      \csname LTb\endcsname
      \put(1876,616){\makebox(0,0){\strut{}1e-04}}%
      \csname LTb\endcsname
      \put(4285,616){\makebox(0,0){\strut{}1e-03}}%
      \csname LTb\endcsname
      \put(6695,616){\makebox(0,0){\strut{}1e-02}}%
    }%
    \gplgaddtomacro\gplfronttext{%
      \csname LTb\endcsname
      \put(-238,2827){\rotatebox{-270}{\makebox(0,0){\strut{}\# evaluations}}}%
      \put(4285,196){\makebox(0,0){\strut{}avg rel error}}%
      \csname LTb\endcsname
      \put(5792,4486){\makebox(0,0)[r]{\strut{}A}}%
      \csname LTb\endcsname
      \put(5792,4066){\makebox(0,0)[r]{\strut{}B}}%
      \csname LTb\endcsname
      \put(5792,3646){\makebox(0,0)[r]{\strut{}C}}%
      \csname LTb\endcsname
      \put(5792,3226){\makebox(0,0)[r]{\strut{}D}}%
    }%
    \gplbacktext
    \put(0,0){\includegraphics[width={360.00bp},height={252.00bp}]{scatter_dim6}}%
    \gplfronttext
  \end{picture}%
\endgroup

%% file: figures/scatter_dim7.tex
\begingroup
  \makeatletter
  \providecommand\color[2][]{%
    \GenericError{(gnuplot) \space\space\space\@spaces}{%
      Package color not loaded in conjunction with
      terminal option `colourtext'%
    }{See the gnuplot documentation for explanation.%
    }{Either use 'blacktext' in gnuplot or load the package
      color.sty in LaTeX.}%
    \renewcommand\color[2][]{}%
  }%
  \providecommand\includegraphics[2][]{%
    \GenericError{(gnuplot) \space\space\space\@spaces}{%
      Package graphicx or graphics not loaded%
    }{See the gnuplot documentation for explanation.%
    }{The gnuplot epslatex terminal needs graphicx.sty or graphics.sty.}%
    \renewcommand\includegraphics[2][]{}%
  }%
  \providecommand\rotatebox[2]{#2}%
  \@ifundefined{ifGPcolor}{%
    \newif\ifGPcolor
    \GPcolortrue
  }{}%
  \@ifundefined{ifGPblacktext}{%
    \newif\ifGPblacktext
    \GPblacktexttrue
  }{}%
  \let\gplgaddtomacro\g@addto@macro
  \gdef\gplbacktext{}%
  \gdef\gplfronttext{}%
  \makeatother
  \ifGPblacktext
    \def\colorrgb#1{}%
    \def\colorgray#1{}%
  \else
    \ifGPcolor
      \def\colorrgb#1{\color[rgb]{#1}}%
      \def\colorgray#1{\color[gray]{#1}}%
      \expandafter\def\csname LTw\endcsname{\color{white}}%
      \expandafter\def\csname LTb\endcsname{\color{black}}%
      \expandafter\def\csname LTa\endcsname{\color{black}}%
      \expandafter\def\csname LT0\endcsname{\color[rgb]{1,0,0}}%
      \expandafter\def\csname LT1\endcsname{\color[rgb]{0,1,0}}%
      \expandafter\def\csname LT2\endcsname{\color[rgb]{0,0,1}}%
      \expandafter\def\csname LT3\endcsname{\color[rgb]{1,0,1}}%
      \expandafter\def\csname LT4\endcsname{\color[rgb]{0,1,1}}%
      \expandafter\def\csname LT5\endcsname{\color[rgb]{1,1,0}}%
      \expandafter\def\csname LT6\endcsname{\color[rgb]{0,0,0}}%
      \expandafter\def\csname LT7\endcsname{\color[rgb]{1,0.3,0}}%
      \expandafter\def\csname LT8\endcsname{\color[rgb]{0.5,0.5,0.5}}%
    \else
      \def\colorrgb#1{\color{black}}%
      \def\colorgray#1{\color[gray]{#1}}%
      \expandafter\def\csname LTw\endcsname{\color{white}}%
      \expandafter\def\csname LTb\endcsname{\color{black}}%
      \expandafter\def\csname LTa\endcsname{\color{black}}%
      \expandafter\def\csname LT0\endcsname{\color{black}}%
      \expandafter\def\csname LT1\endcsname{\color{black}}%
      \expandafter\def\csname LT2\endcsname{\color{black}}%
      \expandafter\def\csname LT3\endcsname{\color{black}}%
      \expandafter\def\csname LT4\endcsname{\color{black}}%
      \expandafter\def\csname LT5\endcsname{\color{black}}%
      \expandafter\def\csname LT6\endcsname{\color{black}}%
      \expandafter\def\csname LT7\endcsname{\color{black}}%
      \expandafter\def\csname LT8\endcsname{\color{black}}%
    \fi
  \fi
    \setlength{\unitlength}{0.0500bp}%
    \ifx\gptboxheight\undefined%
      \newlength{\gptboxheight}%
      \newlength{\gptboxwidth}%
      \newsavebox{\gptboxtext}%
    \fi%
    \setlength{\fboxrule}{0.5pt}%
    \setlength{\fboxsep}{1pt}%
\begin{picture}(7200.00,5040.00)%
    \gplgaddtomacro\gplbacktext{%
      \csname LTb\endcsname
      \put(1708,1540){\makebox(0,0)[r]{\strut{}4.0e+07}}%
      \csname LTb\endcsname
      \put(1708,2184){\makebox(0,0)[r]{\strut{}4.2e+07}}%
      \csname LTb\endcsname
      \put(1708,3149){\makebox(0,0)[r]{\strut{}4.5e+07}}%
      \csname LTb\endcsname
      \put(1708,4115){\makebox(0,0)[r]{\strut{}4.8e+07}}%
      \csname LTb\endcsname
      \put(1708,4759){\makebox(0,0)[r]{\strut{}5.0e+07}}%
      \csname LTb\endcsname
      \put(1876,616){\makebox(0,0){\strut{}1e-04}}%
      \csname LTb\endcsname
      \put(4285,616){\makebox(0,0){\strut{}1e-03}}%
      \csname LTb\endcsname
      \put(6695,616){\makebox(0,0){\strut{}1e-02}}%
    }%
    \gplgaddtomacro\gplfronttext{%
      \csname LTb\endcsname
      \put(-238,2827){\rotatebox{-270}{\makebox(0,0){\strut{}\# evaluations}}}%
      \put(4285,196){\makebox(0,0){\strut{}avg rel error}}%
      \csname LTb\endcsname
      \put(5792,4486){\makebox(0,0)[r]{\strut{}E}}%
      \csname LTb\endcsname
      \put(5792,4066){\makebox(0,0)[r]{\strut{}F}}%
      \csname LTb\endcsname
      \put(5792,3646){\makebox(0,0)[r]{\strut{}G}}%
      \csname LTb\endcsname
      \put(5792,3226){\makebox(0,0)[r]{\strut{}H}}%
    }%
    \gplbacktext
    \put(0,0){\includegraphics[width={360.00bp},height={252.00bp}]{scatter_dim7}}%
    \gplfronttext
  \end{picture}%
\endgroup

%% file: figures/probe_dim6_ae2_uf3_iw15_ws7_res1024_ncycle0_AFDEIM_dt1e-09_Pressure_Probe0.tex
\begingroup
  \makeatletter
  \providecommand\color[2][]{%
    \GenericError{(gnuplot) \space\space\space\@spaces}{%
      Package color not loaded in conjunction with
      terminal option `colourtext'%
    }{See the gnuplot documentation for explanation.%
    }{Either use 'blacktext' in gnuplot or load the package
      color.sty in LaTeX.}%
    \renewcommand\color[2][]{}%
  }%
  \providecommand\includegraphics[2][]{%
    \GenericError{(gnuplot) \space\space\space\@spaces}{%
      Package graphicx or graphics not loaded%
    }{See the gnuplot documentation for explanation.%
    }{The gnuplot epslatex terminal needs graphicx.sty or graphics.sty.}%
    \renewcommand\includegraphics[2][]{}%
  }%
  \providecommand\rotatebox[2]{#2}%
  \@ifundefined{ifGPcolor}{%
    \newif\ifGPcolor
    \GPcolortrue
  }{}%
  \@ifundefined{ifGPblacktext}{%
    \newif\ifGPblacktext
    \GPblacktexttrue
  }{}%
  \let\gplgaddtomacro\g@addto@macro
  \gdef\gplbacktext{}%
  \gdef\gplfronttext{}%
  \makeatother
  \ifGPblacktext
    \def\colorrgb#1{}%
    \def\colorgray#1{}%
  \else
    \ifGPcolor
      \def\colorrgb#1{\color[rgb]{#1}}%
      \def\colorgray#1{\color[gray]{#1}}%
      \expandafter\def\csname LTw\endcsname{\color{white}}%
      \expandafter\def\csname LTb\endcsname{\color{black}}%
      \expandafter\def\csname LTa\endcsname{\color{black}}%
      \expandafter\def\csname LT0\endcsname{\color[rgb]{1,0,0}}%
      \expandafter\def\csname LT1\endcsname{\color[rgb]{0,1,0}}%
      \expandafter\def\csname LT2\endcsname{\color[rgb]{0,0,1}}%
      \expandafter\def\csname LT3\endcsname{\color[rgb]{1,0,1}}%
      \expandafter\def\csname LT4\endcsname{\color[rgb]{0,1,1}}%
      \expandafter\def\csname LT5\endcsname{\color[rgb]{1,1,0}}%
      \expandafter\def\csname LT6\endcsname{\color[rgb]{0,0,0}}%
      \expandafter\def\csname LT7\endcsname{\color[rgb]{1,0.3,0}}%
      \expandafter\def\csname LT8\endcsname{\color[rgb]{0.5,0.5,0.5}}%
    \else
      \def\colorrgb#1{\color{black}}%
      \def\colorgray#1{\color[gray]{#1}}%
      \expandafter\def\csname LTw\endcsname{\color{white}}%
      \expandafter\def\csname LTb\endcsname{\color{black}}%
      \expandafter\def\csname LTa\endcsname{\color{black}}%
      \expandafter\def\csname LT0\endcsname{\color{black}}%
      \expandafter\def\csname LT1\endcsname{\color{black}}%
      \expandafter\def\csname LT2\endcsname{\color{black}}%
      \expandafter\def\csname LT3\endcsname{\color{black}}%
      \expandafter\def\csname LT4\endcsname{\color{black}}%
      \expandafter\def\csname LT5\endcsname{\color{black}}%
      \expandafter\def\csname LT6\endcsname{\color{black}}%
      \expandafter\def\csname LT7\endcsname{\color{black}}%
      \expandafter\def\csname LT8\endcsname{\color{black}}%
    \fi
  \fi
    \setlength{\unitlength}{0.0500bp}%
    \ifx\gptboxheight\undefined%
      \newlength{\gptboxheight}%
      \newlength{\gptboxwidth}%
      \newsavebox{\gptboxtext}%
    \fi%
    \setlength{\fboxrule}{0.5pt}%
    \setlength{\fboxsep}{1pt}%
\begin{picture}(7200.00,5040.00)%
    \gplgaddtomacro\gplbacktext{%
      \csname LTb\endcsname
      \put(1876,1216){\makebox(0,0)[r]{\strut{}8.50e+05}}%
      \csname LTb\endcsname
      \put(1876,1789){\makebox(0,0)[r]{\strut{}9.00e+05}}%
      \csname LTb\endcsname
      \put(1876,2362){\makebox(0,0)[r]{\strut{}9.50e+05}}%
      \csname LTb\endcsname
      \put(1876,2934){\makebox(0,0)[r]{\strut{}1.00e+06}}%
      \csname LTb\endcsname
      \put(1876,3507){\makebox(0,0)[r]{\strut{}1.05e+06}}%
      \csname LTb\endcsname
      \put(1876,4080){\makebox(0,0)[r]{\strut{}1.10e+06}}%
      \csname LTb\endcsname
      \put(1876,4653){\makebox(0,0)[r]{\strut{}1.15e+06}}%
      \csname LTb\endcsname
      \put(2044,616){\makebox(0,0){\strut{}0e+00}}%
      \csname LTb\endcsname
      \put(3373,616){\makebox(0,0){\strut{}1e-05}}%
      \csname LTb\endcsname
      \put(6031,616){\makebox(0,0){\strut{}3e-05}}%
    }%
    \gplgaddtomacro\gplfronttext{%
      \csname LTb\endcsname
      \put(-238,2827){\rotatebox{-270}{\makebox(0,0){\strut{}pressure}}}%
      \put(4369,196){\makebox(0,0){\strut{}time}}%
      \csname LTb\endcsname
      \put(5120,4486){\makebox(0,0)[r]{\strut{}full model}}%
      \csname LTb\endcsname
      \put(5120,4066){\makebox(0,0)[r]{\strut{}AADEIM}}%
      \csname LTb\endcsname
      \put(5120,3646){\makebox(0,0)[r]{\strut{}static}}%
    }%
    \gplbacktext
    \put(0,0){\includegraphics[width={360.00bp},height={252.00bp}]{probe_dim6_ae2_uf3_iw15_ws7_res1024_ncycle0_AFDEIM_dt1e-09_Pressure_Probe0}}%
    \gplfronttext
  \end{picture}%
\endgroup

%% file: figures/probe_dim6_ae2_uf3_iw15_ws7_res1024_ncycle0_AFDEIM_dt1e-09_Pressure_Probe1.tex
\begingroup
  \makeatletter
  \providecommand\color[2][]{%
    \GenericError{(gnuplot) \space\space\space\@spaces}{%
      Package color not loaded in conjunction with
      terminal option `colourtext'%
    }{See the gnuplot documentation for explanation.%
    }{Either use 'blacktext' in gnuplot or load the package
      color.sty in LaTeX.}%
    \renewcommand\color[2][]{}%
  }%
  \providecommand\includegraphics[2][]{%
    \GenericError{(gnuplot) \space\space\space\@spaces}{%
      Package graphicx or graphics not loaded%
    }{See the gnuplot documentation for explanation.%
    }{The gnuplot epslatex terminal needs graphicx.sty or graphics.sty.}%
    \renewcommand\includegraphics[2][]{}%
  }%
  \providecommand\rotatebox[2]{#2}%
  \@ifundefined{ifGPcolor}{%
    \newif\ifGPcolor
    \GPcolortrue
  }{}%
  \@ifundefined{ifGPblacktext}{%
    \newif\ifGPblacktext
    \GPblacktexttrue
  }{}%
  \let\gplgaddtomacro\g@addto@macro
  \gdef\gplbacktext{}%
  \gdef\gplfronttext{}%
  \makeatother
  \ifGPblacktext
    \def\colorrgb#1{}%
    \def\colorgray#1{}%
  \else
    \ifGPcolor
      \def\colorrgb#1{\color[rgb]{#1}}%
      \def\colorgray#1{\color[gray]{#1}}%
      \expandafter\def\csname LTw\endcsname{\color{white}}%
      \expandafter\def\csname LTb\endcsname{\color{black}}%
      \expandafter\def\csname LTa\endcsname{\color{black}}%
      \expandafter\def\csname LT0\endcsname{\color[rgb]{1,0,0}}%
      \expandafter\def\csname LT1\endcsname{\color[rgb]{0,1,0}}%
      \expandafter\def\csname LT2\endcsname{\color[rgb]{0,0,1}}%
      \expandafter\def\csname LT3\endcsname{\color[rgb]{1,0,1}}%
      \expandafter\def\csname LT4\endcsname{\color[rgb]{0,1,1}}%
      \expandafter\def\csname LT5\endcsname{\color[rgb]{1,1,0}}%
      \expandafter\def\csname LT6\endcsname{\color[rgb]{0,0,0}}%
      \expandafter\def\csname LT7\endcsname{\color[rgb]{1,0.3,0}}%
      \expandafter\def\csname LT8\endcsname{\color[rgb]{0.5,0.5,0.5}}%
    \else
      \def\colorrgb#1{\color{black}}%
      \def\colorgray#1{\color[gray]{#1}}%
      \expandafter\def\csname LTw\endcsname{\color{white}}%
      \expandafter\def\csname LTb\endcsname{\color{black}}%
      \expandafter\def\csname LTa\endcsname{\color{black}}%
      \expandafter\def\csname LT0\endcsname{\color{black}}%
      \expandafter\def\csname LT1\endcsname{\color{black}}%
      \expandafter\def\csname LT2\endcsname{\color{black}}%
      \expandafter\def\csname LT3\endcsname{\color{black}}%
      \expandafter\def\csname LT4\endcsname{\color{black}}%
      \expandafter\def\csname LT5\endcsname{\color{black}}%
      \expandafter\def\csname LT6\endcsname{\color{black}}%
      \expandafter\def\csname LT7\endcsname{\color{black}}%
      \expandafter\def\csname LT8\endcsname{\color{black}}%
    \fi
  \fi
    \setlength{\unitlength}{0.0500bp}%
    \ifx\gptboxheight\undefined%
      \newlength{\gptboxheight}%
      \newlength{\gptboxwidth}%
      \newsavebox{\gptboxtext}%
    \fi%
    \setlength{\fboxrule}{0.5pt}%
    \setlength{\fboxsep}{1pt}%
\begin{picture}(7200.00,5040.00)%
    \gplgaddtomacro\gplbacktext{%
      \csname LTb\endcsname
      \put(1876,1041){\makebox(0,0)[r]{\strut{}8.50e+05}}%
      \csname LTb\endcsname
      \put(1876,1685){\makebox(0,0)[r]{\strut{}9.00e+05}}%
      \csname LTb\endcsname
      \put(1876,2329){\makebox(0,0)[r]{\strut{}9.50e+05}}%
      \csname LTb\endcsname
      \put(1876,2973){\makebox(0,0)[r]{\strut{}1.00e+06}}%
      \csname LTb\endcsname
      \put(1876,3617){\makebox(0,0)[r]{\strut{}1.05e+06}}%
      \csname LTb\endcsname
      \put(1876,4261){\makebox(0,0)[r]{\strut{}1.10e+06}}%
      \csname LTb\endcsname
      \put(2044,616){\makebox(0,0){\strut{}0e+00}}%
      \csname LTb\endcsname
      \put(3373,616){\makebox(0,0){\strut{}1e-05}}%
      \csname LTb\endcsname
      \put(6031,616){\makebox(0,0){\strut{}3e-05}}%
    }%
    \gplgaddtomacro\gplfronttext{%
      \csname LTb\endcsname
      \put(-238,2827){\rotatebox{-270}{\makebox(0,0){\strut{}pressure}}}%
      \put(4369,196){\makebox(0,0){\strut{}time}}%
    }%
    \gplbacktext
    \put(0,0){\includegraphics[width={360.00bp},height={252.00bp}]{probe_dim6_ae2_uf3_iw15_ws7_res1024_ncycle0_AFDEIM_dt1e-09_Pressure_Probe1}}%
    \gplfronttext
  \end{picture}%
\endgroup

%% file: figures/probe_dim6_ae2_uf3_iw15_ws7_res1024_ncycle0_AFDEIM_dt1e-09_Pressure_Probe2.tex
\begingroup
  \makeatletter
  \providecommand\color[2][]{%
    \GenericError{(gnuplot) \space\space\space\@spaces}{%
      Package color not loaded in conjunction with
      terminal option `colourtext'%
    }{See the gnuplot documentation for explanation.%
    }{Either use 'blacktext' in gnuplot or load the package
      color.sty in LaTeX.}%
    \renewcommand\color[2][]{}%
  }%
  \providecommand\includegraphics[2][]{%
    \GenericError{(gnuplot) \space\space\space\@spaces}{%
      Package graphicx or graphics not loaded%
    }{See the gnuplot documentation for explanation.%
    }{The gnuplot epslatex terminal needs graphicx.sty or graphics.sty.}%
    \renewcommand\includegraphics[2][]{}%
  }%
  \providecommand\rotatebox[2]{#2}%
  \@ifundefined{ifGPcolor}{%
    \newif\ifGPcolor
    \GPcolortrue
  }{}%
  \@ifundefined{ifGPblacktext}{%
    \newif\ifGPblacktext
    \GPblacktexttrue
  }{}%
  \let\gplgaddtomacro\g@addto@macro
  \gdef\gplbacktext{}%
  \gdef\gplfronttext{}%
  \makeatother
  \ifGPblacktext
    \def\colorrgb#1{}%
    \def\colorgray#1{}%
  \else
    \ifGPcolor
      \def\colorrgb#1{\color[rgb]{#1}}%
      \def\colorgray#1{\color[gray]{#1}}%
      \expandafter\def\csname LTw\endcsname{\color{white}}%
      \expandafter\def\csname LTb\endcsname{\color{black}}%
      \expandafter\def\csname LTa\endcsname{\color{black}}%
      \expandafter\def\csname LT0\endcsname{\color[rgb]{1,0,0}}%
      \expandafter\def\csname LT1\endcsname{\color[rgb]{0,1,0}}%
      \expandafter\def\csname LT2\endcsname{\color[rgb]{0,0,1}}%
      \expandafter\def\csname LT3\endcsname{\color[rgb]{1,0,1}}%
      \expandafter\def\csname LT4\endcsname{\color[rgb]{0,1,1}}%
      \expandafter\def\csname LT5\endcsname{\color[rgb]{1,1,0}}%
      \expandafter\def\csname LT6\endcsname{\color[rgb]{0,0,0}}%
      \expandafter\def\csname LT7\endcsname{\color[rgb]{1,0.3,0}}%
      \expandafter\def\csname LT8\endcsname{\color[rgb]{0.5,0.5,0.5}}%
    \else
      \def\colorrgb#1{\color{black}}%
      \def\colorgray#1{\color[gray]{#1}}%
      \expandafter\def\csname LTw\endcsname{\color{white}}%
      \expandafter\def\csname LTb\endcsname{\color{black}}%
      \expandafter\def\csname LTa\endcsname{\color{black}}%
      \expandafter\def\csname LT0\endcsname{\color{black}}%
      \expandafter\def\csname LT1\endcsname{\color{black}}%
      \expandafter\def\csname LT2\endcsname{\color{black}}%
      \expandafter\def\csname LT3\endcsname{\color{black}}%
      \expandafter\def\csname LT4\endcsname{\color{black}}%
      \expandafter\def\csname LT5\endcsname{\color{black}}%
      \expandafter\def\csname LT6\endcsname{\color{black}}%
      \expandafter\def\csname LT7\endcsname{\color{black}}%
      \expandafter\def\csname LT8\endcsname{\color{black}}%
    \fi
  \fi
    \setlength{\unitlength}{0.0500bp}%
    \ifx\gptboxheight\undefined%
      \newlength{\gptboxheight}%
      \newlength{\gptboxwidth}%
      \newsavebox{\gptboxtext}%
    \fi%
    \setlength{\fboxrule}{0.5pt}%
    \setlength{\fboxsep}{1pt}%
\begin{picture}(7200.00,5040.00)%
    \gplgaddtomacro\gplbacktext{%
      \csname LTb\endcsname
      \put(1876,1458){\makebox(0,0)[r]{\strut{}9.00e+05}}%
      \csname LTb\endcsname
      \put(1876,2163){\makebox(0,0)[r]{\strut{}9.50e+05}}%
      \csname LTb\endcsname
      \put(1876,2868){\makebox(0,0)[r]{\strut{}1.00e+06}}%
      \csname LTb\endcsname
      \put(1876,3573){\makebox(0,0)[r]{\strut{}1.05e+06}}%
      \csname LTb\endcsname
      \put(1876,4278){\makebox(0,0)[r]{\strut{}1.10e+06}}%
      \csname LTb\endcsname
      \put(2044,616){\makebox(0,0){\strut{}0e+00}}%
      \csname LTb\endcsname
      \put(3373,616){\makebox(0,0){\strut{}1e-05}}%
      \csname LTb\endcsname
      \put(6031,616){\makebox(0,0){\strut{}3e-05}}%
    }%
    \gplgaddtomacro\gplfronttext{%
      \csname LTb\endcsname
      \put(-238,2827){\rotatebox{-270}{\makebox(0,0){\strut{}pressure}}}%
      \put(4369,196){\makebox(0,0){\strut{}time}}%
    }%
    \gplbacktext
    \put(0,0){\includegraphics[width={360.00bp},height={252.00bp}]{probe_dim6_ae2_uf3_iw15_ws7_res1024_ncycle0_AFDEIM_dt1e-09_Pressure_Probe2}}%
    \gplfronttext
  \end{picture}%
\endgroup

%% file: figures/probe_dim6_ae2_uf3_iw15_ws7_res1024_ncycle0_AFDEIM_dt1e-09_Velocity_Probe0.tex
\begingroup
  \makeatletter
  \providecommand\color[2][]{%
    \GenericError{(gnuplot) \space\space\space\@spaces}{%
      Package color not loaded in conjunction with
      terminal option `colourtext'%
    }{See the gnuplot documentation for explanation.%
    }{Either use 'blacktext' in gnuplot or load the package
      color.sty in LaTeX.}%
    \renewcommand\color[2][]{}%
  }%
  \providecommand\includegraphics[2][]{%
    \GenericError{(gnuplot) \space\space\space\@spaces}{%
      Package graphicx or graphics not loaded%
    }{See the gnuplot documentation for explanation.%
    }{The gnuplot epslatex terminal needs graphicx.sty or graphics.sty.}%
    \renewcommand\includegraphics[2][]{}%
  }%
  \providecommand\rotatebox[2]{#2}%
  \@ifundefined{ifGPcolor}{%
    \newif\ifGPcolor
    \GPcolortrue
  }{}%
  \@ifundefined{ifGPblacktext}{%
    \newif\ifGPblacktext
    \GPblacktexttrue
  }{}%
  \let\gplgaddtomacro\g@addto@macro
  \gdef\gplbacktext{}%
  \gdef\gplfronttext{}%
  \makeatother
  \ifGPblacktext
    \def\colorrgb#1{}%
    \def\colorgray#1{}%
  \else
    \ifGPcolor
      \def\colorrgb#1{\color[rgb]{#1}}%
      \def\colorgray#1{\color[gray]{#1}}%
      \expandafter\def\csname LTw\endcsname{\color{white}}%
      \expandafter\def\csname LTb\endcsname{\color{black}}%
      \expandafter\def\csname LTa\endcsname{\color{black}}%
      \expandafter\def\csname LT0\endcsname{\color[rgb]{1,0,0}}%
      \expandafter\def\csname LT1\endcsname{\color[rgb]{0,1,0}}%
      \expandafter\def\csname LT2\endcsname{\color[rgb]{0,0,1}}%
      \expandafter\def\csname LT3\endcsname{\color[rgb]{1,0,1}}%
      \expandafter\def\csname LT4\endcsname{\color[rgb]{0,1,1}}%
      \expandafter\def\csname LT5\endcsname{\color[rgb]{1,1,0}}%
      \expandafter\def\csname LT6\endcsname{\color[rgb]{0,0,0}}%
      \expandafter\def\csname LT7\endcsname{\color[rgb]{1,0.3,0}}%
      \expandafter\def\csname LT8\endcsname{\color[rgb]{0.5,0.5,0.5}}%
    \else
      \def\colorrgb#1{\color{black}}%
      \def\colorgray#1{\color[gray]{#1}}%
      \expandafter\def\csname LTw\endcsname{\color{white}}%
      \expandafter\def\csname LTb\endcsname{\color{black}}%
      \expandafter\def\csname LTa\endcsname{\color{black}}%
      \expandafter\def\csname LT0\endcsname{\color{black}}%
      \expandafter\def\csname LT1\endcsname{\color{black}}%
      \expandafter\def\csname LT2\endcsname{\color{black}}%
      \expandafter\def\csname LT3\endcsname{\color{black}}%
      \expandafter\def\csname LT4\endcsname{\color{black}}%
      \expandafter\def\csname LT5\endcsname{\color{black}}%
      \expandafter\def\csname LT6\endcsname{\color{black}}%
      \expandafter\def\csname LT7\endcsname{\color{black}}%
      \expandafter\def\csname LT8\endcsname{\color{black}}%
    \fi
  \fi
    \setlength{\unitlength}{0.0500bp}%
    \ifx\gptboxheight\undefined%
      \newlength{\gptboxheight}%
      \newlength{\gptboxwidth}%
      \newsavebox{\gptboxtext}%
    \fi%
    \setlength{\fboxrule}{0.5pt}%
    \setlength{\fboxsep}{1pt}%
\begin{picture}(7200.00,5040.00)%
    \gplgaddtomacro\gplbacktext{%
      \csname LTb\endcsname
      \put(868,1624){\makebox(0,0)[r]{\strut{}$0$}}%
      \csname LTb\endcsname
      \put(868,2361){\makebox(0,0)[r]{\strut{}$10$}}%
      \csname LTb\endcsname
      \put(868,3097){\makebox(0,0)[r]{\strut{}$20$}}%
      \csname LTb\endcsname
      \put(868,3834){\makebox(0,0)[r]{\strut{}$30$}}%
      \csname LTb\endcsname
      \put(868,4571){\makebox(0,0)[r]{\strut{}$40$}}%
      \csname LTb\endcsname
      \put(1036,616){\makebox(0,0){\strut{}0e+00}}%
      \csname LTb\endcsname
      \put(2653,616){\makebox(0,0){\strut{}1e-05}}%
      \csname LTb\endcsname
      \put(5887,616){\makebox(0,0){\strut{}3e-05}}%
    }%
    \gplgaddtomacro\gplfronttext{%
      \csname LTb\endcsname
      \put(98,2827){\rotatebox{-270}{\makebox(0,0){\strut{}velocity}}}%
      \put(3865,196){\makebox(0,0){\strut{}time}}%
      \csname LTb\endcsname
      \put(4475,4486){\makebox(0,0)[r]{\strut{} full model}}%
      \csname LTb\endcsname
      \put(4475,4066){\makebox(0,0)[r]{\strut{} AADEIM}}%
      \csname LTb\endcsname
      \put(4475,3646){\makebox(0,0)[r]{\strut{}static}}%
    }%
    \gplbacktext
    \put(0,0){\includegraphics[width={360.00bp},height={252.00bp}]{probe_dim6_ae2_uf3_iw15_ws7_res1024_ncycle0_AFDEIM_dt1e-09_Velocity_Probe0}}%
    \gplfronttext
  \end{picture}%
\endgroup

%% file: figures/probe_dim6_ae2_uf3_iw15_ws7_res1024_ncycle0_AFDEIM_dt1e-09_Velocity_Probe1.tex
\begingroup
  \makeatletter
  \providecommand\color[2][]{%
    \GenericError{(gnuplot) \space\space\space\@spaces}{%
      Package color not loaded in conjunction with
      terminal option `colourtext'%
    }{See the gnuplot documentation for explanation.%
    }{Either use 'blacktext' in gnuplot or load the package
      color.sty in LaTeX.}%
    \renewcommand\color[2][]{}%
  }%
  \providecommand\includegraphics[2][]{%
    \GenericError{(gnuplot) \space\space\space\@spaces}{%
      Package graphicx or graphics not loaded%
    }{See the gnuplot documentation for explanation.%
    }{The gnuplot epslatex terminal needs graphicx.sty or graphics.sty.}%
    \renewcommand\includegraphics[2][]{}%
  }%
  \providecommand\rotatebox[2]{#2}%
  \@ifundefined{ifGPcolor}{%
    \newif\ifGPcolor
    \GPcolortrue
  }{}%
  \@ifundefined{ifGPblacktext}{%
    \newif\ifGPblacktext
    \GPblacktexttrue
  }{}%
  \let\gplgaddtomacro\g@addto@macro
  \gdef\gplbacktext{}%
  \gdef\gplfronttext{}%
  \makeatother
  \ifGPblacktext
    \def\colorrgb#1{}%
    \def\colorgray#1{}%
  \else
    \ifGPcolor
      \def\colorrgb#1{\color[rgb]{#1}}%
      \def\colorgray#1{\color[gray]{#1}}%
      \expandafter\def\csname LTw\endcsname{\color{white}}%
      \expandafter\def\csname LTb\endcsname{\color{black}}%
      \expandafter\def\csname LTa\endcsname{\color{black}}%
      \expandafter\def\csname LT0\endcsname{\color[rgb]{1,0,0}}%
      \expandafter\def\csname LT1\endcsname{\color[rgb]{0,1,0}}%
      \expandafter\def\csname LT2\endcsname{\color[rgb]{0,0,1}}%
      \expandafter\def\csname LT3\endcsname{\color[rgb]{1,0,1}}%
      \expandafter\def\csname LT4\endcsname{\color[rgb]{0,1,1}}%
      \expandafter\def\csname LT5\endcsname{\color[rgb]{1,1,0}}%
      \expandafter\def\csname LT6\endcsname{\color[rgb]{0,0,0}}%
      \expandafter\def\csname LT7\endcsname{\color[rgb]{1,0.3,0}}%
      \expandafter\def\csname LT8\endcsname{\color[rgb]{0.5,0.5,0.5}}%
    \else
      \def\colorrgb#1{\color{black}}%
      \def\colorgray#1{\color[gray]{#1}}%
      \expandafter\def\csname LTw\endcsname{\color{white}}%
      \expandafter\def\csname LTb\endcsname{\color{black}}%
      \expandafter\def\csname LTa\endcsname{\color{black}}%
      \expandafter\def\csname LT0\endcsname{\color{black}}%
      \expandafter\def\csname LT1\endcsname{\color{black}}%
      \expandafter\def\csname LT2\endcsname{\color{black}}%
      \expandafter\def\csname LT3\endcsname{\color{black}}%
      \expandafter\def\csname LT4\endcsname{\color{black}}%
      \expandafter\def\csname LT5\endcsname{\color{black}}%
      \expandafter\def\csname LT6\endcsname{\color{black}}%
      \expandafter\def\csname LT7\endcsname{\color{black}}%
      \expandafter\def\csname LT8\endcsname{\color{black}}%
    \fi
  \fi
    \setlength{\unitlength}{0.0500bp}%
    \ifx\gptboxheight\undefined%
      \newlength{\gptboxheight}%
      \newlength{\gptboxwidth}%
      \newsavebox{\gptboxtext}%
    \fi%
    \setlength{\fboxrule}{0.5pt}%
    \setlength{\fboxsep}{1pt}%
\begin{picture}(7200.00,5040.00)%
    \gplgaddtomacro\gplbacktext{%
      \csname LTb\endcsname
      \put(1036,1503){\makebox(0,0)[r]{\strut{}$-20$}}%
      \csname LTb\endcsname
      \put(1036,2237){\makebox(0,0)[r]{\strut{}$0$}}%
      \csname LTb\endcsname
      \put(1036,2971){\makebox(0,0)[r]{\strut{}$20$}}%
      \csname LTb\endcsname
      \put(1036,3705){\makebox(0,0)[r]{\strut{}$40$}}%
      \csname LTb\endcsname
      \put(1036,4439){\makebox(0,0)[r]{\strut{}$60$}}%
      \csname LTb\endcsname
      \put(1204,616){\makebox(0,0){\strut{}0e+00}}%
      \csname LTb\endcsname
      \put(2773,616){\makebox(0,0){\strut{}1e-05}}%
      \csname LTb\endcsname
      \put(5911,616){\makebox(0,0){\strut{}3e-05}}%
    }%
    \gplgaddtomacro\gplfronttext{%
      \csname LTb\endcsname
      \put(98,2827){\rotatebox{-270}{\makebox(0,0){\strut{}velocity}}}%
      \put(3949,196){\makebox(0,0){\strut{}time}}%
    }%
    \gplbacktext
    \put(0,0){\includegraphics[width={360.00bp},height={252.00bp}]{probe_dim6_ae2_uf3_iw15_ws7_res1024_ncycle0_AFDEIM_dt1e-09_Velocity_Probe1}}%
    \gplfronttext
  \end{picture}%
\endgroup

%% file: figures/probe_dim6_ae2_uf3_iw15_ws7_res1024_ncycle0_AFDEIM_dt1e-09_Velocity_Probe2.tex
\begingroup
  \makeatletter
  \providecommand\color[2][]{%
    \GenericError{(gnuplot) \space\space\space\@spaces}{%
      Package color not loaded in conjunction with
      terminal option `colourtext'%
    }{See the gnuplot documentation for explanation.%
    }{Either use 'blacktext' in gnuplot or load the package
      color.sty in LaTeX.}%
    \renewcommand\color[2][]{}%
  }%
  \providecommand\includegraphics[2][]{%
    \GenericError{(gnuplot) \space\space\space\@spaces}{%
      Package graphicx or graphics not loaded%
    }{See the gnuplot documentation for explanation.%
    }{The gnuplot epslatex terminal needs graphicx.sty or graphics.sty.}%
    \renewcommand\includegraphics[2][]{}%
  }%
  \providecommand\rotatebox[2]{#2}%
  \@ifundefined{ifGPcolor}{%
    \newif\ifGPcolor
    \GPcolortrue
  }{}%
  \@ifundefined{ifGPblacktext}{%
    \newif\ifGPblacktext
    \GPblacktexttrue
  }{}%
  \let\gplgaddtomacro\g@addto@macro
  \gdef\gplbacktext{}%
  \gdef\gplfronttext{}%
  \makeatother
  \ifGPblacktext
    \def\colorrgb#1{}%
    \def\colorgray#1{}%
  \else
    \ifGPcolor
      \def\colorrgb#1{\color[rgb]{#1}}%
      \def\colorgray#1{\color[gray]{#1}}%
      \expandafter\def\csname LTw\endcsname{\color{white}}%
      \expandafter\def\csname LTb\endcsname{\color{black}}%
      \expandafter\def\csname LTa\endcsname{\color{black}}%
      \expandafter\def\csname LT0\endcsname{\color[rgb]{1,0,0}}%
      \expandafter\def\csname LT1\endcsname{\color[rgb]{0,1,0}}%
      \expandafter\def\csname LT2\endcsname{\color[rgb]{0,0,1}}%
      \expandafter\def\csname LT3\endcsname{\color[rgb]{1,0,1}}%
      \expandafter\def\csname LT4\endcsname{\color[rgb]{0,1,1}}%
      \expandafter\def\csname LT5\endcsname{\color[rgb]{1,1,0}}%
      \expandafter\def\csname LT6\endcsname{\color[rgb]{0,0,0}}%
      \expandafter\def\csname LT7\endcsname{\color[rgb]{1,0.3,0}}%
      \expandafter\def\csname LT8\endcsname{\color[rgb]{0.5,0.5,0.5}}%
    \else
      \def\colorrgb#1{\color{black}}%
      \def\colorgray#1{\color[gray]{#1}}%
      \expandafter\def\csname LTw\endcsname{\color{white}}%
      \expandafter\def\csname LTb\endcsname{\color{black}}%
      \expandafter\def\csname LTa\endcsname{\color{black}}%
      \expandafter\def\csname LT0\endcsname{\color{black}}%
      \expandafter\def\csname LT1\endcsname{\color{black}}%
      \expandafter\def\csname LT2\endcsname{\color{black}}%
      \expandafter\def\csname LT3\endcsname{\color{black}}%
      \expandafter\def\csname LT4\endcsname{\color{black}}%
      \expandafter\def\csname LT5\endcsname{\color{black}}%
      \expandafter\def\csname LT6\endcsname{\color{black}}%
      \expandafter\def\csname LT7\endcsname{\color{black}}%
      \expandafter\def\csname LT8\endcsname{\color{black}}%
    \fi
  \fi
    \setlength{\unitlength}{0.0500bp}%
    \ifx\gptboxheight\undefined%
      \newlength{\gptboxheight}%
      \newlength{\gptboxwidth}%
      \newsavebox{\gptboxtext}%
    \fi%
    \setlength{\fboxrule}{0.5pt}%
    \setlength{\fboxsep}{1pt}%
\begin{picture}(7200.00,5040.00)%
    \gplgaddtomacro\gplbacktext{%
      \csname LTb\endcsname
      \put(1036,1205){\makebox(0,0)[r]{\strut{}$-40$}}%
      \csname LTb\endcsname
      \put(1036,1777){\makebox(0,0)[r]{\strut{}$-20$}}%
      \csname LTb\endcsname
      \put(1036,2349){\makebox(0,0)[r]{\strut{}$0$}}%
      \csname LTb\endcsname
      \put(1036,2921){\makebox(0,0)[r]{\strut{}$20$}}%
      \csname LTb\endcsname
      \put(1036,3493){\makebox(0,0)[r]{\strut{}$40$}}%
      \csname LTb\endcsname
      \put(1036,4065){\makebox(0,0)[r]{\strut{}$60$}}%
      \csname LTb\endcsname
      \put(1036,4637){\makebox(0,0)[r]{\strut{}$80$}}%
      \csname LTb\endcsname
      \put(1204,616){\makebox(0,0){\strut{}0e+00}}%
      \csname LTb\endcsname
      \put(2773,616){\makebox(0,0){\strut{}1e-05}}%
      \csname LTb\endcsname
      \put(5911,616){\makebox(0,0){\strut{}3e-05}}%
    }%
    \gplgaddtomacro\gplfronttext{%
      \csname LTb\endcsname
      \put(98,2827){\rotatebox{-270}{\makebox(0,0){\strut{}velocity}}}%
      \put(3949,196){\makebox(0,0){\strut{}time}}%
    }%
    \gplbacktext
    \put(0,0){\includegraphics[width={360.00bp},height={252.00bp}]{probe_dim6_ae2_uf3_iw15_ws7_res1024_ncycle0_AFDEIM_dt1e-09_Velocity_Probe2}}%
    \gplfronttext
  \end{picture}%
\endgroup

%% file: figures/probe_dim6_ae2_uf3_iw15_ws7_res1024_ncycle0_AFDEIM_dt1e-09_Temperature_Probe0.tex
\begingroup
  \makeatletter
  \providecommand\color[2][]{%
    \GenericError{(gnuplot) \space\space\space\@spaces}{%
      Package color not loaded in conjunction with
      terminal option `colourtext'%
    }{See the gnuplot documentation for explanation.%
    }{Either use 'blacktext' in gnuplot or load the package
      color.sty in LaTeX.}%
    \renewcommand\color[2][]{}%
  }%
  \providecommand\includegraphics[2][]{%
    \GenericError{(gnuplot) \space\space\space\@spaces}{%
      Package graphicx or graphics not loaded%
    }{See the gnuplot documentation for explanation.%
    }{The gnuplot epslatex terminal needs graphicx.sty or graphics.sty.}%
    \renewcommand\includegraphics[2][]{}%
  }%
  \providecommand\rotatebox[2]{#2}%
  \@ifundefined{ifGPcolor}{%
    \newif\ifGPcolor
    \GPcolortrue
  }{}%
  \@ifundefined{ifGPblacktext}{%
    \newif\ifGPblacktext
    \GPblacktexttrue
  }{}%
  \let\gplgaddtomacro\g@addto@macro
  \gdef\gplbacktext{}%
  \gdef\gplfronttext{}%
  \makeatother
  \ifGPblacktext
    \def\colorrgb#1{}%
    \def\colorgray#1{}%
  \else
    \ifGPcolor
      \def\colorrgb#1{\color[rgb]{#1}}%
      \def\colorgray#1{\color[gray]{#1}}%
      \expandafter\def\csname LTw\endcsname{\color{white}}%
      \expandafter\def\csname LTb\endcsname{\color{black}}%
      \expandafter\def\csname LTa\endcsname{\color{black}}%
      \expandafter\def\csname LT0\endcsname{\color[rgb]{1,0,0}}%
      \expandafter\def\csname LT1\endcsname{\color[rgb]{0,1,0}}%
      \expandafter\def\csname LT2\endcsname{\color[rgb]{0,0,1}}%
      \expandafter\def\csname LT3\endcsname{\color[rgb]{1,0,1}}%
      \expandafter\def\csname LT4\endcsname{\color[rgb]{0,1,1}}%
      \expandafter\def\csname LT5\endcsname{\color[rgb]{1,1,0}}%
      \expandafter\def\csname LT6\endcsname{\color[rgb]{0,0,0}}%
      \expandafter\def\csname LT7\endcsname{\color[rgb]{1,0.3,0}}%
      \expandafter\def\csname LT8\endcsname{\color[rgb]{0.5,0.5,0.5}}%
    \else
      \def\colorrgb#1{\color{black}}%
      \def\colorgray#1{\color[gray]{#1}}%
      \expandafter\def\csname LTw\endcsname{\color{white}}%
      \expandafter\def\csname LTb\endcsname{\color{black}}%
      \expandafter\def\csname LTa\endcsname{\color{black}}%
      \expandafter\def\csname LT0\endcsname{\color{black}}%
      \expandafter\def\csname LT1\endcsname{\color{black}}%
      \expandafter\def\csname LT2\endcsname{\color{black}}%
      \expandafter\def\csname LT3\endcsname{\color{black}}%
      \expandafter\def\csname LT4\endcsname{\color{black}}%
      \expandafter\def\csname LT5\endcsname{\color{black}}%
      \expandafter\def\csname LT6\endcsname{\color{black}}%
      \expandafter\def\csname LT7\endcsname{\color{black}}%
      \expandafter\def\csname LT8\endcsname{\color{black}}%
    \fi
  \fi
    \setlength{\unitlength}{0.0500bp}%
    \ifx\gptboxheight\undefined%
      \newlength{\gptboxheight}%
      \newlength{\gptboxwidth}%
      \newsavebox{\gptboxtext}%
    \fi%
    \setlength{\fboxrule}{0.5pt}%
    \setlength{\fboxsep}{1pt}%
\begin{picture}(7200.00,5040.00)%
    \gplgaddtomacro\gplbacktext{%
      \csname LTb\endcsname
      \put(1204,1248){\makebox(0,0)[r]{\strut{}$1800$}}%
      \csname LTb\endcsname
      \put(1204,1960){\makebox(0,0)[r]{\strut{}$2000$}}%
      \csname LTb\endcsname
      \put(1204,2672){\makebox(0,0)[r]{\strut{}$2200$}}%
      \csname LTb\endcsname
      \put(1204,3383){\makebox(0,0)[r]{\strut{}$2400$}}%
      \csname LTb\endcsname
      \put(1204,4095){\makebox(0,0)[r]{\strut{}$2600$}}%
      \csname LTb\endcsname
      \put(1372,616){\makebox(0,0){\strut{}0e+00}}%
      \csname LTb\endcsname
      \put(2893,616){\makebox(0,0){\strut{}1e-05}}%
      \csname LTb\endcsname
      \put(5935,616){\makebox(0,0){\strut{}3e-05}}%
    }%
    \gplgaddtomacro\gplfronttext{%
      \csname LTb\endcsname
      \put(-70,2827){\rotatebox{-270}{\makebox(0,0){\strut{}temperature}}}%
      \put(4033,196){\makebox(0,0){\strut{}time}}%
      \csname LTb\endcsname
      \put(4811,2009){\makebox(0,0)[r]{\strut{} full model}}%
      \csname LTb\endcsname
      \put(4811,1589){\makebox(0,0)[r]{\strut{} AADEIM}}%
      \csname LTb\endcsname
      \put(4811,1169){\makebox(0,0)[r]{\strut{}static}}%
    }%
    \gplbacktext
    \put(0,0){\includegraphics[width={360.00bp},height={252.00bp}]{probe_dim6_ae2_uf3_iw15_ws7_res1024_ncycle0_AFDEIM_dt1e-09_Temperature_Probe0}}%
    \gplfronttext
  \end{picture}%
\endgroup

%% file: figures/probe_dim6_ae2_uf3_iw15_ws7_res1024_ncycle0_AFDEIM_dt1e-09_Temperature_Probe1.tex
\begingroup
  \makeatletter
  \providecommand\color[2][]{%
    \GenericError{(gnuplot) \space\space\space\@spaces}{%
      Package color not loaded in conjunction with
      terminal option `colourtext'%
    }{See the gnuplot documentation for explanation.%
    }{Either use 'blacktext' in gnuplot or load the package
      color.sty in LaTeX.}%
    \renewcommand\color[2][]{}%
  }%
  \providecommand\includegraphics[2][]{%
    \GenericError{(gnuplot) \space\space\space\@spaces}{%
      Package graphicx or graphics not loaded%
    }{See the gnuplot documentation for explanation.%
    }{The gnuplot epslatex terminal needs graphicx.sty or graphics.sty.}%
    \renewcommand\includegraphics[2][]{}%
  }%
  \providecommand\rotatebox[2]{#2}%
  \@ifundefined{ifGPcolor}{%
    \newif\ifGPcolor
    \GPcolortrue
  }{}%
  \@ifundefined{ifGPblacktext}{%
    \newif\ifGPblacktext
    \GPblacktexttrue
  }{}%
  \let\gplgaddtomacro\g@addto@macro
  \gdef\gplbacktext{}%
  \gdef\gplfronttext{}%
  \makeatother
  \ifGPblacktext
    \def\colorrgb#1{}%
    \def\colorgray#1{}%
  \else
    \ifGPcolor
      \def\colorrgb#1{\color[rgb]{#1}}%
      \def\colorgray#1{\color[gray]{#1}}%
      \expandafter\def\csname LTw\endcsname{\color{white}}%
      \expandafter\def\csname LTb\endcsname{\color{black}}%
      \expandafter\def\csname LTa\endcsname{\color{black}}%
      \expandafter\def\csname LT0\endcsname{\color[rgb]{1,0,0}}%
      \expandafter\def\csname LT1\endcsname{\color[rgb]{0,1,0}}%
      \expandafter\def\csname LT2\endcsname{\color[rgb]{0,0,1}}%
      \expandafter\def\csname LT3\endcsname{\color[rgb]{1,0,1}}%
      \expandafter\def\csname LT4\endcsname{\color[rgb]{0,1,1}}%
      \expandafter\def\csname LT5\endcsname{\color[rgb]{1,1,0}}%
      \expandafter\def\csname LT6\endcsname{\color[rgb]{0,0,0}}%
      \expandafter\def\csname LT7\endcsname{\color[rgb]{1,0.3,0}}%
      \expandafter\def\csname LT8\endcsname{\color[rgb]{0.5,0.5,0.5}}%
    \else
      \def\colorrgb#1{\color{black}}%
      \def\colorgray#1{\color[gray]{#1}}%
      \expandafter\def\csname LTw\endcsname{\color{white}}%
      \expandafter\def\csname LTb\endcsname{\color{black}}%
      \expandafter\def\csname LTa\endcsname{\color{black}}%
      \expandafter\def\csname LT0\endcsname{\color{black}}%
      \expandafter\def\csname LT1\endcsname{\color{black}}%
      \expandafter\def\csname LT2\endcsname{\color{black}}%
      \expandafter\def\csname LT3\endcsname{\color{black}}%
      \expandafter\def\csname LT4\endcsname{\color{black}}%
      \expandafter\def\csname LT5\endcsname{\color{black}}%
      \expandafter\def\csname LT6\endcsname{\color{black}}%
      \expandafter\def\csname LT7\endcsname{\color{black}}%
      \expandafter\def\csname LT8\endcsname{\color{black}}%
    \fi
  \fi
    \setlength{\unitlength}{0.0500bp}%
    \ifx\gptboxheight\undefined%
      \newlength{\gptboxheight}%
      \newlength{\gptboxwidth}%
      \newsavebox{\gptboxtext}%
    \fi%
    \setlength{\fboxrule}{0.5pt}%
    \setlength{\fboxsep}{1pt}%
\begin{picture}(7200.00,5040.00)%
    \gplgaddtomacro\gplbacktext{%
      \csname LTb\endcsname
      \put(1204,1503){\makebox(0,0)[r]{\strut{}$2300$}}%
      \csname LTb\endcsname
      \put(1204,2187){\makebox(0,0)[r]{\strut{}$2400$}}%
      \csname LTb\endcsname
      \put(1204,2872){\makebox(0,0)[r]{\strut{}$2500$}}%
      \csname LTb\endcsname
      \put(1204,3556){\makebox(0,0)[r]{\strut{}$2600$}}%
      \csname LTb\endcsname
      \put(1204,4241){\makebox(0,0)[r]{\strut{}$2700$}}%
      \csname LTb\endcsname
      \put(1372,616){\makebox(0,0){\strut{}0e+00}}%
      \csname LTb\endcsname
      \put(2893,616){\makebox(0,0){\strut{}1e-05}}%
      \csname LTb\endcsname
      \put(5935,616){\makebox(0,0){\strut{}3e-05}}%
    }%
    \gplgaddtomacro\gplfronttext{%
      \csname LTb\endcsname
      \put(-70,2827){\rotatebox{-270}{\makebox(0,0){\strut{}temperature}}}%
      \put(4033,196){\makebox(0,0){\strut{}time}}%
    }%
    \gplbacktext
    \put(0,0){\includegraphics[width={360.00bp},height={252.00bp}]{probe_dim6_ae2_uf3_iw15_ws7_res1024_ncycle0_AFDEIM_dt1e-09_Temperature_Probe1}}%
    \gplfronttext
  \end{picture}%
\endgroup

%% file: figures/probe_dim6_ae2_uf3_iw15_ws7_res1024_ncycle0_AFDEIM_dt1e-09_Temperature_Probe2.tex
\begingroup
  \makeatletter
  \providecommand\color[2][]{%
    \GenericError{(gnuplot) \space\space\space\@spaces}{%
      Package color not loaded in conjunction with
      terminal option `colourtext'%
    }{See the gnuplot documentation for explanation.%
    }{Either use 'blacktext' in gnuplot or load the package
      color.sty in LaTeX.}%
    \renewcommand\color[2][]{}%
  }%
  \providecommand\includegraphics[2][]{%
    \GenericError{(gnuplot) \space\space\space\@spaces}{%
      Package graphicx or graphics not loaded%
    }{See the gnuplot documentation for explanation.%
    }{The gnuplot epslatex terminal needs graphicx.sty or graphics.sty.}%
    \renewcommand\includegraphics[2][]{}%
  }%
  \providecommand\rotatebox[2]{#2}%
  \@ifundefined{ifGPcolor}{%
    \newif\ifGPcolor
    \GPcolortrue
  }{}%
  \@ifundefined{ifGPblacktext}{%
    \newif\ifGPblacktext
    \GPblacktexttrue
  }{}%
  \let\gplgaddtomacro\g@addto@macro
  \gdef\gplbacktext{}%
  \gdef\gplfronttext{}%
  \makeatother
  \ifGPblacktext
    \def\colorrgb#1{}%
    \def\colorgray#1{}%
  \else
    \ifGPcolor
      \def\colorrgb#1{\color[rgb]{#1}}%
      \def\colorgray#1{\color[gray]{#1}}%
      \expandafter\def\csname LTw\endcsname{\color{white}}%
      \expandafter\def\csname LTb\endcsname{\color{black}}%
      \expandafter\def\csname LTa\endcsname{\color{black}}%
      \expandafter\def\csname LT0\endcsname{\color[rgb]{1,0,0}}%
      \expandafter\def\csname LT1\endcsname{\color[rgb]{0,1,0}}%
      \expandafter\def\csname LT2\endcsname{\color[rgb]{0,0,1}}%
      \expandafter\def\csname LT3\endcsname{\color[rgb]{1,0,1}}%
      \expandafter\def\csname LT4\endcsname{\color[rgb]{0,1,1}}%
      \expandafter\def\csname LT5\endcsname{\color[rgb]{1,1,0}}%
      \expandafter\def\csname LT6\endcsname{\color[rgb]{0,0,0}}%
      \expandafter\def\csname LT7\endcsname{\color[rgb]{1,0.3,0}}%
      \expandafter\def\csname LT8\endcsname{\color[rgb]{0.5,0.5,0.5}}%
    \else
      \def\colorrgb#1{\color{black}}%
      \def\colorgray#1{\color[gray]{#1}}%
      \expandafter\def\csname LTw\endcsname{\color{white}}%
      \expandafter\def\csname LTb\endcsname{\color{black}}%
      \expandafter\def\csname LTa\endcsname{\color{black}}%
      \expandafter\def\csname LT0\endcsname{\color{black}}%
      \expandafter\def\csname LT1\endcsname{\color{black}}%
      \expandafter\def\csname LT2\endcsname{\color{black}}%
      \expandafter\def\csname LT3\endcsname{\color{black}}%
      \expandafter\def\csname LT4\endcsname{\color{black}}%
      \expandafter\def\csname LT5\endcsname{\color{black}}%
      \expandafter\def\csname LT6\endcsname{\color{black}}%
      \expandafter\def\csname LT7\endcsname{\color{black}}%
      \expandafter\def\csname LT8\endcsname{\color{black}}%
    \fi
  \fi
    \setlength{\unitlength}{0.0500bp}%
    \ifx\gptboxheight\undefined%
      \newlength{\gptboxheight}%
      \newlength{\gptboxwidth}%
      \newsavebox{\gptboxtext}%
    \fi%
    \setlength{\fboxrule}{0.5pt}%
    \setlength{\fboxsep}{1pt}%
\begin{picture}(7200.00,5040.00)%
    \gplgaddtomacro\gplbacktext{%
      \csname LTb\endcsname
      \put(1204,1422){\makebox(0,0)[r]{\strut{}$2300$}}%
      \csname LTb\endcsname
      \put(1204,2127){\makebox(0,0)[r]{\strut{}$2400$}}%
      \csname LTb\endcsname
      \put(1204,2833){\makebox(0,0)[r]{\strut{}$2500$}}%
      \csname LTb\endcsname
      \put(1204,3539){\makebox(0,0)[r]{\strut{}$2600$}}%
      \csname LTb\endcsname
      \put(1204,4245){\makebox(0,0)[r]{\strut{}$2700$}}%
      \csname LTb\endcsname
      \put(1372,616){\makebox(0,0){\strut{}0e+00}}%
      \csname LTb\endcsname
      \put(2893,616){\makebox(0,0){\strut{}1e-05}}%
      \csname LTb\endcsname
      \put(5935,616){\makebox(0,0){\strut{}3e-05}}%
    }%
    \gplgaddtomacro\gplfronttext{%
      \csname LTb\endcsname
      \put(-70,2827){\rotatebox{-270}{\makebox(0,0){\strut{}temperature}}}%
      \put(4033,196){\makebox(0,0){\strut{}time}}%
    }%
    \gplbacktext
    \put(0,0){\includegraphics[width={360.00bp},height={252.00bp}]{probe_dim6_ae2_uf3_iw15_ws7_res1024_ncycle0_AFDEIM_dt1e-09_Temperature_Probe2}}%
    \gplfronttext
  \end{picture}%
\endgroup

%% file: figures/probe_dim6_ae2_uf3_iw15_ws7_res1024_ncycle0_AFDEIM_dt1e-09_MassFraction_Probe0.tex
\begingroup
  \makeatletter
  \providecommand\color[2][]{%
    \GenericError{(gnuplot) \space\space\space\@spaces}{%
      Package color not loaded in conjunction with
      terminal option `colourtext'%
    }{See the gnuplot documentation for explanation.%
    }{Either use 'blacktext' in gnuplot or load the package
      color.sty in LaTeX.}%
    \renewcommand\color[2][]{}%
  }%
  \providecommand\includegraphics[2][]{%
    \GenericError{(gnuplot) \space\space\space\@spaces}{%
      Package graphicx or graphics not loaded%
    }{See the gnuplot documentation for explanation.%
    }{The gnuplot epslatex terminal needs graphicx.sty or graphics.sty.}%
    \renewcommand\includegraphics[2][]{}%
  }%
  \providecommand\rotatebox[2]{#2}%
  \@ifundefined{ifGPcolor}{%
    \newif\ifGPcolor
    \GPcolortrue
  }{}%
  \@ifundefined{ifGPblacktext}{%
    \newif\ifGPblacktext
    \GPblacktexttrue
  }{}%
  \let\gplgaddtomacro\g@addto@macro
  \gdef\gplbacktext{}%
  \gdef\gplfronttext{}%
  \makeatother
  \ifGPblacktext
    \def\colorrgb#1{}%
    \def\colorgray#1{}%
  \else
    \ifGPcolor
      \def\colorrgb#1{\color[rgb]{#1}}%
      \def\colorgray#1{\color[gray]{#1}}%
      \expandafter\def\csname LTw\endcsname{\color{white}}%
      \expandafter\def\csname LTb\endcsname{\color{black}}%
      \expandafter\def\csname LTa\endcsname{\color{black}}%
      \expandafter\def\csname LT0\endcsname{\color[rgb]{1,0,0}}%
      \expandafter\def\csname LT1\endcsname{\color[rgb]{0,1,0}}%
      \expandafter\def\csname LT2\endcsname{\color[rgb]{0,0,1}}%
      \expandafter\def\csname LT3\endcsname{\color[rgb]{1,0,1}}%
      \expandafter\def\csname LT4\endcsname{\color[rgb]{0,1,1}}%
      \expandafter\def\csname LT5\endcsname{\color[rgb]{1,1,0}}%
      \expandafter\def\csname LT6\endcsname{\color[rgb]{0,0,0}}%
      \expandafter\def\csname LT7\endcsname{\color[rgb]{1,0.3,0}}%
      \expandafter\def\csname LT8\endcsname{\color[rgb]{0.5,0.5,0.5}}%
    \else
      \def\colorrgb#1{\color{black}}%
      \def\colorgray#1{\color[gray]{#1}}%
      \expandafter\def\csname LTw\endcsname{\color{white}}%
      \expandafter\def\csname LTb\endcsname{\color{black}}%
      \expandafter\def\csname LTa\endcsname{\color{black}}%
      \expandafter\def\csname LT0\endcsname{\color{black}}%
      \expandafter\def\csname LT1\endcsname{\color{black}}%
      \expandafter\def\csname LT2\endcsname{\color{black}}%
      \expandafter\def\csname LT3\endcsname{\color{black}}%
      \expandafter\def\csname LT4\endcsname{\color{black}}%
      \expandafter\def\csname LT5\endcsname{\color{black}}%
      \expandafter\def\csname LT6\endcsname{\color{black}}%
      \expandafter\def\csname LT7\endcsname{\color{black}}%
      \expandafter\def\csname LT8\endcsname{\color{black}}%
    \fi
  \fi
    \setlength{\unitlength}{0.0500bp}%
    \ifx\gptboxheight\undefined%
      \newlength{\gptboxheight}%
      \newlength{\gptboxwidth}%
      \newsavebox{\gptboxtext}%
    \fi%
    \setlength{\fboxrule}{0.5pt}%
    \setlength{\fboxsep}{1pt}%
\begin{picture}(7200.00,5040.00)%
    \gplgaddtomacro\gplbacktext{%
      \csname LTb\endcsname
      \put(1204,894){\makebox(0,0)[r]{\strut{}$0$}}%
      \csname LTb\endcsname
      \put(1204,1601){\makebox(0,0)[r]{\strut{}$0.05$}}%
      \csname LTb\endcsname
      \put(1204,2308){\makebox(0,0)[r]{\strut{}$0.1$}}%
      \csname LTb\endcsname
      \put(1204,3016){\makebox(0,0)[r]{\strut{}$0.15$}}%
      \csname LTb\endcsname
      \put(1204,3723){\makebox(0,0)[r]{\strut{}$0.2$}}%
      \csname LTb\endcsname
      \put(1204,4431){\makebox(0,0)[r]{\strut{}$0.25$}}%
      \csname LTb\endcsname
      \put(1372,616){\makebox(0,0){\strut{}0e+00}}%
      \csname LTb\endcsname
      \put(2893,616){\makebox(0,0){\strut{}1e-05}}%
      \csname LTb\endcsname
      \put(5935,616){\makebox(0,0){\strut{}3e-05}}%
    }%
    \gplgaddtomacro\gplfronttext{%
      \csname LTb\endcsname
      \put(98,2827){\rotatebox{-270}{\makebox(0,0){\strut{}species mass fraction}}}%
      \put(4033,196){\makebox(0,0){\strut{}time}}%
      \csname LTb\endcsname
      \put(4811,4486){\makebox(0,0)[r]{\strut{} full model}}%
      \csname LTb\endcsname
      \put(4811,4066){\makebox(0,0)[r]{\strut{} AADEIM}}%
      \csname LTb\endcsname
      \put(4811,3646){\makebox(0,0)[r]{\strut{}static}}%
    }%
    \gplbacktext
    \put(0,0){\includegraphics[width={360.00bp},height={252.00bp}]{probe_dim6_ae2_uf3_iw15_ws7_res1024_ncycle0_AFDEIM_dt1e-09_MassFraction_Probe0}}%
    \gplfronttext
  \end{picture}%
\endgroup

%% file: figures/probe_dim6_ae2_uf3_iw15_ws7_res1024_ncycle0_AFDEIM_dt1e-09_MassFraction_Probe1.tex
\begingroup
  \makeatletter
  \providecommand\color[2][]{%
    \GenericError{(gnuplot) \space\space\space\@spaces}{%
      Package color not loaded in conjunction with
      terminal option `colourtext'%
    }{See the gnuplot documentation for explanation.%
    }{Either use 'blacktext' in gnuplot or load the package
      color.sty in LaTeX.}%
    \renewcommand\color[2][]{}%
  }%
  \providecommand\includegraphics[2][]{%
    \GenericError{(gnuplot) \space\space\space\@spaces}{%
      Package graphicx or graphics not loaded%
    }{See the gnuplot documentation for explanation.%
    }{The gnuplot epslatex terminal needs graphicx.sty or graphics.sty.}%
    \renewcommand\includegraphics[2][]{}%
  }%
  \providecommand\rotatebox[2]{#2}%
  \@ifundefined{ifGPcolor}{%
    \newif\ifGPcolor
    \GPcolortrue
  }{}%
  \@ifundefined{ifGPblacktext}{%
    \newif\ifGPblacktext
    \GPblacktexttrue
  }{}%
  \let\gplgaddtomacro\g@addto@macro
  \gdef\gplbacktext{}%
  \gdef\gplfronttext{}%
  \makeatother
  \ifGPblacktext
    \def\colorrgb#1{}%
    \def\colorgray#1{}%
  \else
    \ifGPcolor
      \def\colorrgb#1{\color[rgb]{#1}}%
      \def\colorgray#1{\color[gray]{#1}}%
      \expandafter\def\csname LTw\endcsname{\color{white}}%
      \expandafter\def\csname LTb\endcsname{\color{black}}%
      \expandafter\def\csname LTa\endcsname{\color{black}}%
      \expandafter\def\csname LT0\endcsname{\color[rgb]{1,0,0}}%
      \expandafter\def\csname LT1\endcsname{\color[rgb]{0,1,0}}%
      \expandafter\def\csname LT2\endcsname{\color[rgb]{0,0,1}}%
      \expandafter\def\csname LT3\endcsname{\color[rgb]{1,0,1}}%
      \expandafter\def\csname LT4\endcsname{\color[rgb]{0,1,1}}%
      \expandafter\def\csname LT5\endcsname{\color[rgb]{1,1,0}}%
      \expandafter\def\csname LT6\endcsname{\color[rgb]{0,0,0}}%
      \expandafter\def\csname LT7\endcsname{\color[rgb]{1,0.3,0}}%
      \expandafter\def\csname LT8\endcsname{\color[rgb]{0.5,0.5,0.5}}%
    \else
      \def\colorrgb#1{\color{black}}%
      \def\colorgray#1{\color[gray]{#1}}%
      \expandafter\def\csname LTw\endcsname{\color{white}}%
      \expandafter\def\csname LTb\endcsname{\color{black}}%
      \expandafter\def\csname LTa\endcsname{\color{black}}%
      \expandafter\def\csname LT0\endcsname{\color{black}}%
      \expandafter\def\csname LT1\endcsname{\color{black}}%
      \expandafter\def\csname LT2\endcsname{\color{black}}%
      \expandafter\def\csname LT3\endcsname{\color{black}}%
      \expandafter\def\csname LT4\endcsname{\color{black}}%
      \expandafter\def\csname LT5\endcsname{\color{black}}%
      \expandafter\def\csname LT6\endcsname{\color{black}}%
      \expandafter\def\csname LT7\endcsname{\color{black}}%
      \expandafter\def\csname LT8\endcsname{\color{black}}%
    \fi
  \fi
    \setlength{\unitlength}{0.0500bp}%
    \ifx\gptboxheight\undefined%
      \newlength{\gptboxheight}%
      \newlength{\gptboxwidth}%
      \newsavebox{\gptboxtext}%
    \fi%
    \setlength{\fboxrule}{0.5pt}%
    \setlength{\fboxsep}{1pt}%
\begin{picture}(7200.00,5040.00)%
    \gplgaddtomacro\gplbacktext{%
      \csname LTb\endcsname
      \put(1036,896){\makebox(0,0)[r]{\strut{}$0$}}%
      \csname LTb\endcsname
      \put(1036,1669){\makebox(0,0)[r]{\strut{}$0.2$}}%
      \csname LTb\endcsname
      \put(1036,2441){\makebox(0,0)[r]{\strut{}$0.4$}}%
      \csname LTb\endcsname
      \put(1036,3214){\makebox(0,0)[r]{\strut{}$0.6$}}%
      \csname LTb\endcsname
      \put(1036,3986){\makebox(0,0)[r]{\strut{}$0.8$}}%
      \csname LTb\endcsname
      \put(1036,4759){\makebox(0,0)[r]{\strut{}$1$}}%
      \csname LTb\endcsname
      \put(1204,616){\makebox(0,0){\strut{}0e+00}}%
      \csname LTb\endcsname
      \put(2773,616){\makebox(0,0){\strut{}1e-05}}%
      \csname LTb\endcsname
      \put(5911,616){\makebox(0,0){\strut{}3e-05}}%
    }%
    \gplgaddtomacro\gplfronttext{%
      \csname LTb\endcsname
      \put(98,2827){\rotatebox{-270}{\makebox(0,0){\strut{}species mass fraction}}}%
      \put(3949,196){\makebox(0,0){\strut{}time}}%
    }%
    \gplbacktext
    \put(0,0){\includegraphics[width={360.00bp},height={252.00bp}]{probe_dim6_ae2_uf3_iw15_ws7_res1024_ncycle0_AFDEIM_dt1e-09_MassFraction_Probe1}}%
    \gplfronttext
  \end{picture}%
\endgroup

%% file: figures/probe_dim6_ae2_uf3_iw15_ws7_res1024_ncycle0_AFDEIM_dt1e-09_MassFraction_Probe2.tex
\begingroup
  \makeatletter
  \providecommand\color[2][]{%
    \GenericError{(gnuplot) \space\space\space\@spaces}{%
      Package color not loaded in conjunction with
      terminal option `colourtext'%
    }{See the gnuplot documentation for explanation.%
    }{Either use 'blacktext' in gnuplot or load the package
      color.sty in LaTeX.}%
    \renewcommand\color[2][]{}%
  }%
  \providecommand\includegraphics[2][]{%
    \GenericError{(gnuplot) \space\space\space\@spaces}{%
      Package graphicx or graphics not loaded%
    }{See the gnuplot documentation for explanation.%
    }{The gnuplot epslatex terminal needs graphicx.sty or graphics.sty.}%
    \renewcommand\includegraphics[2][]{}%
  }%
  \providecommand\rotatebox[2]{#2}%
  \@ifundefined{ifGPcolor}{%
    \newif\ifGPcolor
    \GPcolortrue
  }{}%
  \@ifundefined{ifGPblacktext}{%
    \newif\ifGPblacktext
    \GPblacktexttrue
  }{}%
  \let\gplgaddtomacro\g@addto@macro
  \gdef\gplbacktext{}%
  \gdef\gplfronttext{}%
  \makeatother
  \ifGPblacktext
    \def\colorrgb#1{}%
    \def\colorgray#1{}%
  \else
    \ifGPcolor
      \def\colorrgb#1{\color[rgb]{#1}}%
      \def\colorgray#1{\color[gray]{#1}}%
      \expandafter\def\csname LTw\endcsname{\color{white}}%
      \expandafter\def\csname LTb\endcsname{\color{black}}%
      \expandafter\def\csname LTa\endcsname{\color{black}}%
      \expandafter\def\csname LT0\endcsname{\color[rgb]{1,0,0}}%
      \expandafter\def\csname LT1\endcsname{\color[rgb]{0,1,0}}%
      \expandafter\def\csname LT2\endcsname{\color[rgb]{0,0,1}}%
      \expandafter\def\csname LT3\endcsname{\color[rgb]{1,0,1}}%
      \expandafter\def\csname LT4\endcsname{\color[rgb]{0,1,1}}%
      \expandafter\def\csname LT5\endcsname{\color[rgb]{1,1,0}}%
      \expandafter\def\csname LT6\endcsname{\color[rgb]{0,0,0}}%
      \expandafter\def\csname LT7\endcsname{\color[rgb]{1,0.3,0}}%
      \expandafter\def\csname LT8\endcsname{\color[rgb]{0.5,0.5,0.5}}%
    \else
      \def\colorrgb#1{\color{black}}%
      \def\colorgray#1{\color[gray]{#1}}%
      \expandafter\def\csname LTw\endcsname{\color{white}}%
      \expandafter\def\csname LTb\endcsname{\color{black}}%
      \expandafter\def\csname LTa\endcsname{\color{black}}%
      \expandafter\def\csname LT0\endcsname{\color{black}}%
      \expandafter\def\csname LT1\endcsname{\color{black}}%
      \expandafter\def\csname LT2\endcsname{\color{black}}%
      \expandafter\def\csname LT3\endcsname{\color{black}}%
      \expandafter\def\csname LT4\endcsname{\color{black}}%
      \expandafter\def\csname LT5\endcsname{\color{black}}%
      \expandafter\def\csname LT6\endcsname{\color{black}}%
      \expandafter\def\csname LT7\endcsname{\color{black}}%
      \expandafter\def\csname LT8\endcsname{\color{black}}%
    \fi
  \fi
    \setlength{\unitlength}{0.0500bp}%
    \ifx\gptboxheight\undefined%
      \newlength{\gptboxheight}%
      \newlength{\gptboxwidth}%
      \newsavebox{\gptboxtext}%
    \fi%
    \setlength{\fboxrule}{0.5pt}%
    \setlength{\fboxsep}{1pt}%
\begin{picture}(7200.00,5040.00)%
    \gplgaddtomacro\gplbacktext{%
      \csname LTb\endcsname
      \put(1036,896){\makebox(0,0)[r]{\strut{}$0$}}%
      \csname LTb\endcsname
      \put(1036,1669){\makebox(0,0)[r]{\strut{}$0.2$}}%
      \csname LTb\endcsname
      \put(1036,2441){\makebox(0,0)[r]{\strut{}$0.4$}}%
      \csname LTb\endcsname
      \put(1036,3214){\makebox(0,0)[r]{\strut{}$0.6$}}%
      \csname LTb\endcsname
      \put(1036,3986){\makebox(0,0)[r]{\strut{}$0.8$}}%
      \csname LTb\endcsname
      \put(1036,4759){\makebox(0,0)[r]{\strut{}$1$}}%
      \csname LTb\endcsname
      \put(1204,616){\makebox(0,0){\strut{}0e+00}}%
      \csname LTb\endcsname
      \put(2773,616){\makebox(0,0){\strut{}1e-05}}%
      \csname LTb\endcsname
      \put(5911,616){\makebox(0,0){\strut{}3e-05}}%
    }%
    \gplgaddtomacro\gplfronttext{%
      \csname LTb\endcsname
      \put(98,2827){\rotatebox{-270}{\makebox(0,0){\strut{}species mass fraction}}}%
      \put(3949,196){\makebox(0,0){\strut{}time}}%
    }%
    \gplbacktext
    \put(0,0){\includegraphics[width={360.00bp},height={252.00bp}]{probe_dim6_ae2_uf3_iw15_ws7_res1024_ncycle0_AFDEIM_dt1e-09_MassFraction_Probe2}}%
    \gplfronttext
  \end{picture}%
\endgroup